\documentclass[12pt]{amsart}
\usepackage[table]{xcolor}
\usepackage[margin=1.25in]{geometry}

\usepackage[utf8]{inputenc}
\usepackage{tikz}
\usetikzlibrary{patterns}
\usepackage{amssymb,latexsym,amsmath,extarrows}
\usepackage[pagewise]{lineno}

\usepackage{listings}
\usepackage{makecell}

\usepackage{graphicx,mathrsfs,comment}
\usepackage{hyperref,url}

\newcommand{\x}{{\bf{x}}}
\def\set4{\mathcal I}
\def\tup14{(1,2,3,4)}
\def\bv{{\mathbf v}}

\makeatletter
\newcommand\@avprod[2]{%
  {\sbox0{$\m@th#1\prod$}%
   \vphantom{\usebox0}%
   \ooalign{%
     \hidewidth
     \smash{\vrule height\dimexpr\ht0+1pt\relax depth\dimexpr\dp0+1pt\relax}%
     \hidewidth\cr
     $\m@th#1\prod$\cr
   }%
  }%
}
\newcommand{\avprod}{\mathop{\mathpalette\@avprod\relax}\displaylimits}

\usepackage{amstext}
\usepackage{bbm}

\numberwithin{equation}{section}

\newtheorem{theorem}{Theorem}
\newtheorem{lemma}{Lemma}[section]

\newtheorem{proposition}[lemma]{Proposition}
\newtheorem{remark}[lemma]{Remark}
\newtheorem{example}[lemma]{Example}

\newtheorem{definition}[lemma]{Definition}
\newtheorem{corollary}[lemma]{Corollary}

\newtheorem{claim}[lemma]{Claim}

\newtheorem{condition}{Condition}

\newtheorem{observation}{Observation}

\usepackage{esint}

\usepackage{enumitem}

\newcommand{\eps}{\epsilon}
\newcommand{\bn}{\mathbf{n}}
\renewcommand{\bv}{\mathbf{v}}

\newcommand{\cC}{\mathcal C}

\renewcommand{\bn}{{\bf n}}

\renewcommand{\cC}{{\mathcal C}}

\newcommand{\R}{\mathbb{R}}

\newcommand{\p}{\mathbf p}
\renewcommand{\x}{\mathbf x}
\renewcommand{\max}{\rm{max}}
\renewcommand{\set}{\rm{set}}

\newcommand{\xixi}{\xi_0 + \xi}
\newcommand{\sym}{\mathrm{Sym}}
\newcommand{\Bphi}{B_{\epsilon_\phi}}

\begin{document}

\title[]{New Curved Kakeya Estimates}

\author{Arian Nadjimzadah} \address{Arian Nadjimzadah\\  Deparment of Mathematics, University of California, Los Angeles, USA}\email{anad@math.ucla.edu}
\thanks{The author was supported in part by NSF grant DMS-2347850. }

\maketitle

\begin{abstract}

We prove that for an open class of translation-invariant H\"ormander-type phase functions, the corresponding curved Kakeya sets in $\R^3$ have Hausdorff dimension at least $(13-\sqrt{13})/4 = 2.348\ldots$.
Previous methods, such as polynomial partitioning, have attained at best dimension $2+1/3$ for curved Kakeya sets in $\R^3$. The proof combines Wolff's classical hairbrush argument with a new Kakeya estimate for 3-parameter families of curves satisfying stable geometric conditions that we call coniness and twistiness. The latter estimate extends ideas of Katz, Wu, and Zahl from the study of $\mathrm{SL}_2$-Kakeya sets.  
\end{abstract}

\section{Introduction}\label{sec: 1}

\subsection{Curved Kakeya in $\R^3$ beyond the $2+1/3$ threshold}

As first discovered by Bourgain \cite{BourgainSeveralVariables}, curved Kakeya problems arise naturally from the wavepacket geometry of variable-coefficient oscillatory integral operators. This connection was further developed by Wisewell in \cite{Wisewell}, and we refer the reader there for details. 
Let $\phi : \R^3 \times \R^2 \to \R$ be a smooth phase function satisfying H\"ormander's non-degeneracy conditions \cite{HormanderOriginal}, and write $\x = (x,t) \in \R^2 \times \R$. For $\eps_\phi > 0$ sufficiently small and $(\xi,v) \in B_{\eps_\phi}^2 \times B_{\eps_{\phi}}^2$, the associated $\phi$-curve is 
\begin{align}
    \ell_{\xi,v} := \{\x \in B_{\eps_\phi}^3 : \nabla_\xi \phi(\x, \xi) = v\}.
\end{align}
We write $\mathcal C(\phi)$ for this 4-parameter family and call $\xi$ the direction parameter. In the \emph{translation-invariant} case 
\begin{align}
    \phi(x,t,\xi) = \langle x,\xi \rangle + \psi(t,\xi),
\end{align}
H\"ormander's non-degeneracy conditions are equivalent to $\partial_t \nabla_\xi^2 \psi(t,\xi) \neq 0$, and
the curves have the explicit parameterization 
\begin{align}
    \ell_{\xi,v}(t) = (v - \nabla_\xi \psi(t,\xi),t).
\end{align}
By a change of variables, we may assume that 
a translation-invariant phase function $\phi(\x,\xi) = \langle x, \xi \rangle + \psi(t,\xi)$ satisfies one of the two alternatives:
\begin{align}
    \partial_t \nabla_\xi^2 \psi(0,0) &= I, \ &\text{(positive-definite)} \\
    \partial_t \nabla_\xi^2 \psi(0,0) &=I_-:= \begin{pmatrix}
        1 & 0 \\
        0 & -1
    \end{pmatrix}. &\text{(non-positive-definite)}
\end{align}

A compact set $K \subset \R^3$ is a $\phi$-\emph{Kakeya set} if, for every $\xi \in B_{\eps_{\phi}
}^2$, $K$ contains a curve $\ell_{\xi,v} \in \mathcal C(\phi)$. Let $d_{\mathrm{set}}(\phi)$ denote the infimum of the Hausdorff dimensions of such sets. A first step to understanding the $L^p$ estimates for H\"ormander-type oscillatory integral operators is to better understand the behavior of $\phi$-Kakeya sets.

The straight-line Kakeya problem and the restriction problem in $\R^3$ are obtained from the phase $\phi_{\mathrm{parab}}(\x,\xi) = \langle x,\xi\rangle + t |\xi|^2$. Recently Wang and Zahl proved $d_{\mathrm{set}}(\phi_{\mathrm{parab}})=3$, resolving the Kakeya conjecture in $\R^3$ \cite{wangKakeya}. 

Much less is known in the curved setting, though. There are certain $\phi$ so that $\phi$-Kakeya sets can lie in a surface, giving the worst possible bound $d_{\mathrm{set}}(\phi)=2$, and for generic $\phi$, small portions of $\delta$-discretized $\phi$-Kakeya sets can compress near a surface \cite{BourgainSeveralVariables}. There has been a long line of work in quantifying how poorly $\phi$-Kakeya sets can behave, in order to better understand the corresponding oscillatory integral operator, including but certainly not limited to \cite{minicozziSogge, soggeChaoticCurvature,  GuthSharpOscillatory, 
HormanderDichotomy, beyondUniversalEstimates, DaiOscillatory, dai2026localgeometryoscillatoryintegrals}.

In all cases, excluding those $\phi$ for which $\mathcal C(\phi)$ is simply diffeomorphic to a family of lines, the best unconditional estimates available were $d_{\mathrm{set}}(\phi) \geq 2+1/3$. 
For the very special class of phase functions satisfying \emph{Bourgain's condition}, the threshold $2+1/3$ can be obtained by a polynomial partitioning argument \cite{HickmanRogersZhang, DaiOscillatory}. See \cite{nadjimzadah2025bourgainsconditionstickykakeya} for an example of a (non-translation-invariant) phase function satisfying Bourgain's condition whose curves are not diffeomorphic to a family of lines. Earlier, Sogge attained Minkowski dimension $2+1/3$ for certain families of curves by using geometric methods \cite{soggeChaoticCurvature}. Theorem \ref{thm: main phicurved kakeya} surpasses the $2+1/3$ threshold for curved Kakeya sets in $\R^3$. 

\begin{theorem}[Main Theorem on $\phi$-curved Kakeya]\label{thm: main phicurved kakeya}
    Let $\phi(x,t,\xi) = \langle x, \xi \rangle + \psi(t,\xi)$ be a translation-invariant non-positive-definite phase function. Define 
    \begin{align}
        B &=  \tfrac{1}{2}\partial_t^2 \nabla_\xi^2 \psi(0,0), \\
        C &= \tfrac{1}{6} \partial_t^3 \nabla_\xi^2 \psi(0,0). 
    \end{align}
    Suppose that the matrices $B$ and $C$ satisfy
    \begin{align}\label{eq: main kakeya open condition}
        |B_{12}(B_{22} - B_{11}) + C_{12}| > \frac{1}{2} |(B_{22} + B_{11})(B_{22} - B_{11}) + C_{11} + C_{22}|.
    \end{align}
    Then 
    \begin{align}
        d_{\rm{set}}(\phi) \geq \frac{13 - \sqrt{13}}{4} \geq 2.348.
    \end{align}
    In fact the following stronger discretized statement holds. For all $\epsilon > 0$ there is an $M = M(\epsilon)$ and $\delta_0 = \delta_0(\epsilon)$ so that the following holds for all $\delta \in (0, \delta_0]$. Let $L \subset C(\phi)$ be a $\delta$-direction separated collection of curves. Let $E \subset \R^3$ be a union of $\delta$-cubes and suppose that $|\ell \cap E| \geq \lambda$ for each $\ell \in L$. Then 
    \begin{align}
        |E| \geq \delta^\epsilon \lambda^{M} \delta^{\frac{a + 1}{2}} (\delta^2 \# L)^{3/4}, 
    \end{align}
    where $a = (\sqrt{13} - 3)/2 \approx 0.303$ is the positive root of $a^2  + 3a - 1$.  
\end{theorem}

The hypothesis on $\phi$ in Theorem \ref{thm: main phicurved kakeya} is stable under translation-invariant perturbations. A simple example of a phase function satisfying the hypothesis is

\begin{align}
    \phi(x,t,\xi) = \langle x,\xi\rangle + \frac{1}{2}\langle \left( t I_- + t^2
    \begin{pmatrix}
        0 & 1 \\
        1 & 1
    \end{pmatrix} \right)
    \xi, \xi\rangle.
\end{align}

The proof of Theorem \ref{thm: main phicurved kakeya} combines Wolff's hairbrush argument \cite{wolffHairbrush} with a new \emph{stable} curved Kakeya estimate for 3-parameter families of curves (Theorem \ref{thm: CT 3 param kak}) in order to estimate the hairbrushes of curves. We will see that \eqref{eq: main kakeya open condition} is a natural condition that guarantees the hairbrushes of curves satisfy the conditions of Theorem \ref{thm: CT 3 param kak}: \emph{coniness} and \emph{twistiness}. 

\subsection{The 3-parameter curved Kakeya estimate}

The main new input is an incidence estimate for general 3-parameter families of curves. Let 
\begin{align}
    X : B_2^3 \times B_2^1 \to \R^2, \qquad\ell_\p(t) := (X(\p,t),t),
\end{align}
and write $\mathcal C(X) = \{\ell_\p : \p \in B_2^3\}$. We defer the formal definitions used in the theorem to Section \ref{sec: 2}, but give their geometric idea here. 

We say that $\mathcal C(X)$ is \emph{coney} (formally Condition \ref{cond: coney}) if the directions of the 1-parameter family of curves in $\mathcal C(X)$ through a point are not contained in a plane.

We now describe the \emph{twistiness} condition (formally Condition \ref{cond: twisty}). 
The curves intersecting a fixed curve in $\mathcal C(X)$ (at a relatively small angle) fill out a certain twisting \emph{grain}. The \emph{grain direction} can be projected out to give a plane curve. Applying the same procedure for each fixed curve in $\mathcal C(X)$ gives a $3$-parameter family of plane curves. We say that $\mathcal C(X)$ is \emph{twisty} if this family of curves forbids 2nd-order tangency in the sense of \cite{zahlPlanar}.
Twistiness is motivated by the structure of $\mathrm{SL_2}$-Kakeya sets in \cite{KWZ}. 

\begin{theorem}[Coney-Twisty 3-parameter Kakeya]\label{thm: CT 3 param kak}
    Suppose that $X$ is $\mathfrak c$-coney, $\tau$-twisty, and satisfies the basic conditions. For all $\epsilon > 0$ there is $M = M(\epsilon)$ and $\delta_0 = \delta_0(\epsilon, C_0,\ldots, C_{O_\epsilon(1)})$ so that the following holds for all $\delta \in (0, \delta_0]$. 
    Suppose that $L \subset \mathcal C(X)$ satisfies the two-dimensional ball condition
    \begin{align}
        \# (L \cap B) \leq (r/ \delta)^2 \text{ for all $r \in [\delta, 1]$, and all balls $B$ of radius $r$.}
    \end{align}
    Let $E \subset \R^3$ be a union of $\delta$-cubes, and suppose $|\ell \cap E| \geq \lambda$ for each $\ell \in L$. Then 
    \begin{align}
        |E| \geq  \delta^\epsilon \mathfrak c^2 (\tau \lambda)^M \delta^a (\delta^2 \#L), 
    \end{align}
    where $a = (\sqrt{13} - 3)/2 \approx 0.303$ is the positive root of $a^2 + 3a - 1$.
\end{theorem}

Theorem \ref{thm: CT 3 param kak} has implications to restricted projections, which we record in Section \ref{sec: 2}. The new feature of Theorem \ref{thm: CT 3 param kak} is that its hypotheses are stable under small perturbations of the family of curves. 
The proof uses the techniques in \cite{KWZ} when the typical angle between curves are small, and the curved multilinear Kakeya theorem \cite[Theorem 6]{BourgainGuth} when the typical angle is large. Below we give a sketch of the proof of Theorem \ref{thm: main phicurved kakeya} by using Theorem \ref{thm: CT 3 param kak}.

\subsection{Sketch of the proof of Theorem \ref{thm: main phicurved kakeya}}\label{subsec: sketch of pf}

We give a sketch of the proof of Theorem \ref{thm: main phicurved kakeya} by tracking the example $\phi_A(\x,\xi) = \langle x,\xi \rangle + \tfrac{1}{2}\langle A(t)\xi, \xi \rangle$, where $\det A'(0) \neq 0$.
There are two main ingredients. 
\begin{enumerate}
    \item Wolff's classical hairbrush argument lets us extract a hairbrush configuration $H \subset L$ to study. 
    \item Theorem \ref{thm: CT 3 param kak} lets us estimate the union of curves in $H$.
\end{enumerate}

For the sketch, we consider the special case $\lambda = 1$ and $\# L = \delta^{-2}$. Write $\mathbb T$ for the collection of $\delta$-neighborhoods of curves in $L$. We can locate a multiplicity $\mu$ between tubes in $\mathbb T$. By double counting, 
\begin{align}\label{eq: hairbrush sketch1}
    \Big |\bigcup_{T \in \mathbb T} T\Big | \gtrsim 1/\mu. 
\end{align}
We also need an estimate which gets better with large $\mu$. Since the average multiplicity is $\mu$, there is a large subset of tubes $\mathbb T' \subset \mathbb T$ which are incident to $\mu \delta^{-1}$ many other tubes. Take one of these tubes $T_0 \subset \mathbb T$, and for the sketch assume that the tubes in the hairbrush $\mathbb H$ around $T_0$ make an angle $\sim 1$ with $T_0$. If, for instance, $T_0$ is the tubes with center curve $\ell_{0,0}$, then the center curves in $\mathbb H$ are parameterized like 
\begin{align}\label{eq: intro: hairbrush param}
    \ell_{\xi, \theta}(t) = ((A(t) - A(\theta))\xi, t).
\end{align}
The parameter $\theta \in \R$ is the height at which the curve hits $\ell_{0,0}$, and $\xi \in \R^2$ is the curve's direction (these are distinct from $\ell_{\xi,v} \in \mathcal C(\phi)$, where $v \in \R^2$). Unfortunately Wolff's original argument to bound $\bigcup_{T \in \mathbb H} T$ does not work. This was first observed by Wisewell in the special case of quadratic phase functions: 
\begin{observation}[Wisewell \protect{\cite[Proposition 14]{Wisewell}}]
    Even in the special case of $A(t) = tA + t^2 B$, Wolff's hairbrush argument works only when $B$ is a multiple of $A$. Equivalently, $\mathcal C(\phi_A)$ consists of lines up to a change of variables. 
\end{observation}
The problem is that Wolff's argument uses the rare property that the curves in a hairbrush can be organized into essentially disjoint surfaces. We instead can only rely on a property of curves which is stable under perturbations. These are Conditions \ref{cond: coney} and \ref{cond: twisty}. When $\phi$ satisfying the hypothesis of Theorem \ref{thm: main phicurved kakeya}, Theorem \ref{thm: CT 3 param kak} says 
\begin{align}\label{eq: hairbrush sketch2}
    \Big |\bigcup_{T \in \mathbb H} T \Big | \gtrsim \delta^a ( \delta^2 \# \mathbb H) \geq \delta^{a + 1} \mu, 
\end{align}
where $a \approx 0.303$.
By averaging the bounds \eqref{eq: hairbrush sketch1} and \eqref{eq: hairbrush sketch2}, we find 
\begin{align}
    \Big |\bigcup_{T \in \mathbb T} T\Big | \gtrsim \delta^{(a + 1)/2},
\end{align}
corresponding to dimension $3 - (a + 1)/2 \approx 2.348$. 

We made simplifying assumptions at various points in this sketch. First, we cannot assume $\lambda = 1$. To handle this, we will use the two-ends reduction which introduces an isotropic rescaling of the tubes. We also cannot assume the typical hairbrush angle is $\sim 1$. To handle this, we need a radial rescaling of a hairbrush and the bush argument to estimate the union of narrow hairbrushes. These technicalities require the conditions for Theorem \ref{thm: CT 3 param kak} to be fairly flexible.

\subsection{Outline of the article}
Section \ref{sec: 2} gives the formal geometry of 3-parameter families, and records two further 3-parameter Kakeya estimates and various consequences for non-linear projections. Section \ref{sec: 3} introduces the discretization and uniformization tools. Sections \ref{sec: 4}--\ref{sec: 6} prove Theorem \ref{thm: CT 3 param kak}. Sections \ref{sec: 7} and \ref{sec: 8} prove Theorem \ref{thm: main phicurved kakeya}. Appendix \ref{appendix: A} contains the proofs of the auxiliary 3-parameter Kakeya estimates, and Appendix \ref{appendix: B} records the symbolic Mathematica calculations used to check the coney and twisty conditions in special cases.

\subsection{Acknowledgments}

The author would like to thank his advisors Terence Tao and Hong Wang for their support and many helpful conversations throughout the course of this project, and for suggestions that greatly improved the presentation in this article. The author would also like to thank the anonymous referees whose suggestions further improved the presentation.

\section{Geometry of 3-parameter families of curves}\label{sec: 2}

We now give the formal definitions used in Theorem \ref{thm: CT 3 param kak}. 

\begin{definition}[3-parameter family of curves]\label{def: 3 param family}
    Let $X(\p,t) : B_2^3 \times B_2^1 \to \R^2$ be a smooth map which we call the \emph{defining function}. The point $\p \in B_2^3$ determines the curve segment $\ell_\p$ with parameterization 
    \begin{align}
        \ell_\p(t) = (X(\p,t),t) 
    \end{align}
    and $t \in B_2^1$. 
    Write $\mathcal C(X)$ for the $3$-parameter family of these curves.
\end{definition}

For example, the hairbrush family of curves \eqref{eq: intro: hairbrush param} has defining function 
\begin{align}
    X(\p,t) = (A(t) - A(\theta))\xi, 
\end{align}
where $\p = (\xi,\theta)$.
The geometry of curves intersecting a fixed curve will play a central role in our arguments. This leads to the definition of \emph{spread curves}.

\begin{figure}
    \centering

\tikzset{every picture/.style={line width=0.75pt}} 

\begin{tikzpicture}[x=0.75pt,y=0.75pt,yscale=-1,xscale=1]

\draw (327.42,229.19) node  {\includegraphics[width=218pt,height=190pt]{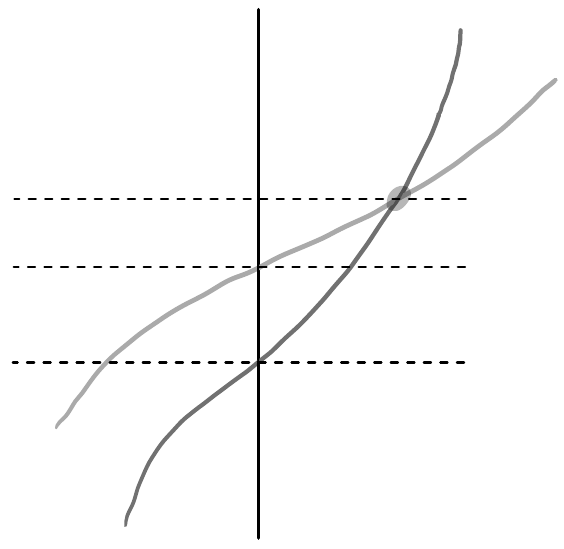}};
\draw [line width=1.5]    (389.08,196.43) -- (417.55,181.31) ;
\draw [shift={(421.08,179.43)}, rotate = 152.02] [fill={rgb, 255:red, 0; green, 0; blue, 0 }  ][line width=0.08]  [draw opacity=0] (11.61,-5.58) -- (0,0) -- (11.61,5.58) -- cycle    ;

\draw (472,110.4) node [anchor=north west][inner sep=0.75pt]  [font=\large]  {$\ell _{\hat{\xi }_{\xi ,\theta ,t}( \theta ') ,\theta '}$};
\draw (159,185.4) node [anchor=north west][inner sep=0.75pt]    {$t$};
\draw (156,263.4) node [anchor=north west][inner sep=0.75pt]    {$\theta $};
\draw (155,219.4) node [anchor=north west][inner sep=0.75pt]    {$\theta '$};
\draw (404,79.4) node [anchor=north west][inner sep=0.75pt]  [font=\large]  {$\ell _{\xi ,\theta }$};
\draw (297,70.4) node [anchor=north west][inner sep=0.75pt]  [font=\large]  {$\ell _{0,0}$};
\draw (430,182.4) node [anchor=north west][inner sep=0.75pt]    {$(\omega_{\xi,\theta,t}( \theta'),1)$};
\end{tikzpicture}

    \caption{An illustration of the spread curve $\omega_{\xi, \theta,t}$ evaluated at a point $\theta'$ (for the special case of curves in a hairbrush). The stem line $\ell_{0,0}$ is hit by $\ell_{\xi,\theta}$ (at height $\theta)$. The lightest curve is the unique curve hitting $\ell_{0,0}$ at height $\theta'$ and $\ell_{\xi,\theta}$ at height $t$. It is determined by the direction $\hat \xi_{\xi, \theta, t}(\theta')$, and height $\theta'$ along the stem (in the general setting of Definition \ref{def: spread curves}, this is $\hat p_{\p,t}(\theta')$ and $\theta'$). 
    The vector shown is $(\omega_{\xi,\theta,t}(\theta'), 1)$, the direction of the lightest curve $\ell_{\hat \xi_{\xi, \theta,t}}(\theta')$ at height $t$.}
    \label{fig: spread curves}
\end{figure}

\begin{definition}[Spread curves]\label{def: spread curves}
     Write $\p = (p,\theta) \in \R^2 \times \R$ and define the curves $\omega_{\p, t}(\theta') = \partial_t X(\hat p_{\p,t}(\theta'), \theta', t)$, where $\hat p_{\p,t}$ solves
    \begin{align}
        X(\hat p_{\p,t}(\theta'),\theta', t) = X(\p, t) 
    \end{align}
    near $\theta$.
\end{definition}
Figure \ref{fig: spread curves} is a picture of Definition \ref{def: spread curves} for the hairbrush family example given above (with $p$ replaced by $\xi$ and $\hat p$ replaced by $\hat \xi$ to match notation for $\phi$-curves in a hairbrush).

Without some mild assumptions, it is not guaranteed that these spread curves exist. We will describe some more structural properties below as well, and we will collect all the mild assumptions needed to guarantee that they exist at the end.
Next we can formalize whether the curves through a point spread into a plane or a cone (with quantified curvature). The next two conditions use the spread curves to do this. 

\begin{condition}[Planiness]\label{cond: planey}\hfill
    \begin{itemize}
        \item (Planey at a point): The map $X$ is \emph{planey at $(\p,t)$} if 
        \begin{align}
            \det(\dot \omega_{\p,t}(\theta), \ddot \omega_{\p,t}(\theta)) = 0.
        \end{align}
        \item (Planey): The map $X$ is \emph{planey} if it is planey at each point $(\p,t) \in B_2^3 \times B_2^1$. 
    \end{itemize}
\end{condition}

\begin{condition}[$\mathfrak c$-coney] \label{cond: coney}\hfill
    \begin{itemize}
        \item ($\mathfrak c$-coney at a point): The map $X$ is \emph{$\mathfrak c$-coney at $(\p,t)$} if 
        \begin{align}
            |\det(\dot \omega_{\p,t}(\theta), \ddot \omega_{\p,t}(\theta))| \geq \mathfrak c. 
        \end{align}
        \item ($\mathfrak c$-coney): The map $X$ is \emph{$\mathfrak c$-coney} if it is $\mathfrak c$-coney at each point $(\p,t) \in B_2^3 \times B_2^1$. 
    \end{itemize}
\end{condition}
If $X$ is $\Omega(1)$-coney, then it is $\Omega(1)$-coney under small perturbations of $X$. We therefore say $\Omega(1)$-coniness is a \emph{stable} property. 

The family $\mathcal C(X)$ enjoys an explicit graininess property. Fix $\ell_{\p} \in \mathcal C(X)$, a scale $\delta > 0$, and $\rho \in [\delta, \delta^{1/2}]$. Then the union of curves passing through $\ell_{\p}$ at angle $\sim r$ inside $N_\rho(\ell_\p)$ 
\begin{align}\label{eq: union at angle r}
    N_\rho(\ell_\p) \cap \bigcup_{\angle(\ell', \ell_\p) \sim r} \ell',
\end{align}
is contained in a $\delta \times \rho \times 1$ (possibly twisted) grain, as long as $\rho$ is small enough depending on $\delta$ and $r$. One expects that \eqref{eq: union at angle r} is contained in a single grain, because the family of curves $\ell'$ through a point on $\ell_\p$ is a 1-parameter family, and the family varies smoothly along $\ell_\p$. This will be made rigorous in Section \ref{sec: 5}. The $\rho$-direction of the grain is determined by 
\begin{align}
    \gamma_\p(t) := \dot \omega_{\p,t}(\theta). 
\end{align}
When $X$ is planey $\rho$ can be taken larger, so the grain becomes wider. When $X$ is coney, $\rho$ cannot be taken too large. We will find that the $\rho$-direction of nearby grains are almost the same. We can therefore describe the grains in $N_\rho(\ell)$ by the following grain projection map.

\begin{definition}[Grain Projection]\label{def: grain projection}
    Fix $\p \in B_2^3$ and define the grain projection $\pi_\p : N_\rho(\ell_\p) \to N_\rho(B_2^1)$ by 
    \begin{align}
        \pi_{\p}(x,t) = (t, (x - X(\p,t)) \cdot \gamma_{\p}(t)^\perp). 
    \end{align}
\end{definition}
The grain through $\ell_\p$ which was described above is $\pi_\p^{-1}(N_\delta(B_2^1))$. The map $\pi_\p$ ``projects out'' the common grain directions. 
Definition \ref{def: grain projection} is a generalization of the \emph{twisted projections} for $\rm{S}L_2$-lines in \cite{KWZ}. The authors in \cite{KWZ} needed control over the overlap of nearby grains. Equivalently, they needed control over the grain projections of curves in $N_{\rho}(\ell_\p)$. Below we define the family of plane curves in our general setting.

\begin{definition}
    Define the family of plane curves 
    \begin{align}
        \mathcal F(\p, X) = \{\pi_\p(\ell_{\p'}) : \p' \in B(\p, \rho)\}. 
    \end{align}
    As long as $\rho \leq \delta^{1/2}$ (which we will always assume), Taylor's theorem says $\mathcal F(\p,X)$ can be taken to be
    \begin{align}
        \{\mathrm{Graph}(h_\p(\p',\cdot)) : p'\in B(\p, \rho)\}
    \end{align}
    up to error $\delta$, where 
    \begin{align}
        h_\p(\p', \cdot) &= \langle \nabla_\p X(\p,t) (\p' - \p), \gamma_\p(t)^\perp \rangle. 
    \end{align}
\end{definition}

As long as $\mathcal F(\p,X)$ forbids 2nd-order tangency in the sense of \cite[Definition 1.8]{zahlPlanar}, then the planar maximal function estimate \cite[Theorem 1.9]{zahlPlanar} holds. To quantify forbidding 2nd-order tangency, we introduce the \emph{tangency matrix}: 

\begin{definition}\label{def: tang matrix}
    Define the map
    \begin{align}
        F_{\p,t}(\p') &= (h_\p(\p',t), \partial_t h_\p(\p',t), \partial_t^2 h_\p(\p',t)),
    \end{align}
    as in \cite[Definition 1.1]{zahlPlanar}.
    This is given by $M_{\p}(t) \cdot (\p' - \p)$, 
    where $M_{\p,t}$ is the \emph{tangency matrix} 
    \begin{align}\label{eq: tangency matrix}
        M_{\p}(t) &= 
        \begin{pmatrix}
        \nabla_\p X(\p,t)^T \gamma_\p(t)^\perp \\
        \partial_t[\nabla_\p X(\p,t)^T \gamma_\p(t)^\perp] \\
        \partial_t^2[\nabla_\p X(\p,t)^T \gamma_\p(t)^\perp] 
    \end{pmatrix}.
    \end{align}
\end{definition}
We may now state the $\tau$-twisty property of $X$: 
\begin{condition}[$\tau$-twisty property]\label{cond: twisty}
    \begin{itemize}\hfill
        \item ($\tau$-twisty at a point): The map $X$ is $\tau$-twisty at $(\p,t)$ if 
        \begin{align}
            |\det M_{\p}(t)| \geq \tau. 
        \end{align}
        \item ($\tau$-twisty): The map $X$ is $\tau$-twisty if it is $\tau$-twisty at every point $(\p,t) \in B_2^3 \times B_2^1$.
    \end{itemize}
\end{condition}
When $X$ is $\Omega(1)$-twisty, then small perturbations of $X$ are also $\Omega(1)$-twisty. So $\Omega(1)$-twistiness is a stable property. It is also meaningful to ask what happens when the $3 \times 3$ matrix $M_{\p}(t)$ has rank 1 or 2. In the former case, we say $X$ is $1$-flat at $(\p,t)$, and in the latter we say $X$ is $2$-flat at $(\p,t)$. When these properties hold for all $(\p,t)$, we say $X$ is $1$-flat and $2$-flat respectively. Considering $\pi_\p$ as a map from space curves to plane curves: 
\begin{itemize}
    \item Twistiness says the projections are injective.
    \item $1$-flat says the image of the curves is a 1-parameter family (and the preimages lie in a 2-parameter family). 
    \item $2$-flat says the image of curves is a 2-parameter family (and the preimages lie in a 1-parameter family).
\end{itemize}

Below is a collection of the mild nondegeneracy conditions needed to make all of the conditions above effective.

\begin{condition}[Basic Conditions]\label{cond: basic}
    We say $X$ satisfies the basic conditions if there exists a sequence of constants $(C_i)_{i=0}^\infty$ with $C_0 = C_1 = C_2$ so that the following holds. For all $(\p=(p,\theta),t) \in B_2^3 \times B_2^1$, 
    \begin{enumerate}
        \item \label{cond: basic nondeg}(Nondegeneracy) The entries of $\nabla_\p X$ and $(\partial_{p_i,p_j} X)_{1\leq i,j \leq 2}$ are $\leq C$ and 
        \begin{align}
            |\det \nabla_p X(\p,t)| \geq C^{-1}. 
        \end{align}
        \item \label{cond: basic angle}(Angle faithfulness)
        \begin{align}
            C^{-1} \leq |\gamma_\p(t)| \leq C.
        \end{align}
        \item \label{cond: basic coniness}(Relevant to coniness)
        \begin{align}
            |\partial_t^i X(\p,t)| \leq C_i.
        \end{align}
        \begin{align}
            |\ddot \omega_{\p,t}(\theta)|, |\dddot \omega_{\p,t}(\theta)|  \leq C.
        \end{align}
        \item \label{cond: basic twistiness}(Relevant to twistiness)
        For each $i \geq 0$, 
        \begin{align}
            |\partial_t^i[\nabla_{\p} X(\p,t)^T \gamma_\p(t)^\perp]| \leq C_i.
        \end{align}
    \end{enumerate}
\end{condition}

\begin{remark}
    If one does not care how the constants $C_i$ depend on $X$, then by compactness it is enough to check 
    \begin{enumerate}
        \item $\det \nabla_p X \neq 0$ and 
        \item $\gamma_{\p}(t) \neq 0$.
    \end{enumerate}
    However in our application to Theorem \ref{thm: main phicurved kakeya}, the more precise control is crucial. 
\end{remark}

Below we check the conditions we developed on a few familiar examples. 
\begin{example}[Straight hairbrush is planey and $1$-flat]\label{ex: straight hairbrush}
    The curves in the classical hairbrush may be parameterized by 
    \begin{align}
        \ell_{\xi,\theta}(t) &= ((\theta-t) \xi, t),
    \end{align}
    so the defining map is $X(\xi, \theta,t) = (\theta-t)\xi$. Assume that $1/2 \leq |t- \theta| \leq 1$ and $1/2 \leq |\xi| \leq 1$.
    The spread-curves are 
    \begin{align}
        \omega_{\xi,\theta,t}(\theta') &= \frac{t-\theta}{\theta'-t} \xi.
    \end{align}
    It is straightforward to check that $X$ is planey, as expected. 
    Also $\gamma_{\xi,\theta}(t) = (\theta-t)^{-1} \xi$. This corresponds to the grains around the line in direction $(\xi, 1)$ running parallel to the plane $\mathrm{Span}((\xi,0), (0,0,1))$. 
    It is also straightforward to compute the tangency matrix: 
    \begin{align}
        M_{\xi,\theta}(t) &= 
        \begin{pmatrix}
            \xi^\perp & 0  \\
            0 & 0 \\
            0 & 0
        \end{pmatrix},
    \end{align}
    which has rank 1. This corresponds to the fact that lines in a fixed $\delta^{1/2}$-tube in a hairbrush are organized into parallel planes. Of course, we know this is true globally for a hairbrush of straight lines.
\end{example}

We also check the conditions on the family of $SL_2$-lines considered by Katz, Wu, and Zahl. 
\begin{example}[$\mathrm{SL}_2$-family is planey and $\Omega(1)$-twisty]\label{ex: SL2 family}
    Define the $\mathrm{SL}_2$-family of lines by 
    \begin{align}
        \mathcal L_{\mathrm{SL}_2} = \{\ell_{a,b,c,d} : ad - bc = 1\},
    \end{align}
    where $\ell_{a,b,c,d} = [(a,b,0) + \mathrm{span}(c,d,1)] \cap B_1^3$. 
    Even though $\mathcal L_{\mathrm{SL}_2}$ is not globally given by a defining function, we can cover it by charts and put a defining function on each chart. After dealing with this technicality, $\mathcal L_{\mathrm{SL}_2}$ is planey and $\Omega(1)$-twisty. The family in fact enjoys a stronger form of planiness \cite[Section 1.2.2]{KWZ} and a stronger form of twistiness \cite[Lemma 2.5]{KWZ}. 
\end{example}

We will prove in Section \ref{sec: 8} that phase functions satisfying the conditions of Theorem \ref{thm: main phicurved kakeya} are coney and twisty. As a very special case: 

\begin{example}[A stable collection of coney and twisty hairbrushes]
    The family of hairbrush curves 
    \begin{align}
        \ell_{\xi,\theta}(t) = (((t -\theta)A +(t^2 - \theta^2)B)\xi, t),
    \end{align}
    with $A = I_-$, $B$ symmetric, $B_{11} \neq B_{22}$, and $|B_{12}| > \frac{1}{2} |B_{11} + B_{22}|$
    is $\Omega(1)$-coney and $\Omega(1)$-twisty. 
\end{example}

\subsection{Further results on 3-parameter Kakeya and projection consequences}\label{subsec: further results}

The proof of Theorem \ref{thm: CT 3 param kak} uses both coniness and twistiness. We record two further variants. One uses only coniness, and the other is a mild generalization of the result in \cite{KWZ} using planiness and twistiness. These variants are not needed for Theorem \ref{thm: main phicurved kakeya}, but they may be useful elsewhere, and they provide a useful comparison.

\begin{theorem}[Coney 3-parameter Kakeya]\label{thm: C 3 param kak}
    Suppose that $X$ is $\mathfrak c$-coney and satisfies the basic conditions. For all $\epsilon > 0$ there is $\delta_0 = \delta_0(\epsilon, C_0,\ldots, C_{O_\epsilon}(1))$ so that the following holds for all $\delta \in (0,\delta_0]$. Suppose that $L \subset \mathcal C(X)$ satisfies the two-dimensional ball condition 
    \begin{align}
        \# (L \cap B) \leq (r/ \delta)^2 \text{ for all $r \in [\delta, 1]$, and all balls $B$ of radius $r$.}
    \end{align}
    Let $E \subset \R^3$ be a union of $\delta$-cubes, and suppose $|\ell \cap E| \geq \lambda$ for each $\ell \in L$. Then 
    \begin{align}
        |E| \geq \delta^\epsilon \mathfrak c^{3/2} \lambda^2 \delta^{1/2} (\delta^2 \#L). 
    \end{align}
\end{theorem}

The bound is worse than Theorem \ref{thm: CT 3 param kak}, but the dependence on $\lambda$ is much better.

\begin{theorem}[Planey-Twisty 3-parameter Kakeya (minor generalization of \protect{\cite[Theorem 1.2]{KWZ}})]\label{thm: PT 3 param kak}

    Suppose that $X$ is planey, $\tau$-twisty, and satisfies the basic conditions. For all $\epsilon > 0$ there is $M = M(\epsilon)$ and $\delta_0 = \delta_0(\epsilon, C_0,\ldots, C_{O_\epsilon(1)})$ so that the following holds for all $\delta \in (0,\delta_0]$. Suppose that $L \subset \mathcal C(X)$ satisfies the two-dimensional ball condition.
    Let $E \subset \R^3$ be a union of $\delta$-cubes, and suppose $|\ell \cap E| \geq \lambda$ for each $\ell \in L$. Then 
    \begin{align}
        |E| \geq \delta^\epsilon (\tau \lambda)^M  (\delta^2 \# L).
    \end{align}
\end{theorem}
The proof of Theorem \ref{thm: C 3 param kak} and an outline of the proof of Theorem \ref{thm: PT 3 param kak} are contained in Appendix \ref{appendix: A}.

By the same discretization argument that Katz, Wu, and Zahl used to convert their sharp $\rm{SL}_2$-Kakeya result (\cite[Theorem 1.2]{KWZ}) to the linear restricted projection result (\cite[Theorem 1.3]{KWZ}, first proved using decoupling in \cite{ganRestrictedProjectionsPlanes}), we get non-linear restricted projection results as corollaries of Theorems \ref{thm: CT 3 param kak}, \ref{thm: C 3 param kak}, and \ref{thm: PT 3 param kak}. 

\begin{corollary}[Corresponding Restricted Projection Theorems]
    Let $F_t : \R^3 \to \R^2$ be a family of maps, and define $X(\p,t) = F_t(\p)$. Let $A$ be a Borel subset of $\R^3$ with $\dim(A) = 2$. Then 
    \begin{itemize}
        \item If $X$ is $\Omega(1)$-coney and $\Omega(1)$-twisty, then 
        \begin{align}
            \dim F_t(A) \geq 1.697  \text{ for a.e. $t \in [0,1]$}.
        \end{align}
        \item If $X$ is $\Omega(1)$-coney,
        then 
        \begin{align}
            \dim F_t(A) \geq 1.5 \text{ for a.e. $t \in [0,1]$}.
        \end{align}
        \item If $X$ is planey and $\Omega(1)$-twisty, then 
        \begin{align}
            \dim F_t(A) = 2  \text{ for a.e. $t \in [0,1]$}.
        \end{align}
    \end{itemize}
\end{corollary}

\begin{remark}
    The first bound is most interesting, seeing as $\Omega(1)$-coniness and $\Omega(1)$-twistiness are stable properties. In order to prove a stable nonlinear projection result, it is tempting to apply Shmerkin's approximation method from \cite{shmerkinNonlinearBourgain} directly to the linear restricted projection theorem in \cite{ganRestrictedProjectionsPlanes}. The worst case example for that approach is a well-spaced collection of $\delta^{-2}$ points
\begin{align}
    A = (\delta^{2/3} \mathbb Z)^3 \cap B_1^3.
\end{align}
Let $m = |\log \delta|$. The slope function for $A$ (as in \cite{shmerkinNonlinearBourgain}) is 
\begin{align}
    f(x) = 
    \begin{cases}
        3x, & x \in [1,(2/3)m] \\
        2, & x \in ((2/3)m,m].
    \end{cases}
\end{align}
Cut $[1,m]$ into intervals $[1,x]$, $[x,y] \ni (2/3)m$, $[y,m]$. Dimension 2 is critical for the linear restricted projection theorem, so it is optimal to choose $x,y$ such that $(f(y) - f(x))/(y-x) = 2$. To linearize $F_t$ on $[x,y]$, we need $y \leq 2x$. We should thus take $y = 2x$ to make full use of the linearization. These equations determine $x = (2/5)m$, $y=(4/5)m$. We have $m \leq 2y$ so we may linearize on $[y,m]$. Cut $[1,x]$ into $O(\log m)$ many intervals and linearize on each. Applying the linear restricted projection theorem on each interval, we find
\begin{align}
    \log(|F_t(A)|_\delta) \gtrsim 2 (2/5 - 0)m + 2(4/5 - 2/5)m + 0(1-4/5)m = 1.6m,
\end{align}
where $|F_t(A)|_\delta$ is the number of $\delta$-balls needed to cover $F_t(A)$.
This corresponds to $\dim F_t(A) \geq 1.6$ for a.e. $t \in [0,1]$, which is worse than our result. 
\end{remark}

\section{Notation, Discretization, and Pigeonholing/Uniformization Tools}\label{sec: 3}

\subsection{Asymptotic notation}
The relation $A \lesssim B$ means that $A \leq C B$, where $C$ depends only on finitely many of the constants $C_i$ in Condition \ref{cond: basic} (when $X$ is relevant), or finitely many derivatives of $\phi$ (when $\phi$ is relevant). We write $A \lesssim_\epsilon B$ if $C$ might additionally depend on $\epsilon$. For instance, $C$ might depend on $\epsilon,C_0,\ldots, C_{O_\epsilon(1)}$. When we write $A \lessapprox^\delta B$, we mean that 
\begin{align}
    A \lesssim 2^{100|\log \delta|/\log |\log \delta|} B.
\end{align}
We will only need the full strength of this notation in the proofs of Theorem \ref{thm: CT 3 param kak} and \ref{thm: C 3 param kak}, 
and elsewhere the reader can replace it with $A \lesssim |\log \delta|^{100} B$ if they wish. 
When the small parameter $\delta$ is clear from context, we will write $A \lessapprox B$. We will adopt the notation from \cite{KWZ} and write $A = O(B)$ to mean $A \lesssim B$, $A = \Omega(B)$ to mean $A \gtrsim B$, $A = O^*(B)$ to mean $A \lessapprox^\delta B$, and 
$A = \Omega^*(B)$ to mean $A \gtrapprox^\delta B$. 

\subsection{Distances and angles in $\cC(X)$}

We give $\cC(X)$ the metric induced by $\R^3$: 
\begin{definition}\label{def: metric on C(X)}
    The distance between $\ell_{\p_1}$ and $\ell_{\p_2}$ in $\cC (X)$ is 
    \begin{align}
        d(\ell_{\p_1}, \ell_{\p_2}) = |\p_1 - \p_2|. 
    \end{align}
\end{definition}
Thanks to Part 1 of Condition \ref{cond: basic}, this metric reflects the physical distance between curves as expected.
\begin{lemma}\label{lem: metric on C(X)}
    There is a constant $C = C(C_0)$ and $r_0 = r_0(C_0)$ so that the following hold. 
    \begin{enumerate}
        \item If $d(\ell_{\p_1}, \ell_{\p_2}) \leq r$ then  $\ell_{\p_1} \subset N_{Cr}(\ell_{\p_2})$. Conversely if $r \le  r_0$ and $\ell_{\p_1} \subset N_{r}(\ell_{\p_2})$ then $d(\ell_{\p_1}, \ell_{\p_2}) \leq Cr$. 
        \item If $r \leq r_0$ and $|\ell_\p(t) - \x| \leq r$, then there is a $\p' \in B(\p, Cr)$ such that 
        \begin{align}
            \ell_{\p'}(t) = \x. 
        \end{align}
    \end{enumerate}
\end{lemma}

\begin{proof}
    By the fundamental theorem of calculus, the first direction of (1) is true. The second direction of (1) follows from a quantitative version of the implicit function theorem, such as the one found in these notes \cite{LiveraniImplicit}. 
    Item (2) also follows from the implicit function theorem. These use all of Condition \ref{cond: basic} Part \ref{cond: basic nondeg}. We leave the straightforward details to the reader. 
\end{proof}

The parameter $\theta$ in $\p = (p,\theta)$ gives a way to measure angles of intersecting curves. Let $\bn_{\p}(t) \in S^2$ be the direction of $\ell_\p$ at position $t$. 
\begin{lemma}\label{lem: angles in C(X)}
    Take $\ell_{\p_1}, \ell_{\p_2} \in \cC(X)$ which intersect at height $t$. Then 
    \begin{align}
        \angle(\bn_{\p_1}(t), \bn_{\p_2}(t)) \sim |\theta_1 - \theta_2|.
    \end{align}
\end{lemma}

\begin{proof}
    By the mean value theorem, $\omega_{p,\theta_1,t}(\theta_2) - \omega_{p,\theta_1,t}(\theta_1) \sim (\theta_2 - \theta_1) \eta$, where $|\eta| \sim 1$ by Condition \ref{cond: basic} Part 2. Therefore 
    \begin{align}
        |(\omega_{p,\theta_1,t}(\theta_2), 1) \times (\omega_{p,\theta_2,t}(\theta_1),1)| &= |\theta_2 - \theta_1| |(-\eta^2,-\eta^1, \omega_{p,\theta_1,t}^1(\theta_1)\eta^2 - \omega^2_{p,\theta_1,t}(\theta_1)\eta^1)| \\ &\sim |\theta_2 - \theta_1|. 
    \end{align}
    Therefore 
    \begin{align}
        \sin (\angle(\bn_{\p_1}(t), \bn_{\p_2}(t))) &= \frac{|(\omega_{p,\theta_1,t}(\theta_2), 1) \times (\omega_{p,\theta_1,t}(\theta_1), 1)|}{|(\omega_{p,\theta_1,t}(\theta_2), 1)||(\omega_{p,\theta_1,t}(\theta_1), 1)|} \\&\sim |\theta_2 - \theta_1|,
    \end{align}
    which proves the claim. 
\end{proof}

\subsection{Distances and angles in $\cC(\phi)$}

The metric on $\cC(\phi)$ is induced by $\R^4$: 
\begin{definition}
    The distance between $\ell_{\xi_1,v_1}$ and $\ell_{\xi_2,v_2}$ in $\cC(\phi)$ is 
    \begin{align}
        d(\ell_{\xi_1,v_1}, \ell_{\xi_2,v_2}) = |(\xi_1,v_1) - (\xi_2,v_2)|. 
    \end{align}
\end{definition}
Similarly, the metric reflects distances between curves in space as expected: 
\begin{lemma}\label{lem: metric on C(phi)}
    There is a constant $C = C(\phi)$ and $r_0 = r_0(\phi)$ so that the following hold. 
    \begin{enumerate}
        \item If $d(\ell_{\xi_1,v_1}, \ell_{\xi_2,v_2}) \leq r$ then  $\ell_{\xi_1,v_1} \in N_{Cr}(\ell_{\xi_2,v_2})$. Conversely if $r \leq  r_0$ and $\ell_{\xi_1,v_1} \in N_{r}(\ell_{\xi_2,v_2})$ then $d(\ell_{\xi_1,v_1}, \ell_{\xi_2,v_2}) \leq Cr$. 
        \item If $r \leq r_0$ and $|\ell_{\xi,v}(t) - \x| \leq r$, then there is a $(\xi',v') \in B((\xi,v), Cr)$ such that 
        \begin{align}
            \ell_{\xi',v'}(t) = \x. 
        \end{align}
    \end{enumerate}
\end{lemma}

\begin{proof}
    To run the implicit function theorem, we only need the phase function $\phi$ to satisfy the H\"ormander conditions.
\end{proof}

Write $\bn_{\xi,v}(t) \in S^2$ for the direction of $\ell_{\xi,v}$ at position $t$. 
It is straightforward to show that the parameter $\xi$ describes the angle: 
\begin{lemma}
    If $\ell_{\xi_1,v_1}, \ell_{\xi_2,v_2} \in \cC(\phi)$ intersect at height $t$ then 
    \begin{align}
        \angle(\bn_{\xi_1,v_1}(t), \bn_{\xi_2,v_2}(t)) \sim |\xi_1 - \xi_2|. 
    \end{align}
\end{lemma}

\subsection{Shadings}

Tile $\R^3$ by grid-aligned $\delta$-cubes. 
In both $\cC(X)$ and $\cC(\phi)$, define a shading of a curve $\ell$ to be a union $Y(\ell)$ of cubes intersecting $\ell$. The shading is $\lambda$-dense if $|Y(\ell)| \geq \lambda \delta^2$. 
Write $(L,Y)$ for a collection of curves with their associated $\lambda$-dense shadings. Define 
\begin{align}
    E_{(L,Y)} = \bigcup_{\ell \in L} Y(\ell). 
\end{align}
When the shading $Y$ of $L$ is clear from context, we write $E_L$. A collection of curves $L$ is $\delta$-separated if the corresponding points in $\R^3$ or $\R^4$ are $\delta$-separated, for $\cC(X)$ or $\cC(\phi)$ respectively. Recall from the introduction that $L \subset \cC(\phi)$ is $\delta$-direction separated if the points $\xi(\ell) \in \R^2$ are $\delta$-separated. Here we defined $\xi(\ell_{\xi,v}) = \xi$.

\subsection{Quantitative Broad-Narrow Reduction}
We restrict our attention here to $\cC(X)$. Define 
\begin{align}
    L(\x) = \{\ell \in L : \x \in Y(\ell)\}.
\end{align}
Suppose that $\ell_{\p_1}, \ell_{\p_2} \in L(\x)$. In light of Lemmas \ref{lem: metric on C(X)} and \ref{lem: angles in C(X)}, we define the angle between $\ell_{\p_1}, \ell_{\p_2} \in L(\x)$ by $|\theta_2 - \theta_1|$. For $\ell = \ell_{p,\theta}$, define $\theta(\ell) = \theta$. Define for an interval $J \subset B_1^1$ 
\begin{align}
    L_J(\x) = \{\ell \in L(\x) : \theta(\ell) \in J\}. 
\end{align}
The next lemma identifies how the curves spread in $L(\x)$. 
\begin{lemma}[Quantitative Broad-Narrow Reduction \protect{\cite[Lemma 1.27]{wangwurestriction}}]\label{lem: broad narrow}
    Let $(L,Y)_\delta$ be a collection of curves in $C(X)$ with their associated shading by $\delta$-cubes. Suppose that $\# L(\x) \geq \delta^{-\epsilon}$. 
    Then there is a scale $r = |\log \delta|^{-j} \in [\delta, 1]$ and an interval $J \subset B^1$ with $|J| \sim r$ such that the following hold. 
    \begin{itemize}
        \item $\# L_J(\x) \gtrsim 2^{-100 |\log \delta| / \log |\log \delta|} \#L(\x)$,
        \item There are 3 subsets $L_1(\x), L_2(\x), L_3(\x) \subset L_J(\x)$ which are $r/|\log \delta|$ $\theta$-separated and $\# L_j(\x) \gtrsim |\log \delta|^{-1} \# L_J(\x)$. 
    \end{itemize}
\end{lemma}

\begin{proof}
    The proof is the same as in \cite[Lemma 1.27]{wangwurestriction}, except we need to be more precise about losses due to pigeonholing. For example, a loss of the form $c_\epsilon \delta^\epsilon$ is unacceptable.
    
    Set $\rho_j = |\log \delta|^{-j}$ with $\rho_n \approx \delta$. Then $n \sim |\log \delta| / \log |\log \delta|$. 
    Run the following algorithm to locate the correct scale. 
    Partition $\omega_0 = B^1$ into $\sim |\log \delta|$ many finitely overlapping $\rho_1$-intervals $\omega_1 \in \Omega_1$. If there is an $\omega_1$ such that $\# L_{\omega_1}(\x) \geq 100^{-1} \# L(\x)$, repeat the procedure in this smaller interval $\omega_1$. We continue in this way until the first $j$ where $\# L_{\omega_{j+1}}(\x) < (100)^{-1} \# L_{\omega_j}(\x)$. In this case, there is a set $\Omega_{j+1}'$ of at least $50$ intervals $\omega_{j+1}$ such that
    \begin{align}
        \# L_{\omega_{j+1}}(\x) \geq (2 |\log \delta|)^{-1} \# L_{\omega_j}. 
    \end{align}
    We may extract 3 $\rho_{j+1}$-separated intervals from $\Omega_{j+1}'$ and consider the corresponding collection of $\rho_{j+1}$-$\theta$-separated curves $L_1(\x)$, $L_2(\x)$, $L_3(\x)$. Set $J = \omega_j$. Since $j \lesssim |\log \delta| / \log |\log \delta|$, we have $\# L_J(\x) \geq (100)^{-|\log \delta| / \log |\log \delta|} \# L(\x)$.
\end{proof}
Applying Lemma \ref{lem: broad narrow} to a point $\x \in E_L$, we may define 
\begin{align}
    r_L(\x) = r \text{ and } J_L(\x) = J. 
\end{align}
Both functions $r_L$ and $J_L$ are constant on $\delta$-cubes since the shadings are unions of grid-aligned $\delta$-cubes.

\subsection{Regular shadings}
It is convenient to work with regular shadings $Y(\ell)$.
\begin{definition}[Regular Shading \cite{KWZ}]
    A shading $Y(\ell)$ is regular if for each $\x \in Y(\ell)$ and each $\delta \leq r \leq 1$, 
    \begin{align}
        |Y(\ell) \cap B_r| \geq \frac{r |Y(\ell)|}{C|\log \delta|}.
    \end{align}
\end{definition}
The next Lemma says that every shading has a large regular refinement. 
\begin{lemma}[Regular refinements \protect{\cite[Lemma 2.7]{KWZ}}]\label{lem: regular refinement}
    A shading $Y(\ell)$ has a regular refinement $Y'(\ell) \subset Y(\ell)$ with $|Y'(\ell)| \geq \frac{1}{2} |Y(\ell)|$. 
\end{lemma}

\begin{proof}
    Straighten $\ell$ to a line, apply the result verbatim from \cite[Lemma 2.7]{KWZ}, and undo the straightening. 
\end{proof}

\subsection{From arbitrary measurable shadings to shadings by cubes}
We call $(L,Y)_\delta$ a \emph{measurable} $\lambda$-dense shading if for each $\ell \in L$, $Y(\ell) \subset N_{O(\delta)}(\ell)$ is a measurable set and $|Y(\ell)| \geq \lambda \delta^2$. Similarly define $E_{(L,Y)} = \bigcup_{\ell \in L} Y(\ell)$.  Measurable shadings will appear in the proof of Theorem \ref{thm: main phicurved kakeya}, and the following lemma tells us that we can replace them with our standard shadings by $\delta$-cubes.

\begin{lemma}[Reduction to shadings by cubes]\label{lem: meas to cubes}
    Let $(L,Y)_\delta$ be a set of curves in $\R^3$ with their associated $\lambda$-dense measurable shadings. There is an absolute constant $C_0$ such that, if 
    \begin{align}\label{eq: nontrivial density assumption}
        \frac{\lambda}{\delta} \geq C_0 \log(\#L + 2),
    \end{align}
    then there exists a $\Omega(\lambda)$-dense shading $(L, \tilde Y)$ by grid-aligned $\delta$-cubes so that 
    \begin{align}
        |E_{(L,\tilde Y)}| &\sim |E_{(L,Y)}|, \text{ and } \label{eq: meas to cubes: same vol} \\
        E_{(L, \tilde Y)} &\subset N_{O(\delta)}(E_{(L,Y)}).\label{eq: meas to cubes: containment}
    \end{align}
\end{lemma}

The proof is a standard randomization argument which we give below. 
\begin{proof}
    Let $\mathcal Q$ be a partition of $\R^3$ into axis-parallel $\delta$-cubes and set $\mathcal Q(\ell) := \{Q \in \mathcal Q : |Y(\ell) \cap Q| > 0\}$. Define the probabilities 
    \begin{align}
        p_Q :=\frac{|E_{(L,Y)} \cap Q|}{|Q|}
    \end{align}
    and consider a
    random set of cubes $ \tilde {\mathcal Q}$ where each $Q \in \mathcal Q$ is chosen independently with probability $p_Q$. Define the union of $\delta$-cubes
    $\tilde E := \bigcup_{Q \in \tilde {\mathcal Q}} Q$.
    Set 
    $\tilde {\mathcal Q}(\ell) := \tilde {\mathcal Q} \cap \mathcal Q(\ell)$ and 
    define the shadings by $\delta$-cubes $\tilde Y(\ell) := \bigcup_{Q \in \tilde {\mathcal Q}(\ell)} Q$.
    We have $E_{(L, \tilde Y)} = \tilde E$.

    Next define the random variables $X := \delta^{-3} |\tilde E|$ and for $\ell \in L$, $X_\ell := \delta^{-3} |\tilde Y(\ell)|$. We compute 
    \begin{align}
        \mathbb E X &= \delta^{-3} \sum_{Q \in \mathcal Q} |E \cap Q| = \delta^{-3} |E|, \\
        \mathbb E X_\ell &\geq \delta^{-3} \sum_{Q \in \mathcal Q(\ell)} |Y(\ell) \cap Q| = \delta^{-3} |Y(\ell)| \geq \frac{\lambda}{\delta}.
    \end{align}
    Also note that $E \supset Y(\ell)$ for any $\ell$, so $\mathbb E X \geq \lambda/\delta$. We would now like to show that with nonzero probability, $X \sim \delta^{-3} |E|$ and for all $\ell \in L$, $X_\ell \gtrsim \lambda /\delta$. 
    
    We use the following standard Chernoff estimate: if $X$ is a sum of independent Bernoulli random variables with mean $\mu$, then 
    \begin{align}
        \mathbb P(X < \mu/2) \leq e^{-\mu/8}, \qquad \mathbb P(X > 2\mu) \leq e^{-\mu/3}.
    \end{align}
    By applying these Chernoff estimates, the union bound, and the above computations for $\mathbb EX, \mathbb EX_\ell$, we find 
    \begin{align}
        \mathbb P(\delta^{-3}|E|/2 \leq X \leq 2\delta^{-3}|E| &\text{ and for all $\ell \in L$, } X_\ell \geq \lambda/(2\delta)) \\ &\geq 1 - (\#L +2)e^{-\lambda/(8\delta)} \\
        &\geq 1/2.
    \end{align}
    The last inequality follows from the assumption \eqref{eq: nontrivial density assumption} and choosing $C_0$ large enough. Fix such a realization $\tilde Q$ satisfying the properties. Then for each $\ell \in L$
    \begin{align}
        |\tilde Y(\ell)| = \delta^3 X_\ell \gtrsim \lambda \delta^2,
    \end{align}
    so $(L,\tilde Y)$ is $\Omega(\lambda)$-dense. Similarly, $|\tilde E| \sim |E|$. It is also clear from the construction that $\tilde E \subset N_{O(\delta)}(E_{(L,Y)})$. This finishes the proof of the lemma.  
\end{proof}

\section{Main Tools in the Proof of Theorem \ref{thm: CT 3 param kak}}\label{sec: 4}

Depending on whether the broad-narrow angle $r$ is large or small, we rely on different tools. 

\subsection{Large angle tools}
When $r$ is large, we rely on the following 3-linear Kakeya inequality due to Bourgain and Guth \cite{BourgainGuth}.
\begin{theorem}[Multilinear Kakeya \protect{\cite[Theorem 6] {BourgainGuth}}]\label{thm: multilinear Kakeya}
    Let $\mathbb T$ be a collection of $\delta$-neighborhoods of $C^\infty$-curves. Assume that for each $k \geq 0$, the $C^k$ norms of the curves are uniformly bounded. Then for any $\epsilon > 0$, 
    \begin{align}\label{eq: multilinear Kakeya}
        \int (\sum_{T_1,T_2,T_3 \in \mathbb T} \chi_{T_1} \chi_{T_2} \chi_{T_3} |\det(\bn_1(\x),  \bn_2(\x), \bn_3(\x))|)^{1/2} \lesssim_\epsilon \delta^{-\epsilon} \delta^3 (\# \mathbb T)^{3/2}. 
    \end{align}
\end{theorem}

\begin{remark}
    As originally stated, \cite[Theorem 6]{BourgainGuth} has no $\epsilon$-losses and applies only to algebraic curves of bounded degree. However this can be upgraded to smooth curves by using the Jackson type theorem \cite[Theorem 2]{jacksonapprox}, and tracking the degree $D$ in the proof. This is outlined in \cite{BourgainGuth}.
\end{remark}
The $\mathfrak c$-coniness condition is designed to work with Theorem \ref{thm: multilinear Kakeya}: 
\begin{lemma}[$\mathfrak c$-coniness means transverse]\label{lem: coniness means transverse}
    Suppose that $\mathcal C(X)$ satisfies the $\mathfrak c$-coniness condition. Let $\ell_{\p_1},\ell_{\p_2},\ell_{\p_3} \in \mathcal C(X)$ intersect at height $t$ and
    \begin{align}
        r/K \leq |\theta_i - \theta_j| \leq r
    \end{align}
    for all $i \neq j$. As long as $r \leq c(C_0) K^{-3} \mathfrak c$, then 
    \begin{align}
        |\det(\bn_{\p_1}(t), \bn_{\p_2}(t), \bn_{\p_3}(t))| \gtrsim \mathfrak c (r/K)^3. 
    \end{align}
\end{lemma}

\begin{proof}
    Define the curve $\bv(\theta) = (\omega_{\p_1,t}(\theta), 1)$. We have 
    \begin{align}
        \bn_{\p_i}(t) = \frac{\bv(\theta_i)}{|\bv(\theta_i)|} 
    \end{align}
    for $i = 1,2,3$. Condition \ref{cond: basic} guarantees $|\bv(\theta)| \lesssim 1$, so it is enough to show 
    \begin{align}\label{eq: 1.lem: coney means transverse}
        |\det(\bv(\theta_1), \bv(\theta_2), \bv(\theta_3))| \gtrsim  \mathfrak c (r/K)^3.
    \end{align}
    By Taylor expanding $\bv(\theta)$ in $\theta$, we get 
    \begin{align}
        \bv(\theta) = \bv(\theta_1) + (\theta - \theta_1) \dot \bv(\theta_1)  + \frac{(\theta - \theta_1)^2}{2}\ddot \bv(\theta_1) + O(r^3). 
    \end{align}
    Now compute the LHS of \eqref{eq: 1.lem: coney means transverse}:
    \begin{align}
        \det(\bv(\theta_1), \bv(\theta_2), \bv(\theta_3)) &= 
        \frac{(\theta_2 - \theta_1)(\theta_3 - \theta_1)(\theta_3 - \theta_2)}{2}\det(\bv(\theta_1), \dot \bv(\theta_1), \ddot \bv(\theta_1))
        + O(r^4).
    \end{align}
    By $\mathfrak c$-coniness, 
    \begin{align}
        |\det(\bv(\theta_1), \bv(\theta_2), \bv(\theta_3))| \gtrsim \mathfrak c (r/K)^3 
    \end{align}
    as long as $r \leq c K^{-3} \mathfrak c$, for a constant $c$ depending on $C_0$ in Condition \ref{cond: basic}. 
\end{proof}

\subsection{Small Angle Tools}

When $r$ is small, we run the argument in \cite{KWZ}. Though instead of using \cite[Theorem 1.7]{PYZ}, it is cleaner to use a version of \cite[Theorem 1.9]{zahlPlanar}. 

\begin{theorem}[Tweak of \protect{\cite[Theorem 1.9]{zahlPlanar}}]\label{thm: zahl planar tweak}
    Let $(C_i)_{i\geq 0}$ be constants, $I \subset \R$ an interval, and $\mathcal F \subset C^\infty(I)$ a family of curves satisfying 
    \begin{itemize}
        \item (Uniformly Smooth) For each $i \geq 0$, $\sup_{f \in \mathcal F} \| f^{(i)} \|_\infty \leq C_i$.
        \item  (Forbidding $2$nd order tangency) There exists a constant $c > 0$ such that for all $f,g \in \mathcal F$ we have
        \begin{align}\label{eq: forbid 2nd order tang}
            \inf_{t \in I} \sum_{i=0}^2 |f^{(i)}(t) - g^{(i)}(t)| \geq c \|f - g\|_{C^2(I)}. 
        \end{align}
    \end{itemize}
    Fix $K \geq 1$ and let $F \subset \mathcal F$ satisfy the one-dimensional ball condition 
    \begin{align}\label{eq: K 2dim ball}
        \# (F \cap B) \leq K(r/\delta) \text{ for all balls $B \subset C^2(I)$ of radius $r$}.
    \end{align}
    Let $f^\delta$ be the $\delta$-neighborhood of the graph of $f$, and let $Y(f) \subset f^\delta$ satisfy $|Y(f)| \geq \lambda \delta$. Then 
    \begin{align}\label{eq: thm zahl planar conclusion}
        \Big |\bigcup_{f \in F } Y(f)\Big | \geq c_\epsilon c^3 \delta^\epsilon \lambda^3 K^{-1} (\delta \# F). 
    \end{align}
    The constant $c_\epsilon$ depends on the numbers $\epsilon, |I|, C_0,\ldots C_{O_{\epsilon}(1)}$.
\end{theorem}

\begin{proof}[Deducing Theorem \ref{thm: zahl planar tweak} from the proof of \protect{\cite[Theorem 1.9]{zahlPlanar}}]
    Only a small change needs to be made on page 34 of \cite{zahlPlanar}. At the top of page 34, the author concludes that: ``$\ell \lesssim c^{-O(1)}$''. In the case of $k = 2$, forbidding 2nd order tangency, this is $\ell \lesssim c^{-3}$.
    In equation (5.9) at the bottom of the same page, one needs to leave the factor $c^{-3}$ explicit at the end of the calculation.
    
    The conclusion of Theorem 1.9 from \cite{zahlPlanar} is written as an upper bound on an $L^{3/2}$ norm, though \eqref{eq: thm zahl planar conclusion} is equivalent by standard arguments. The condition \eqref{eq: K 2dim ball} is also different from \cite{zahlPlanar}. The case for general $K$ is obtained by randomly sampling each curve $f \in F$ with probability $K^{-1}$. This adjustment was needed in \cite{KWZ} too, where such a comment can also be found. 
\end{proof}

The $\tau$-twistiness condition is, of course, designed to work with Theorem \ref{thm: zahl planar tweak}:

\begin{lemma}\label{lem: twist implies 2nd ord tangency}
    Suppose that $X$ is $\tau$-twisty and satisfies the basic conditions.
    \begin{enumerate}
        \item Then for each $\p \in B_1^3$, the family $\mathcal F(\p, X)$ forbids 2nd-order tangency in the sense of \eqref{eq: forbid 2nd order tang} with $c = \Omega(\tau)$.
        \item Let $A \subset B(\p, \delta^{1/2})$ be a collection of points satisfying the one-dimensional ball condition 
        \begin{align}\label{eq: 1 twist implies 2nd ord tangency}
            \#(A \cap B_r) \leq K(r/\delta) \text{ for all balls $B_r \subset \R^3$}.
        \end{align}
        Then the set $F = \{f_{\p'} : \p' \in A\} \subset \mathcal F(\p,X)$ satisfies the one-dimensional ball condition 
        \begin{align}
            \#(F \cap B_r) \leq (\tau^{-3}K) (r/\delta) \text{ for all balls $B_r \subset C^2(B_1^1)$.}
        \end{align}
    \end{enumerate}
\end{lemma}

\begin{proof}
    To show $\mathcal F(\p,X)$ forbids 2nd-order tangency, it is enough to show 
    \begin{align}
        \inf_t |F_{\p,t}(\p_1') - F_{\p,t}(\p_2')| \gtrsim \tau \cdot \sup_t |F_{\p,t}(\p_1') - F_{\p,t}(\p_2')|
    \end{align}
    for any $\p_1',\p_2'$. Recall that $F_{\p,t}(\p') = M_\p(t) \cdot (\p'-\p)$. 
    By the basic conditions and $\tau$-twistiness, the singular values of $M_{\p}(t)$ are $\lesssim 1$ and $\gtrsim \tau$. Thus 
    \begin{align}\label{eq: twist implies 2nd ord tangency}
        \tau |\p_2' - \p_1'| \lesssim |M_\p(t) (\p_2'- \p_1')| \lesssim |\p_2' - \p_1'|,
    \end{align}
    which proves claim (1). 
    For claim (2), the first inequality in \eqref{eq: twist implies 2nd ord tangency} shows that 
    \begin{align}
        \tau |\p_2' - \p_1'| \lesssim \|f_{\p'_2} -f_{\p'_1}\|_{C^2(B_1^1)}.
    \end{align}
    Thus a $C^2$-ball of radius $r$ has parameters inside $B_{O(r/\tau)} \subset \R^3$. By covering $B_{O(r/\tau)}$ by approximately $\tau^{-3}$ balls $B_i$ of radius $r$ and using \eqref{eq: 1 twist implies 2nd ord tangency}, 
    \begin{align}
        \#(F \cap B_r) &\lesssim \sum_i \#(A \cap B_i) \\
        &\lesssim (\tau^{-3} K) (r/\delta).
    \end{align}
\end{proof}

\section{Grain Structure of $\cC(X)$}\label{sec: 5}

We will build up the structural properties of the $\delta$-tubes in a fatter $\rho$-tube and eventually apply Theorem \ref{thm: zahl planar tweak}. Most of the steps here have analogues in \cite{KWZ}, and we indicate when this is the case.   
Fix a curve $\ell_\p \in \mathcal C(X)$ and define
\begin{align}
    \mathcal C(\p, X) = \{\ell' \in \mathcal C(X) : \ell' \text{ intersects } \ell_p \}.
\end{align}
Recall from  Definition \ref{def: spread curves} that $(\hat p_{\p, s}(\theta'), \theta')$ is the parameter of a curves passing through $\ell_{\p}$ at height $s$ and angle parameter $\theta'$. We may therefore parameterize $\mathcal C(\p, X)$ by $s$ and $\theta'$, and write $\ell_{s,\theta'}$ as a slight abuse of notation. We first apply Taylor's theorem to isolate the leading order terms of $\ell_{s,\theta'}(t)$ that we can analyze. 

\begin{lemma}[Expansion of $\ell_{s,\theta'}(t)$]\label{lem: taylor expand C(X,p)}
    Let $X$ satisfy the basic conditions. The curve $\ell_{s,\theta'} \in \mathcal C(X,\p)$ has the expansion 
    \begin{align} 
        \ell_{s, \theta'}(t) - \ell_{\p}(t) = (&(t-s)(\theta'-\theta)\gamma_\p(s_0) \\&+ O(|\theta'-\theta|^2 |t - s| + |\theta'-\theta| |t - s|^2 +|\theta'-\theta| |t -s| |s - s_0|),0). \nonumber 
    \end{align}
\end{lemma}

\begin{remark}
    If $X$ is in addition planey, as is the case for $\mathcal L_{\rm{SL}_2}$, we get a stronger version of Lemma \ref{lem: taylor expand C(X,p)}. This is Lemma \ref{lem: planey taylor expand C(X,p)} in Appendix \ref{appendix: A}. 
\end{remark}

\begin{proof}
    Define $Y(\theta', s, t) = X(\hat p_{\p,s}(\theta'), \theta', t) - X(\p,t)$. Then $Y$ has the following properties:
    \begin{enumerate}
        \item $Y(\theta, s, t) = 0$ for all $s,t \in B_1^1$. 
        \item $Y(\theta', s, s) = 0$ for all $\theta', s \in B_1^1$. 
    \end{enumerate}
    Taylor expand $Y$ in $t$ around $s$ to get 
    \begin{align}
        Y(\theta', s, t) = \partial_t Y(\theta', s, s)(t - s) + R(\theta', s, t),
    \end{align}
    where $R(\theta', s, t) = \int_{s}^t \partial_t^2 Y(\theta', s, a) (t-a)da$. Expanding $\partial_t^2 Y(\theta',s, a)$ to constant order in $\theta'$ near $\theta$, 
    \begin{align}
        R(\theta', s, t) &=  \int_s^t O(|\theta'-\theta|) (t - a) da\\
        &= O(|\theta' - \theta| |t - s|^2).
    \end{align}
    Taylor expand $\partial_t Y(\theta', s, s)$ in $\theta'$ around $\theta$ to get 
    \begin{align}\label{eq: splice to planey version}
        Y(\theta', s, t) = \partial_\theta \partial_t Y(\theta, s, s) (t - s)(\theta'-\theta) + R_2(\theta', s)(t-s) + O(|\theta' - \theta| |t-s|^2), 
    \end{align}
    where $R_2(\theta', s) =  \int_{\theta}^{\theta'} \partial_\theta^2 \partial_t Y(a, s, s)(\theta' - a) da$. 
   
    By the basic conditions we have 
    \begin{align}
        R_2(\theta', s) = O(|\theta' - \theta|^2),
    \end{align}
    so 
    \begin{align}
        Y(\theta', s, t) = \gamma_\p(s)(t-s)(\theta' - \theta) + O(|\theta' - \theta|^2 |t-s| + |\theta' - \theta| |t - s|^2).
    \end{align} 
    Expanding $\gamma_\p(s) = \gamma_\p(s_0) + O(|s -s_0|)$ we get 
    \begin{align}
        Y(\theta', s, t) =& \gamma_p(s_0) (t - s) (\theta' - \theta) \\&+ O(|\theta' - \theta |^2 |t-s| + |\theta' - \theta| |t - s|^2 + |\theta' - \theta| |t -s| |s - s_0|). \nonumber
    \end{align}
\end{proof}

We now use Lemma \ref{lem: taylor expand C(X,p)} to give the grain structure. Fix $r \geq \delta$, 
\begin{align}
    \rho \in [\delta, \min((\delta r)^{1/2}, \delta / r)],
\end{align}
and assume that $|\theta' - \theta| \sim r$. As a consequence of Lemma \ref{lem: taylor expand C(X,p)}, 
\begin{align}
    |\ell_{s,\theta'}(t) - \ell_\p(t)| \sim |t-s| |\theta' - \theta| \sim r|t-s|
\end{align}
as long as $|\theta'-\theta|,|s-s_0|,|t-s| \ll 1$. 
So $\ell_{s,\theta'}$ is inside $N_\rho(\ell_\p)$ for only $t$ in the range 
\begin{align}
    |t - s| \lesssim \rho /r.
\end{align}
Fix $s_0 \in B_1^1$ and consider $|s - s_0| \leq (\delta / r)^{1/2}$. For this range of $s$ and inside $N_\rho(\ell_{\p})$, we have by Lemma \ref{lem: taylor expand C(X,p)} that 
\begin{align}\label{eq: nearby curves delta}
    \ell_{s, \theta}(t) - \ell_p(t) = ((t - s) (\theta' - \theta) \gamma_p(s_0) + O(\delta), t). 
\end{align}
The image of the RHS of \eqref{eq: nearby curves delta} for the range $|t -s| \leq \rho / r$ is contained in a $\delta \times \rho  \times (\delta /r)^{1/2}$-slab. After adding $\ell_\p(t)$, the slab is sheared slightly. These slabs are much more conveniently describe as preimages of certain rectangles under the twisted projection $\pi_\p : N_\rho(\ell_\p) \to N_\rho(B_1^1)$ defined in Section \ref{sec: 2}. 
\begin{definition}[Prism decomposition]\label{def: prism decomposition}
    Tile $B^1 \times [-\rho, \rho]$ by approximately $(r / \delta)^{1/2}\cdot (\rho / \delta)$ translates of $[-(\delta/ r)^{1/2},(\delta /r)^{1/2}]\times [-\delta, \delta]$, and call the family $\mathcal R$. The \emph{prism decomposition} $\mathcal P_\rho(\ell_\p)$ of $N_\rho(\ell_\p)$ is 
    \begin{align}
        \mathcal P_\rho(\ell_\p) = \{P = \pi_p^{-1}(R) : R \in \mathcal R\}. 
    \end{align}
    Each $P \in \mathcal P_\rho(\ell_\p)$ has dimensions $\delta \times \rho \times d$, where $d = (\delta/ r)^{1/2}$. We denote the sub-collection generated by the rectangles intersecting the $x$-axis by $\mathcal P_\rho^0(\ell_\p)$. 
\end{definition}

Cut $N_\rho(\ell)$ into segments $U \in \mathcal U$ of length $(\delta /r)^{1/2}$.
Equation \eqref{eq: nearby curves delta} shows that for each $U \in \mathcal U$, 
\begin{align}
    U \cap (\bigcup_{\ell_{s,\theta} \in \mathcal C(X,p), |\theta| \sim r} \ell_{s,\theta}) \subset P
\end{align}
for exactly one $P \in \mathcal P^0(T_\rho)$. 
The next lemma says that the prism decompositions of nearby curves agree on their overlap. 
\begin{lemma}[Compatibility of Prism Decompositions]\label{lem: prism compatibility}
    There are constants $c >0$ and $C = O(1)$ so that the following holds. Fix a curve $\ell_\p \in \cC(X)$ and let 
    $\mathcal P = \mathcal P_\rho(\ell_\p)$ be its prism decomposition. Consider $\ell_{\p'} \in \cC(X)$ 
    satisfying $d(\ell_\p, \ell_{\p'}) \leq c\rho$, and fix $P' = \pi^{-1}(R') \in \mathcal P^0_\rho(\ell_{\p'})$, where
    \begin{align}
        R' = ([-d, d] \times [-\delta, \delta]) + (t_0,0). 
    \end{align}
    Then there is a $y = O(\rho) \in \R$ such that 
    \begin{align}
        C^{-1} \pi_\p^{-1}(R' + (0,y)) \subset \pi_{\p'}^{-1}(R') \subset C \pi_\p^{-1}(R' + (0,y)). 
    \end{align}
    Thus the prism decompositions of $\ell_\p$ and $\ell_{\p'}$ agree on their overlap. 
\end{lemma}

\begin{proof}
    Write 
    \begin{align}
        q_\p(x,t) = (x - X(\p,t)) \cdot \gamma_\p(t)^\perp,
    \end{align}
    so that $\pi_\p(x,t)=(t,q_\p(x,t))$, and set 
    \begin{align}
        y=(X(\p',t_0)-X(\p,t_0))\cdot \gamma_\p(t_0)^\perp.
    \end{align}
    Then $y = O(\rho)$.
    Take $(x,t) \in P' = \pi_{\p'}^{-1}(R')$. In particular, $|t-t_0| \leq d$, $|q_{\p'}(x,t)|\leq \delta$, and $|x-X(\p',t)| \lesssim \rho$. We have 
    \begin{align}\label{eq: q expr}
        q_\p(x,t) - q_{\p'}(x,t) = (&X(\p',t)-X(\p,t)) \cdot \gamma_\p(t)^\perp \\
        &+ (x-X(\p',t)) \cdot(\gamma_\p(t)^\perp-\gamma_{\p'}(t)^\perp).
    \end{align}
    For $F(t)  = (X(\p',t)-X(\p,t)) \cdot \gamma_\p(t)^\perp$, the basic conditions and Taylor's theorem give 
    \begin{align}
        F(t) = F(t_0) + O(|\p'-\p| |t-t_0|) = y + O(c\rho d)/
    \end{align}
    The second term in \eqref{eq: q expr} is $O(\rho|\p'-\p|)= O(c \rho^2)$. Thus 
    \begin{align}
        q_\p(x,t) = q_{\p'}(x,t) + y + O(c(\rho d + \rho^2)) = q_{\p'}(x,t) + y + O(c\delta),
    \end{align}
    since $\rho d \leq \delta$ and $\rho^2 \leq \delta$. Then with $c$ sufficiently small, 
    \begin{align}
        P' \subset C \pi_{\p}^{-1}(R' + (0,y)).
    \end{align}
    A similar calculation gives the other inclusion. 
\end{proof}

\begin{figure}
    \centering

\tikzset{every picture/.style={line width=0.75pt}} 

\begin{tikzpicture}[x=0.75pt,y=0.75pt,yscale=-1,xscale=1]

\draw (310,235) node  {\includegraphics[width=255pt,height=135pt]{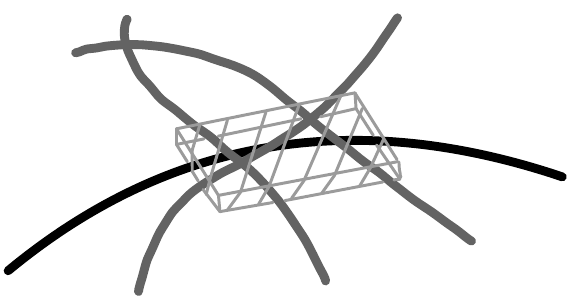}};
\draw    (272,277) -- (377,258) ;
\draw [shift={(377,258)}, rotate = 169.74] [color={rgb, 255:red, 0; green, 0; blue, 0 }  ][line width=0.75]    (0,5.59) -- (0,-5.59)(10.93,-3.29) .. controls (6.95,-1.4) and (3.31,-0.3) .. (0,0) .. controls (3.31,0.3) and (6.95,1.4) .. (10.93,3.29)   ;
\draw [shift={(272,277)}, rotate = 349.74] [color={rgb, 255:red, 0; green, 0; blue, 0 }  ][line width=0.75]    (0,5.59) -- (0,-5.59)(10.93,-3.29) .. controls (6.95,-1.4) and (3.31,-0.3) .. (0,0) .. controls (3.31,0.3) and (6.95,1.4) .. (10.93,3.29)   ;
\draw    (355,197) -- (380,238) ;
\draw [shift={(380,238)}, rotate = 238.63] [color={rgb, 255:red, 0; green, 0; blue, 0 }  ][line width=0.75]    (0,5.59) -- (0,-5.59)(10.93,-3.29) .. controls (6.95,-1.4) and (3.31,-0.3) .. (0,0) .. controls (3.31,0.3) and (6.95,1.4) .. (10.93,3.29)   ;
\draw [shift={(355,197)}, rotate = 58.63] [color={rgb, 255:red, 0; green, 0; blue, 0 }  ][line width=0.75]    (0,5.59) -- (0,-5.59)(10.93,-3.29) .. controls (6.95,-1.4) and (3.31,-0.3) .. (0,0) .. controls (3.31,0.3) and (6.95,1.4) .. (10.93,3.29)   ;
\draw    (379,241) -- (379,252) ;
\draw [shift={(379,252)}, rotate = 270] [color={rgb, 255:red, 0; green, 0; blue, 0 }  ][line width=0.75]    (0,2.24) -- (0,-2.24)(4.37,-1.32) .. controls (2.78,-0.56) and (1.32,-0.12) .. (0,0) .. controls (1.32,0.12) and (2.78,0.56) .. (4.37,1.32)   ;
\draw [shift={(379,241)}, rotate = 90] [color={rgb, 255:red, 0; green, 0; blue, 0 }  ][line width=0.75]    (0,2.24) -- (0,-2.24)(4.37,-1.32) .. controls (2.78,-0.56) and (1.32,-0.12) .. (0,0) .. controls (1.32,0.12) and (2.78,0.56) .. (4.37,1.32)   ;
\draw    (321,166) -- (321.94,197) ;
\draw [shift={(322,199)}, rotate = 268.26] [color={rgb, 255:red, 0; green, 0; blue, 0 }  ][line width=0.75]    (10.93,-3.29) .. controls (6.95,-1.4) and (3.31,-0.3) .. (0,0) .. controls (3.31,0.3) and (6.95,1.4) .. (10.93,3.29)   ;
\draw    (171,235) -- (175.74,271.02) ;
\draw [shift={(176,273)}, rotate = 262.5] [color={rgb, 255:red, 0; green, 0; blue, 0 }  ][line width=0.75]    (10.93,-3.29) .. controls (6.95,-1.4) and (3.31,-0.3) .. (0,0) .. controls (3.31,0.3) and (6.95,1.4) .. (10.93,3.29)   ;
\draw    (447,280) -- (422,280) ;
\draw [shift={(420,280)}, rotate = 360] [color={rgb, 255:red, 0; green, 0; blue, 0 }  ][line width=0.75]    (10.93,-3.29) .. controls (6.95,-1.4) and (3.31,-0.3) .. (0,0) .. controls (3.31,0.3) and (6.95,1.4) .. (10.93,3.29)   ;
\draw    (269.38,236) -- (243.38,217) ;
\draw    (236,242) -- (269.38,236) ;
\draw  [draw opacity=0] (243.81,239.76) .. controls (244.11,233.4) and (246.42,227.5) .. (250.15,222.74) -- (273.78,241.22) -- cycle ; \draw   (243.81,239.76) .. controls (244.11,233.4) and (246.42,227.5) .. (250.15,222.74) ;  

\draw (317,147.4) node [anchor=north west][inner sep=0.75pt]    {$P$};
\draw (157,212.4) node [anchor=north west][inner sep=0.75pt]    {$Y( \ell )$};
\draw (449,272.4) node [anchor=north west][inner sep=0.75pt]    {$Y( \ell ')$};
\draw (331,269.4) node [anchor=north west][inner sep=0.75pt]    {$d$};
\draw (369,199.4) node [anchor=north west][inner sep=0.75pt]    {$\rho $};
\draw (385,239.4) node [anchor=north west][inner sep=0.75pt]    {$\delta $};
\draw (209,220.4) node [anchor=north west][inner sep=0.75pt]    {$\sim r$};

\end{tikzpicture}

    \caption{Prism $P$ of dimensions $\delta \times \rho \times d$ around $\ell$. The gray shadings $Y(\ell')$ intersect $Y(\ell)$ at angle $\sim r$, and fill out $P$.}
    \label{fig: prisms filled up}
\end{figure}

The next lemma says that a $\delta$-tube picks up the mass of the prism that contains it. This is analogous to \cite[Lemma 2.10]{KWZ}.
\begin{lemma}\label{lem: prism filled up}
    Let $0 < \delta \ll 1$, $r \in (\delta, 1]$, $\rho \in [\delta, \min((\delta r)^{1/2}, \delta / r)]$. Let $(L,Y)_{\delta}$ be a set of curves with regular, $\lambda$-dense shading. Fix $\ell \in L$ and suppose that $P \in \mathcal P_\rho^0(\ell)$ satisfies 
    \begin{align}\label{eq: prism filled up hypothesis}
        |P \cap Y(\ell) \cap \{\x : r_L(\x) \sim r, \theta(\ell) \in J_L(\x)\}| \geq \alpha d \delta^2. 
    \end{align}
    There is a constant $C = O(1)$ so that 
    \begin{align}
        |CP \cap E_L| \gtrapprox \alpha \lambda^2 |P|. 
    \end{align}
\end{lemma}
See Figure \ref{fig: prisms filled up} for the geometric picture of Lemma \ref{lem: prism filled up}.
We use the same notation as \cite{KWZ} and most of the proof reads exactly the same, though we include it here since the particular grain geometry we developed above features prominently. 
\begin{proof}
    Let $X$ be the set on the LHS of \eqref{eq: prism filled up hypothesis}. Let $\x_1,\ldots, \x_M$, $M \gtrsim \alpha d r \delta^{-1}$ be points in $X$ arranged so that $|\x_i - \x_j| \geq 100 |i - j| \delta /r$. For each index $i$, we can find a curve $\ell_i$ with $|\theta(\ell_i) - \theta(\ell)| \in [r |\log \delta|^{-1}, 10r]$. We see this for fixed $\x = \x_i$ as follows. There is an interval $J = J_L(\x)$ of length $\sim r$ with two $r/|\log \delta|$ $\theta$-separated (non-empty) sets of curves $L_1(\x), L_2(\x) \subset L_J(\x)$. By hypothesis $\theta(\ell) \in J_L(\x)$, so we can choose $\ell'$ in at least one of $L_1(\x)$ or $L_2(\x)$ with $|\theta(\ell') - \theta(\ell)|$ of the claimed size.

    Recall that, after a change of variables, the curves in $P$ are parameterized by 
    \begin{align}
        \ell_i(t) = (\gamma(t_0) \theta (t - t_i), t) + O(\delta),
    \end{align}
    where $\theta \approx r$. If $|t - t_i| \ll \rho /r$ then $\ell_{i}(t) \in P$. So $\ell_i \cap CP$ is up to an $O(\delta)$ error a line segment of length $\rho / r$, for an appropriate choice of $C$. Set $\tilde P = CP$. Since $Y(\ell_i)$ is a regular shading, 
    \begin{align}
        |Y(\ell_i) \cap \tilde P| \gtrsim \lambda (\rho /r) \delta^2.
    \end{align}
    So we can estimate the $L^1$ norm 
    \begin{align}
        \Big \| \sum_{i=1}^M \chi_{Y(\ell_i)}\Big \|_{L^1(\tilde P)} \gtrsim \lambda (\rho / r) \delta^2 M.
    \end{align}
    If $i \neq j$ and $Y(\ell_i) \cap Y(\ell_j) \neq \emptyset$, the law of sines applied to the triangle formed by $\ell, \ell_i$, and $\ell_j$ says that $\angle(\ell_i,\ell_j) \gtrsim \delta r \rho^{-1} |i - j|$. Hence for any $1 \leq i,j \leq M$,
    \begin{align}
        |\tilde P \cap Y(\ell_i) \cap Y(\ell_j)| \leq |P \cap N_\delta(\ell_i) \cap N_\delta(\ell_j)| \lesssim \frac{\delta^2 \rho}{r(|i - j| + 1)}.
    \end{align}
    This lets us estimate the $L^2$ norm
    \begin{align}
        \Big \| \sum_{i=1}^M \chi_{Y(\ell_i)}\Big \|_{L^2(\tilde P)}^2 \lesssim \delta^2 \rho r^{-1} \sum_{1 \leq i,j \leq M} \frac{1}{|i -j| + 1} \lessapprox \delta^2 \rho r^{-1} M. 
    \end{align}
    By Cauchy-Schwarz, 
    \begin{align}
        |\tilde P \cap E_L| &\geq 
        \frac{\|\sum_{i=1}^M \chi_{Y(\ell_i)}\|_{L^1(\tilde P)}^2}{\| \sum_{i=1}^M \chi_{Y(\ell_i)}\|_{L^2(\tilde P)}^2} \\
        &\gtrapprox \alpha \lambda^2 |P|.
    \end{align}
    
\end{proof}

In the next Lemma, we consider the curves in the $\rho$-neighborhood of a fixed curve $\ell_0$. Each curve accumulates a prism by Lemma \ref{lem: prism filled up} and contributes to the total mass of $E_L$ inside $N_\rho(\ell_0)$. This is analogous to Lemma \cite[Lemma 2.12]{KWZ}. 

\begin{lemma}\label{lem: add prisms fat tube}
    Let $0 < \delta \ll 1$, $r \in (\delta, 1]$, $\rho \in (\delta, \min((\delta r)^{1/2}, \delta / r)]$. Let $(L,Y)$ and $(L',Y')$ be sets of curves with regular, $\lambda$-dense shadings. Suppose $L' \subset L \cap B(\ell_0, \rho)$. Also assume that 
    \begin{align}
        Y'(\ell') \subset \{\x : r_L(\x) \sim r, \theta(\ell') \in J_L(\x)\}
    \end{align}
    for all $\ell' \in L'$. Then there exists a curve $\tilde \ell \in B(\ell_0, 2\rho)$ and a $\Omega^*(\lambda^4)$-dense shading $\tilde Y$ of $\tilde \ell$ by $\rho$-cubes. After replacing $(L',Y')$ by a $\Omega^*(1)$-density refinement, we have
    \begin{align}\label{eq: add prisms fat tube}
        |E_L \cap Q| \gtrapprox \lambda^3 \rho^2 \Big |\bigcup_{\ell' \in L'} \pi_{\ell_0}(Y'(\ell'))\Big| \text{ For each $\rho$-cube $Q \subset \tilde Y(\tilde \ell)$.}
    \end{align}
\end{lemma}

\begin{figure}
    \centering

\tikzset{every picture/.style={line width=0.75pt}} 

\begin{tikzpicture}[x=0.75pt,y=0.75pt,yscale=-1,xscale=1]

\draw (337,199) node  {\includegraphics[width=326pt,height=204pt]{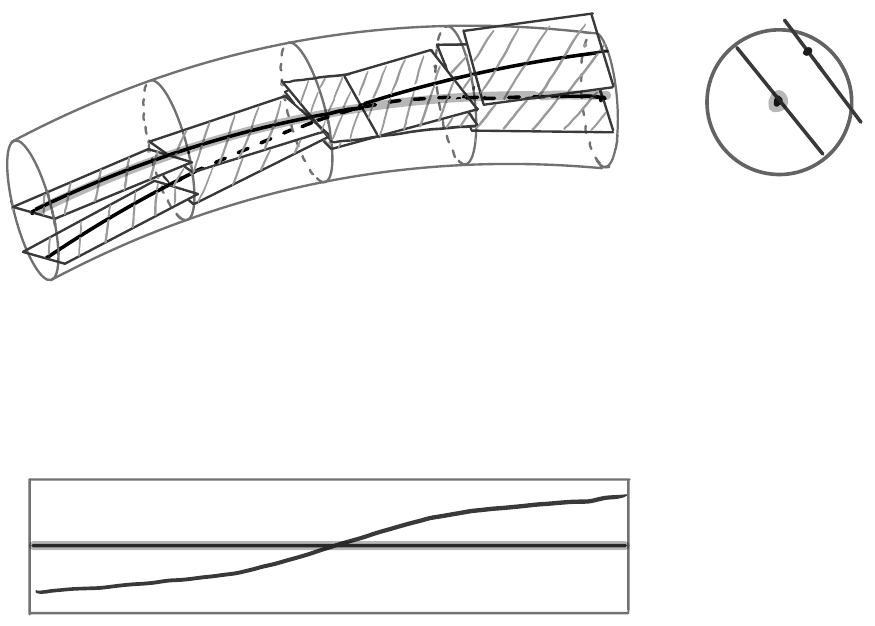}};
\draw    (305,343) -- (431,343) ;
\draw [shift={(431,343)}, rotate = 180] [color={rgb, 255:red, 0; green, 0; blue, 0 }  ][line width=0.75]    (0,5.59) -- (0,-5.59)(10.93,-3.29) .. controls (6.95,-1.4) and (3.31,-0.3) .. (0,0) .. controls (3.31,0.3) and (6.95,1.4) .. (10.93,3.29)   ;
\draw    (266,343) -- (134,343) ;
\draw [shift={(134,343)}, rotate = 360] [color={rgb, 255:red, 0; green, 0; blue, 0 }  ][line width=0.75]    (0,5.59) -- (0,-5.59)(10.93,-3.29) .. controls (6.95,-1.4) and (3.31,-0.3) .. (0,0) .. controls (3.31,0.3) and (6.95,1.4) .. (10.93,3.29)   ;
\draw    (118,291) -- (118.2,271) ;
\draw [shift={(118.2,271)}, rotate = 90.57] [color={rgb, 255:red, 0; green, 0; blue, 0 }  ][line width=0.75]    (0,5.59) -- (0,-5.59)(10.93,-3.29) .. controls (6.95,-1.4) and (3.31,-0.3) .. (0,0) .. controls (3.31,0.3) and (6.95,1.4) .. (10.93,3.29)   ;
\draw    (118,310) -- (118,329) ;
\draw [shift={(118,329)}, rotate = 270] [color={rgb, 255:red, 0; green, 0; blue, 0 }  ][line width=0.75]    (0,5.59) -- (0,-5.59)(10.93,-3.29) .. controls (6.95,-1.4) and (3.31,-0.3) .. (0,0) .. controls (3.31,0.3) and (6.95,1.4) .. (10.93,3.29)   ;
\draw    (482,79) -- (466.88,60.74) ;
\draw [shift={(465.6,59.2)}, rotate = 50.37] [fill={rgb, 255:red, 0; green, 0; blue, 0 }  ][line width=0.08]  [draw opacity=0] (12,-3) -- (0,0) -- (12,3) -- cycle    ;
\draw    (286,176) -- (286.96,251) ;
\draw [shift={(287,254)}, rotate = 269.27] [fill={rgb, 255:red, 0; green, 0; blue, 0 }  ][line width=0.08]  [draw opacity=0] (8.93,-4.29) -- (0,0) -- (8.93,4.29) -- cycle    ;
\draw    (306.38,152) -- (281,152) ;
\draw [shift={(281,152)}, rotate = 360] [color={rgb, 255:red, 0; green, 0; blue, 0 }  ][line width=0.75]    (0,5.59) -- (0,-5.59)(10.93,-3.29) .. controls (6.95,-1.4) and (3.31,-0.3) .. (0,0) .. controls (3.31,0.3) and (6.95,1.4) .. (10.93,3.29)   ;
\draw    (323,152) -- (348.38,152.5) ;
\draw [shift={(348.38,152.5)}, rotate = 181.13] [color={rgb, 255:red, 0; green, 0; blue, 0 }  ][line width=0.75]    (0,5.59) -- (0,-5.59)(10.93,-3.29) .. controls (6.95,-1.4) and (3.31,-0.3) .. (0,0) .. controls (3.31,0.3) and (6.95,1.4) .. (10.93,3.29)   ;
\draw    (150,92) -- (161.46,133.07) ;
\draw [shift={(162,135)}, rotate = 254.41] [color={rgb, 255:red, 0; green, 0; blue, 0 }  ][line width=0.75]    (10.93,-3.29) .. controls (6.95,-1.4) and (3.31,-0.3) .. (0,0) .. controls (3.31,0.3) and (6.95,1.4) .. (10.93,3.29)   ;
\draw    (280,36) -- (280.95,74) ;
\draw [shift={(281,76)}, rotate = 268.57] [color={rgb, 255:red, 0; green, 0; blue, 0 }  ][line width=0.75]    (10.93,-3.29) .. controls (6.95,-1.4) and (3.31,-0.3) .. (0,0) .. controls (3.31,0.3) and (6.95,1.4) .. (10.93,3.29)   ;
\draw    (149,199) -- (149.93,174) ;
\draw [shift={(150,172)}, rotate = 92.12] [color={rgb, 255:red, 0; green, 0; blue, 0 }  ][line width=0.75]    (10.93,-3.29) .. controls (6.95,-1.4) and (3.31,-0.3) .. (0,0) .. controls (3.31,0.3) and (6.95,1.4) .. (10.93,3.29)   ;
\draw    (101,143) -- (101.2,123) ;
\draw [shift={(101.2,123)}, rotate = 90.57] [color={rgb, 255:red, 0; green, 0; blue, 0 }  ][line width=0.75]    (0,5.59) -- (0,-5.59)(10.93,-3.29) .. controls (6.95,-1.4) and (3.31,-0.3) .. (0,0) .. controls (3.31,0.3) and (6.95,1.4) .. (10.93,3.29)   ;
\draw    (101,166) -- (101,187) ;
\draw [shift={(101,187)}, rotate = 270] [color={rgb, 255:red, 0; green, 0; blue, 0 }  ][line width=0.75]    (0,5.59) -- (0,-5.59)(10.93,-3.29) .. controls (6.95,-1.4) and (3.31,-0.3) .. (0,0) .. controls (3.31,0.3) and (6.95,1.4) .. (10.93,3.29)   ;
\draw    (120,102) -- (134.44,151.08) ;
\draw [shift={(135,153)}, rotate = 253.61] [color={rgb, 255:red, 0; green, 0; blue, 0 }  ][line width=0.75]    (10.93,-3.29) .. controls (6.95,-1.4) and (3.31,-0.3) .. (0,0) .. controls (3.31,0.3) and (6.95,1.4) .. (10.93,3.29)   ;
\draw    (427,136) .. controls (445.43,153.52) and (464.05,155.79) .. (485.67,139.99) ;
\draw [shift={(487,139)}, rotate = 142.75] [color={rgb, 255:red, 0; green, 0; blue, 0 }  ][line width=0.75]    (10.93,-3.29) .. controls (6.95,-1.4) and (3.31,-0.3) .. (0,0) .. controls (3.31,0.3) and (6.95,1.4) .. (10.93,3.29)   ;
\draw    (528,140) -- (535,132) ;
\draw [shift={(535,132)}, rotate = 131.19] [color={rgb, 255:red, 0; green, 0; blue, 0 }  ][line width=0.75]    (0,2.24) -- (0,-2.24)(4.37,-1.32) .. controls (2.78,-0.56) and (1.32,-0.12) .. (0,0) .. controls (1.32,0.12) and (2.78,0.56) .. (4.37,1.32)   ;
\draw [shift={(528,140)}, rotate = 311.19] [color={rgb, 255:red, 0; green, 0; blue, 0 }  ][line width=0.75]    (0,2.24) -- (0,-2.24)(4.37,-1.32) .. controls (2.78,-0.56) and (1.32,-0.12) .. (0,0) .. controls (1.32,0.12) and (2.78,0.56) .. (4.37,1.32)   ;
\draw    (386,241) -- (381.26,277.02) ;
\draw [shift={(381,279)}, rotate = 277.5] [color={rgb, 255:red, 0; green, 0; blue, 0 }  ][line width=0.75]    (10.93,-3.29) .. controls (6.95,-1.4) and (3.31,-0.3) .. (0,0) .. controls (3.31,0.3) and (6.95,1.4) .. (10.93,3.29)   ;
\draw    (443,334) -- (416.39,306.44) ;
\draw [shift={(415,305)}, rotate = 46.01] [color={rgb, 255:red, 0; green, 0; blue, 0 }  ][line width=0.75]    (10.93,-3.29) .. controls (6.95,-1.4) and (3.31,-0.3) .. (0,0) .. controls (3.31,0.3) and (6.95,1.4) .. (10.93,3.29)   ;
\draw    (180,198) -- (180,168) ;
\draw [shift={(180,166)}, rotate = 90] [color={rgb, 255:red, 0; green, 0; blue, 0 }  ][line width=0.75]    (10.93,-3.29) .. controls (6.95,-1.4) and (3.31,-0.3) .. (0,0) .. controls (3.31,0.3) and (6.95,1.4) .. (10.93,3.29)   ;

\draw (279,334.4) node [anchor=north west][inner sep=0.75pt]    {$1$};
\draw (112,294) node [anchor=north west][inner sep=0.75pt]    {$\rho $};
\draw (438,34.4) node [anchor=north west][inner sep=0.75pt]    {$\gamma _{\mathbf{p}_{0}}( t)$};
\draw (298,193.4) node [anchor=north west][inner sep=0.75pt]    {$\pi _{\mathbf{p}_{0}}$};
\draw (309,143.4) node [anchor=north west][inner sep=0.75pt]    {$d$};
\draw (537,139.4) node [anchor=north west][inner sep=0.75pt]  [font=\small]  {$\delta $};
\draw (259,12.4) node [anchor=north west][inner sep=0.75pt]    {$N_{\rho }( \ell _{0})$};
\draw (95,148.0) node [anchor=north west][inner sep=0.75pt]    {$\rho $};
\draw (142,201.4) node [anchor=north west][inner sep=0.75pt]    {$\ell '$};
\draw (143,73.4) node [anchor=north west][inner sep=0.75pt]    {$P$};
\draw (111,81.4) node [anchor=north west][inner sep=0.75pt]    {$\ell _{0}$};
\draw (410,151) node [anchor=north west][inner sep=0.75pt]   [align=left] {Cross-section};
\draw (367,216.4) node [anchor=north west][inner sep=0.75pt]    {$N_{\delta }( \pi _{\mathbf{p}_{0}}( \ell '))$};
\draw (445,324.4) node [anchor=north west][inner sep=0.75pt]    {$N_{\delta }( \pi _{\mathbf{p}_{0}}( \ell _{0}))$};
\draw (173,199.4) node [anchor=north west][inner sep=0.75pt]    {$P'$};

\end{tikzpicture}
    \caption{The grain structure of $\mathcal C(X)$ inside $N_\rho(\ell_0)$. The curve $\ell_0$ is black with gray highlighting. Another curve $\ell' \in B(\ell_0, \rho)$ is drawn in black. The prisms in $\mathcal P_\rho(\ell_0)$ intersecting $\ell_0$ and $\ell'$ are labeled $P$ and $P'$ respectively, and are gray and striped. To the right, we show a cross section of $N_\rho(\ell_0)$. The cross sections of the grains are parallel up to error $\delta$, and pointing in direction $\gamma_{\p_0}(t)$. At the bottom is the projection of the grains in $N_\rho(\ell_0)$ under $\pi_{\p_0}$. The prisms corresponding to $\ell_0$ project to the highlighted black line, and the prisms corresponding to $\ell'$ project to the black curve. This captures the geometry used in Lemma \ref{lem: add prisms fat tube}.}
    \label{fig:grain_geometry_diagram}
\end{figure}

\begin{proof}
    The proof is almost exactly the same as in \cite[Lemma 2.12]{KWZ}, and we include below the small adjustments needed to prove our Lemma \ref{lem: add prisms fat tube}.
    \begin{itemize}
        \item The prisms $P$ in \cite{KWZ} are dimension $\delta \times t \times \delta/t$ and ours are dimension $\delta \times \rho \times d$. We therefore replace $t$ by $\rho$ in the proof, and replace the tube segment length $\delta / t$ by $d$.
        \item The application of \cite[Lemma 2.2]{KWZ} is replaced by Lemma \ref{lem: prism filled up}, and we may use it because of the prism compatibility (Lemma \ref{lem: prism compatibility}). 
    \end{itemize}
\end{proof}
In Figure \ref{fig:grain_geometry_diagram}, we give a diagram of the grain structure of $\mathcal C(X)$ around a fixed curve $\ell_0 = \ell_{\p_0}$.

The plane curves on the RHS of \eqref{eq: add prisms fat tube} are contained in the family $\mathcal F(\p_0, X)$, which can be estimated by Theorem \ref{thm: zahl planar tweak} when $X$ is $\tau$-twisty. However the curves must satisfy the one-dimensional ball condition. The next Lemma says that an arbitrary set can be decomposed into pieces which satisfy this condition. This lemma is directly from \cite{KWZ}. 
\begin{lemma}[\protect{\cite[Lemma 3.2]{KWZ}}]\label{lem: 1 dim cond decomp}
    Let $X \subset \R^d$ be $\delta$-separated and finite, and let $K \geq 1$. Then we can partition $X = Y \sqcup Z$, and $Z$ can be further decomposed into $Z = Z_1 \sqcup \cdots \sqcup Z_N$ so that the following holds: 
    \begin{itemize}
        \item 
        \begin{align}
            \sup_{r \geq \delta, x \in \R^d} \frac{\#(B(x,r) \cap Y)}{r/\delta}\leq K. 
        \end{align}
        \item 
        Each $Z_i$ is contained in a ball of radius $r_i \gtrsim \delta K^{1/d}$, and each $Z_i$ has a subset $Z_i'$ of size $\# Z_i' \sim r_i /\delta$ with 
        \begin{align}
            \sup_{r \geq \delta, x \in \R^d} \frac{\#(B(x,r) \cap Z_i')}{r/\delta}\lessapprox 1. 
        \end{align}
    \end{itemize}    
\end{lemma}

Theorem \ref{thm: zahl planar tweak}, Lemma \ref{lem: twist implies 2nd ord tangency}, Lemma \ref{lem: 1 dim cond decomp}, and Lemma \ref{lem: add prisms fat tube} can be combined to give the following estimate. This is analogous to \cite[Proposition 3.3]{KWZ}:

\begin{proposition}[Analogous to \protect{\cite[Proposition 3.3]{KWZ}}]\label{prop: narrow prop}
    Suppose that $X$ is $\tau$-twisty. 
    For all $\epsilon > 0$ the following holds. Let $\delta \in (0,1]$, let $r \in [\delta, 1]$, and let $\rho \in [\delta, \min((\delta r)^{1/2}, \delta/r)]$. Let $(L,Y)$ and $(L',Y')$ be sets of $\delta$-separated curves with regular, $\lambda$-dense shadings. Suppose that $L' \subset L \cap B(\ell_0, \rho)$ and 
    \begin{align}
        Y'(\ell') \subset \{\x : r_L(\x) \sim r, \theta(\ell') \in J_L(\x) \}
    \end{align}
    for each $\ell' \in L'$. Fix $K \geq 1$. Then at least one of the following must happen.
    \begin{enumerate}
        \item[(A)] There is a curve $\tilde \ell \in B(\ell_0, 2\rho)$ and a $\Omega(\lambda^4)$-dense shading $\tilde Y$ of $\tilde \ell$ by $\rho$-cubes with 
        \begin{align}
            |Q \cap E_L| \gtrsim_\epsilon \delta^{\epsilon} (\tau\lambda)^6 K^{-1} ((\delta/\rho) \# L') |Q| \text{ for each $\rho$-cube $Q \subset \tilde Y(\tilde \ell)$.}
        \end{align}
        \item[(B)] There is $\rho' \in [\delta K^{1/10}, \min((\delta r)^{1/2}, r / \delta)]$ and an arrangement $(\tilde L, \tilde Y)_{
        \rho'}$ of curves with $\Omega^*(\lambda^4)$-dense shading. For each $\tilde \ell \in \tilde L$ we have 
        \begin{align}
            |Q \cap E_L| \gtrsim_\epsilon  \delta^\epsilon (\tau \lambda)^6  |Q| \text{ for each $\rho'$-cube $Q \subset \tilde Y(\tilde \ell)$.} 
        \end{align}
        Also $\Omega^*(\# L')$ of the curves $\ell' \in L'$ satisfy $
        \ell' \in B(\tilde \ell, O(\rho'))$, for some $\tilde \ell \in \tilde L$.
    \end{enumerate}
    \begin{proof}
        We give the proof of Case (A) only, since one can make similar adjustments to Case (B) and follow the proof of \cite[Proposition 3.3]{KWZ}. 

        Apply Lemma \ref{lem: 1 dim cond decomp} to $L'$ with parameter $K$; we obtain a partition $L' = L^* \sqcup L^{**}$, where $L^{**} = L_1^{**} \sqcup \cdots \sqcup L_M^{**}$. 

        \medskip
        
        \noindent \textbf{Case (A):} $\# L^* \geq \frac{1}{2} \# L'$.

        Apply part 2 of Lemma \ref{lem: twist implies 2nd ord tangency} to $L^*$ (or more precisely to the corresponding parameters). Then $F$ satisfies the one-dimensional ball condition with $K = O(\tau^{-3} K)$. For each $\ell \in L^*$ define $Y^*(\ell) = Y'(\ell)$. Applying Lemma \ref{lem: add prisms fat tube} with $L$ as above and $(L^*, Y^*)$ in place of $(L',Y')$, we get a curve $\tilde \ell \in B(\ell_0, 2\rho)$ and a $\Omega^*(\lambda^4)$-dense shading $\tilde Y(\tilde \ell)$ by $\rho$-cubes. After replacing $(L^*, Y^*)$ with a $\Omega^*(1)$-density refinement, we have that for each $t$-cube $Q \subset \tilde Y(\tilde \ell)$, 
        \begin{align}
            |E_L \cap Q| &\gtrapprox \lambda^3 \rho^2 \Big | \bigcup_{\ell \in L^*} \pi_{\ell_0}(Y^*(\ell))\Big | \\
            &\gtrapprox_\epsilon (\lambda^3 \rho^2) (\delta^{\epsilon/2} \tau^6 \lambda^3 K^{-1} (\delta \#L^*)) \\
            &\gtrsim_\epsilon \delta^{\epsilon} (\tau \lambda)^6 K^{-1} ((\delta / \rho) \# L') |Q|.
        \end{align}
        In the second line we applied Theorem \ref{thm: zahl planar tweak}, which is justified by part 1 of Lemma \ref{lem: twist implies 2nd ord tangency}.
    \end{proof}
\end{proposition}

\begin{remark}
    Proposition \ref{prop: narrow prop} is the only place where the $\tau$-twistiness condition is used. 
\end{remark}

\section{Proof of Theorem \ref{thm: CT 3 param kak}}\label{sec: 6}

In this section, we combine Proposition \ref{prop: narrow prop} with 3-linear Kakeya (Theorem \ref{thm: multilinear Kakeya}) and induction on scales to prove Theorem \ref{thm: CT 3 param kak}. 
In \cite{KWZ} the multilinear Kakeya theorem was not required. The family $\mathcal L_{\rm{SL}_2}$ is planey, which allows one to take $\rho \in [\delta, (\delta r)^{1/2}]$ instead of in $[\delta, \min((\delta r)^{1/2}, \delta/r)]$ in Proposition \ref{prop: narrow prop}. This makes the argument in \cite{KWZ} efficient. When $r \leq \rho^{a} = \rho^{0.303 \ldots}$ (in a given step of the induction), we run the Katz-Wu-Zahl argument. When $r \geq \rho^{a}$, we apply multilinear Kakeya. When inducting, we need to ensure that the curves at the coarser scale continue to satisfy the two-dimensional ball condition. The following Lemma taken directly from \cite{KWZ} lets us do this. 

\begin{lemma}[\protect{\cite[Lemma 4.1]{KWZ}}]\label{lem: 2 dim cond refinement}
    Let $0 < \delta \leq \rho \leq r \leq 1$ and let $L, \tilde L \subset \mathcal C(X)$. 
    \begin{enumerate}
        \item Suppose that $\tilde L$ is $\rho$-separated and contained in a ball of radius $r$. Then there exists a set $\tilde L' \subset \tilde L$ with $\# \tilde L' \gtrsim (\rho /r) \# \tilde L$ satisfying the two-dimensional ball condition $\# (\tilde L' \cap B(\ell, s)) \leq (s / \rho)^2$. 
        \item Suppose that $L$ obeys the two-dimensional ball condition $\#(L \cap B(\ell, s)) \leq (s / \delta)^2$ and $L \subset \cup_{\tilde \ell \in \tilde L} B(\tilde \ell, \rho)$. Then for all $\epsilon > 0$, there is $c = c(\epsilon) > 0$ and a set $\tilde L' \subset \tilde L$ with $\# \tilde L' \gtrsim c\delta^\epsilon (\delta / \rho)^2 \# L$ that satisfies the two-dimensional ball condition $\#(\tilde L' \cap B(\ell, s)) \leq (s/\rho)^2$. 
    \end{enumerate}
\end{lemma}

Below is the modified version of Theorem \ref{thm: CT 3 param kak} which works well with induction on scales. 
\begin{proposition}\label{prop: inductive CT 3 param kak}
    Let $X$ be $\tau$-twisty, $\mathfrak c$-coney, and satisfy the basic conditions.  
    There exists $c_1 = c_1(\epsilon)$
    and $c_2 = c_2(\epsilon)$ so that the following holds for all $0 < \delta \leq \rho \ll 1$ and $0 < \lambda \leq 1$. 

    Let $(L,Y)_\rho$ be a $\rho$-separated collection of curves satisfying the two-dimensional ball condition $\#(L \cap B_r) \leq (r / \rho)^2$, with $\lambda$-dense shading. Assume in addition that $\theta(L) = \{\theta(\ell) : \ell \in L\}$ is contained in an interval of length $\ll \mathfrak c/|\log \rho|^{3}$. Then  
    \begin{align}\label{eq: induct bound}
        |E_L| \geq c_1 \rho^\epsilon \delta^{2\epsilon} \mathfrak c (\tau \lambda)^W \rho^a (\rho^2 \# L)^{1-c_2},  
    \end{align}
    where 
    \begin{align}
        W = W(\epsilon, \delta, \rho) = \exp(\frac{100}{\epsilon^3} \frac{\log \rho}{\log \delta}),
    \end{align}
    and $a = (\sqrt{13} - 3)/2 \approx 0.303$ is the positive root of $a^2 + 3a - 1$.
\end{proposition}

Let us first see how Proposition \ref{prop: inductive CT 3 param kak} implies Theorem \ref{thm: CT 3 param kak}. 
\begin{proof}[Proposition \ref{prop: inductive CT 3 param kak} implies Theorem \ref{thm: CT 3 param kak}]
    Partition the parameter space $\p = (p,\theta) \in B_1^3$ into $O(|\log \delta|^3/\mathfrak c)$ many $(1 \times 1 \times \mathfrak c / |\log\delta|^3)$-slabs,
    and let $L = L_1 \sqcup \cdots \sqcup L_N$ be the corresponding partition of the set of curves. Choose the largest $L' = L_i$, so $\# L' \gtrsim \mathfrak c |\log \delta|^{-3} \# L$. Apply Proposition \ref{prop: inductive CT 3 param kak} with $(L,Y)$ replaced by $(L',Y)$, $\rho$ replaced by $\delta$, and $\epsilon$ replaced  by $\epsilon / 10$ to get 
    \begin{align}
        |E_{L}| &\geq |E_{L'}| \gtrsim_\epsilon \delta^\epsilon \mathfrak c (\tau \lambda)^W \delta^a (\delta^2 \# L') \\
        &\gtrsim_\epsilon \delta^{\epsilon/ 2} \mathfrak c^2 (\tau \lambda)^W \delta^a (\delta^2 \#L),
    \end{align}
    which gives the claim after taking $\delta$ small enough compared to $\epsilon$. 
\end{proof}

\begin{proof}[Proof of Proposition \ref{prop: inductive CT 3 param kak}] 
    Steps 1 through 7 are essentially the same as steps 1 through 7 in \cite[Proposition 4.2]{KWZ}, with the exception of Step 2, where we use the broad-narrow angle $r_L(\x)$. Step 8 is the new ingredient, where we combine the $\mathfrak c$-coniness assumption with the 3-linear Kakeya theorem.

    \medskip

    \noindent \textbf{Step 1: Setting up the induction} (following \cite[Step 1 pg. 16]{KWZ}):

    If $L$ is empty there is nothing to show. If $L$ is nonempty, we can lower-bound the volume of $E_L$ by the volume of one tube. This gives the estimate
    \begin{align}\label{eq: step 1 1}
        |E_L| \geq \lambda \rho^2. 
    \end{align}
    If $\delta \geq c_1^{1/\epsilon}$ then the induction closes. By choosing $c_1(\epsilon)$ small enough, we may assume that $\delta$ is very small compared to $\epsilon$. For $\delta > 0$ fixed, we prove the result by induction on $\rho$. If $\rho \geq \delta^\epsilon$, we are done by \eqref{eq: step 1 1}. By setting $c_1$ appropriately, we can assume that $\rho > 0$ is so small that 
    \begin{align}
        (2^{100|\log \rho| / \log |\log \rho|})^{-J} \geq \rho^{1/J}, \quad J = \max(c_2^{-2}, \exp(100/\epsilon^3)). 
    \end{align}
    In practice, we may choose $c_2(\epsilon) = \epsilon^6$. As a consequence, expressions like $(2^{100 |\log \rho| /\log |\log \rho|})^{-W(\epsilon, \delta, \rho)}$ are bounded below by $\rho^{c_2^2}$.
    For $\delta$ and $\rho$ fixed, we prove Proposition \ref{prop: inductive CT 3 param kak} by induction on $\# L$. If $\# L \leq \rho^{-\epsilon /2}$ and $L$ is nonempty, we are done by \eqref{eq: step 1 1}. Assume that $\# L > \rho^{-\epsilon / 2}$ from now on.

    \medskip
    
    \noindent \textbf{Step 2: Finding a typical angle $r$} (following \cite[Step 2 pg. 16]{KWZ})

    Apply Lemma \ref{lem: regular refinement} to each $\ell \in L$ and let $Y_1(\ell) \subset Y(\ell)$ be a $\lambda/ 2$-dense regular shading. Set $L_1 = L$. Decompose $Y_1(\ell)$ into $Y_1^{low}(\ell) \sqcup Y_1^{high}(\ell)$, where 
    \begin{align}
        Y_1^{low}(\ell) = \{\x \in Y_1(\ell) : \# L(\x) \leq \rho^{-\epsilon} \}
    \end{align}
    and $Y_1^{high}(\ell) = Y_1(\ell) \setminus Y_1(\ell)^{low}$. If $\sum_{\ell \in L_1} |Y_1^{low}(\ell)| \geq \sum_{\ell \in L_1} |Y_1^{high}(\ell)|$ then 
    \begin{align}
        |E_L| \geq \rho^\epsilon \frac{\lambda}{4} \rho^2 \# L
    \end{align}
    by double counting. If $c_2 \leq 1/4$ the induction closes. Assume that we are in the high multiplicity case and define $Y_2(\ell) = Y_1^{high}(\ell)$. Set $L_2 = L_1$. By swapping the order of summation, we can compute
    \begin{align}
        \sum_r \sum_{\ell \in L_1}& |\{\x \in Y_2(\ell) : r_{L_1}(\x) \sim r \text{ and } \theta(\ell) \in J_{L_1}(\x)\}| 
        \\ &= \int_{\{\x : \#L_1(\x) > \rho^{-\epsilon}\}} \# \{\ell \in L_1(\x) : \theta(\ell) \in J_{L_1}(\x)\} \\
        &\gtrsim \int_{\{\x : \# L_1(\x) > \rho^{-\epsilon}\}} 100^{-|\log \rho| / \log |\log \rho|} \# L_1(\x) d\x \\
        &\gtrapprox \lambda \rho^2 \# L. 
    \end{align}
    Pigeonhole $r = |\log \rho|^{-j}$ so that 
    \begin{align}
        \sum_{\ell \in L_1} |\{\x \in Y_2(\ell) : r_{L_1}(\x) \sim r \text{ and } \theta(\ell) \in J_{L_1}(\x)\}| \gtrapprox \lambda \rho^2 \# L.
    \end{align}
    For each $\ell \in L_1$ define 
    \begin{align}
        Y_3(\ell) = \{\x \in Y_2(\ell) : r_{L_1}(\x) \sim r \text{ and } \theta(\ell) \in J_{L_1}(\x) \}.
    \end{align}
    If we set 
    \begin{align}
        L_3 = \{\ell \in L_2 : |Y_3(\ell)| \gtrapprox \lambda \rho^2\},
    \end{align}
    and choose the implied factors appropriately, we find 
    \begin{align}
        \# L_3 \gtrapprox \# L. 
    \end{align}
    We have $\# L_3(p) \lesssim r/\rho$, so 
    \begin{align}\label{eq: trivial narrow}
        |E_{L}| \geq |E_{L_3}| \gtrsim (\rho / r) \sum_{\ell \in L_3} |Y_3(\ell)| \gtrapprox (\rho / r) \lambda (\rho^2 \# L). 
    \end{align}
    Either $r \geq \rho^{1-a}$ or the induction closes.
    Apply Lemma \ref{lem: regular refinement} to each $\ell \in L_3$ and let $Y_4(\ell)$ be the $\Omega^*(\lambda)$-dense regular shading. Define $L_4 = L_3$. 

    \medskip
    
    \noindent \textbf{Step 3: A large subset is contained in a $1 \times r$ tube} (following \cite[Step 3 pg. 17]{KWZ})
    Let $\mathcal B$ be a set of $O(r^{-3})$ boundedly overlapping balls whose union covers $L_4$. For each $B \in \mathcal B$, let $L_4^B = L_4 \cap B$. There are $O(1)$ balls $B \in \mathcal B$ with $\x \in E_{L_4^B}$. Thus 
    \begin{align}
        \Big |\bigcup_{\ell \in L_4} Y_4(\ell)\Big | \gtrsim \sum_{B \in \mathcal B} \Big |\bigcup_{\ell \in L_4^B} Y_4(\ell)\Big | .
    \end{align}
    Continuing to follow \cite{KWZ}, set $N = \sup_{B \in \mathcal B} \# L_4^B$. 
    If $N \leq \rho^{c_2} \# L$ we can apply the inductive hypothesis to find
    \begin{align}
        \sum_{B \in \mathcal B} &\Big|\bigcup_{\ell \in L_4^B} Y_4(\ell)\Big | \geq c_1 \rho^\epsilon \rho^a \delta^{2\epsilon} \mathfrak{c} (2^{-10 |\log \rho|/\log|\log \rho|}
        \lambda \tau)^{W(\epsilon, \delta, \rho)} \sum_{B \in \mathcal B} (\rho^2 \# L_4^B)^{1-c_2} \\ 
        &\geq c_1 \rho^\epsilon \rho^a \delta^{2\epsilon} \mathfrak{c}(\lambda \tau)^{W(\epsilon, \delta, \rho)} (2^{-10 |\log \rho|/\log|\log \rho|})^{W(\epsilon, \delta, \rho)}
        (\rho^2 \# L_4)^{1-c_2} (\# L_4 / N)^{c_2} \\
        &\geq c_1 \rho^\epsilon \rho^a \delta^{2\epsilon} \mathfrak{c} (\lambda \tau)^{W(\epsilon, \delta, \rho)}  (2^{-10 |\log \rho|/\log|\log \rho|})^{W(\epsilon, \delta, \rho) + 1} (\rho^2 \# L)^{1-c_2}(\# L / N)^{c_2}.
    \end{align}
    The term in the final brackets is at least  $\rho^{c^2_2} (\# L / N)^{c_2} \geq 1$, so the induction closes. Henceforth assume that there is a set $L_5 = L_4^{B_0}$ satisfying 
    \begin{align}\label{eq: size of L}
        \# L_5 \geq \rho^{c_2} \# L.
    \end{align}
    Define $Y_5(\ell) = Y_4(\ell)$ for each $\ell \in L_5$. 

    \medskip 

    \noindent \textbf{Step 4: Clustering at scale $\tilde \rho$} (following \cite[Step 4 pg. 17]{KWZ})

    Let $\mathcal B_0$ be a cover of $\mathcal C(X) \cap B_0$ by boundedly overlapping balls of radius $\tilde \rho = \min((\rho r)^{1/2}, \rho / r)$. Pigeonhole an integer $M$ and subset of balls $\mathcal B \subset \mathcal B_0$ so that 
    \begin{align}
        \sum_{B \in \mathcal B} \# (L_5 \cap B) \geq |\log \rho|^{-1} \# L_5
    \end{align}
    and 
    \begin{align}
        \#(B \cap L_5) \in [M, 2M) \text{ for all $B \in \mathcal B$.}
    \end{align}
    Define 
    \begin{align}
        \tilde L = \{\ell_B : B \in \mathcal B\},
    \end{align}
    where $\ell_B \in \mathcal C(X)$ is the center of $B$. For each $\tilde \ell \in \tilde L$, define $L_5(\tilde \ell) = L_5 \cap B(\tilde \ell, \tilde \rho)$. Let $(L_5(\tilde \ell), Y_5)_{\rho}$ be the restricted pair of curves and their shadings induced from $Y_5$. We have 
    \begin{align}
        M \cdot (\# \tilde L) \gtrsim |\log \rho|^{-1} (\# L_5).
    \end{align}

    \medskip
    
    \noindent \textbf{Step 5: Apply Proposition \ref{prop: narrow prop}} (Following \cite[Step 5 pg. 18]{KWZ})

    For each $\tilde \ell \in \tilde L$, apply Proposition \ref{prop: narrow prop} with 
    \begin{itemize}
        \item $\rho$ in place of $\delta$; $2^{-10 |\log \rho|/\log |\log \rho|} \lambda$ in place of $\lambda$; $r$ as above, and $\tilde \rho$ in place of $\rho$. 
        \item $(L_1,Y_1)$ in place of $(L,Y)$.
        \item $(L_5(\tilde \ell), Y_5)$ in place of $(L', Y')$.
        \item $\tilde \ell$ in place of $\ell_0$.
        \item $K = \rho^{-\epsilon^3}$ and $\epsilon_1 = \epsilon^4 / 40$ in place of $\epsilon$. 
    \end{itemize}

    \noindent \textbf{Case (A)}: If $\tilde \ell$ is of Type (A), there is a shading $\tilde Y(\tilde \ell)$ that is a union of $\tilde \rho$-cubes, with $|\tilde Y(\tilde \ell)| \gtrapprox \lambda^4 {\tilde \rho}^2$, so that for each $\tilde \rho$-cube $Q \subset \tilde Y(\tilde \ell)$, we have 
    \begin{align} \label{eq: case (A)}
        |Q \cap E_L| &\geq |Q \cap E_{L_1}| \geq c_\epsilon \rho^{\epsilon_1}  K^{-1} (\rho^{O(1/\log|\log \rho|)} \lambda \tau)^6 ((\rho / \tilde \rho) M) |Q| \\
        &\geq c_{\epsilon} \rho^{\epsilon^4/20 + \epsilon^3} (\lambda \tau)^6 ((\rho / \tilde \rho) M) |Q|.
    \end{align}

    \noindent \textbf{Case (B)}: If $\tilde \ell$ is of Type (B), there is a diameter $\rho' = \rho'(\tilde \ell) \in [\rho^{1-\epsilon^3/10}, \min((\rho r)^{1/2}, \rho/r)]$ 
    and set $(L^*_{\tilde \ell}, Y_{\tilde \ell}^*)_{\rho'}$ of curves with a $\Omega^*(\lambda^4)$ dense shading. For each $\ell^* \in L_{\tilde \ell}^*$ and each $\rho'$-cube $Q \subset Y_{\tilde \ell}^*(\ell^*)$, we have
    \begin{align}
        |Q \cap E_L| \geq |Q \cap E_{L_1}| \geq c_{\epsilon} \rho^{\epsilon^4/30}  (\lambda \tau)^6 |Q|. 
    \end{align}
    And $\Omega^*(\#L_5(\tilde \ell))$ of the curves satisfy $\ell \in B(\ell^*, O(\rho'))$ for some $\ell^* \in L_{\tilde \ell}^*$. 

    \medskip

    \noindent \textbf{Step 6: Close induction when $r \leq \rho^a$ and Case (A) holds} (Following \cite[Step 6 pg. 18]{KWZ})

    Suppose that at least half the curves from $\tilde L$ are of Type (A), and refined $\tilde L$ to contain only those curves. Apply Lemma \ref{lem: 2 dim cond refinement} with $L_5$ in place of $L$; $\tilde L$ as above; $\rho$ in place of $\delta$; $\tilde \rho$ in place of $\rho$; and $r$ as above. We obtain $\tilde L' \subset \tilde L$ obeying the two-dimensional ball condition and has size
    \begin{align}
        \# \tilde L' \gtrsim (\tilde \rho / r) (\# \tilde L) \gtrapprox_\rho (\tilde \rho / r) N/M. 
    \end{align}
    Apply the inductive hypothesis at scale $\tilde \rho$ to $(\tilde L', \tilde Y)$ with $\tilde \lambda = \Omega^*(\lambda^4)$:
    \begin{align} \label{eq: case (A) 2}
        \Big |\bigcup_{\tilde \ell \in \tilde L'} \tilde Y(\tilde \ell)\Big | &\geq c_1 \mathfrak c\tilde \rho^a \delta^{2\epsilon} \tilde \rho^\epsilon (\tilde \lambda \tau)^{W(\epsilon, \delta, \tilde \rho )}(\tilde \rho^2 \# \tilde L')^{1-c_2} \\
        &\geq \rho^{-\epsilon^2 / 2} c_1 \mathfrak c \rho^\epsilon \delta^{2\epsilon} \tilde \rho^a (\lambda\tau)^{\frac{1}{2} W(\epsilon, \delta, \rho)} (\tilde \rho^2 (\tilde \rho/ r) N / M)^{1- c_2}. 
    \end{align}
    The set on the LHS of 
    \eqref{eq: case (A) 2} is a union of $\tilde \rho$-cubes $Q$, and \eqref{eq: case (A)} gives a lower bound for $|Q \cap E_L|$ in terms of $|Q|$. Therefore 
    \begin{align}
        |E_L| &\geq (c_\epsilon \rho^{\epsilon^4 / 20 +\epsilon^3} \lambda^6 (\rho / \tilde \rho) M) \Big |\bigcup_{\tilde \ell \in \tilde L'} \tilde Y(\tilde \ell)\Big | \\
        &\geq (c_\epsilon \rho^{\epsilon^4 / 20 +\epsilon^3} \lambda^6 (\rho / \tilde \rho) M) (\rho^{-\epsilon^2 / 2} c_1 \mathfrak c \rho^\epsilon \delta^{2\epsilon} \tilde \rho^a (\lambda \tau)^{\frac{1}{2} W(\epsilon, \delta, \rho)} (\tilde \rho^2 (\tilde \rho/ r) N / M)^{1- c_2}) \\
        &\geq (c_\epsilon \rho^{\epsilon_1 + 4c_2 + \epsilon^3 - \epsilon^2 /2})(c_1 \mathfrak c \rho^\epsilon \delta^{2\epsilon} (\lambda \tau)^{W(\epsilon, \delta, \rho)} \tilde \rho^{2 + a}/(\rho r))(\rho^2 \# L)^{1-c_2}. 
    \end{align}
    The first term in brackets is $\geq 1$ as long as $\epsilon$ is sufficiently small ($\epsilon < 1/100$ works) and $c_1(\epsilon)$ is taken sufficiently small, which forces $\rho$ to be sufficiently small. 
    When $r \leq \rho^{1/3}$, $\tilde \rho = (\rho r)^{1/2}$. Then we get $\tilde \rho^{2 + a} / (\rho r) = (\rho r)^{a/2} \geq \rho^a$, so the induction closes. When $r \geq \rho^{1/3}$ we have $\tilde \rho = \rho / r$. Since $a^2 + 3a - 1 = 0$ by assumption, 
    \begin{align}
        \tilde \rho^{2 + a} / (\rho r) \geq \rho^{-a^2 - 3a + 1} \rho^a = \rho^a,
    \end{align}
    and the induction closes.

    \medskip 
    
    \noindent \textbf{Step 7: Close induction in case (B)}(Following \cite[Step 7 pg. 19]{KWZ})

    This case is almost identical to \cite{KWZ}. One only needs to track the extra factor $\mathfrak c \tau^W \rho^a$ in their argument. We omit this easy extension. 

    \medskip

    \noindent \textbf{Step 8: Close induction in case $r \geq \rho^a$}

    Set $f = \sum_{\ell \in L_5} \chi_{Y_5(\ell)}$. 
    By H\"older's inequality, 
    \begin{align}\label{eq: pf broad 1}
        |E_L| \geq |E_{L_5}| \gtrapprox \frac{(\lambda \rho^2 \# L_5)^3}{(\int f^{3/2})^2}.
    \end{align}
    Fix $\x \in E_{L_5}$. Since $Y_5(\ell) \subset \{\x \in Y_1(\ell) : r_{L_1}(\x) \sim r\}$, there is an interval $J$ of length $r$ and 3 subsets $L_1'(\x),L_1''(\x), L_1'''(\x) \subset (L_1)_J(\x)$ which are $r/|\log \rho|$ $\theta$-separated and have size $\gtrapprox \# L_1(\x)$. 
    By assumption $\theta(L)$ is contained in an interval of length $\ll \mathfrak c / |\log \rho|^3$, so $r \ll \mathfrak c/|\log \rho|^3$.
    Lemma \ref{lem: coniness means transverse} gives 
    \begin{align}
        |\det(\bn_{\ell'}(\x), \bn_{\ell''}(\x), \bn_{\ell'''}(x))| \gtrapprox \mathfrak{c} r^3
    \end{align}
    whenever $\ell' \in L_1'(\x), \ell'' \in L_1''(\x)$, and $\ell''' \in L_1'''(\x)$. 
    Let $T_\ell = N_{C\rho}(\ell)$ for a sufficiently large constant $C$. It follows that 
    \begin{align}
        f(\x)^3  &\lessapprox \mathfrak c^{-1} r^{-3} \sum_{\ell', \ell'', \ell''' \in L_1} \chi_{T_{\ell'}}(\x) \chi_{T_{\ell''}}(\x) \chi_{T_{\ell'''}}(\x)|\bn_{\ell'}(\x) \wedge \bn_{\ell''}(\x) \wedge \bn_{\ell'''}(\x)|. \label{eq: pf broad 2}
    \end{align}
    Integrate the square root of \eqref{eq: pf broad 2} and apply 3-linear Kakeya (Theorem \ref{thm: multilinear Kakeya}) with $\epsilon$ replaced by $\epsilon_1/2$: 
    \begin{align}\label{eq: step 8 1}
        \int f^{3/2} \lessapprox_{\epsilon_1} \rho^{-\epsilon_1/2} \mathfrak c^{-1/2} r^{-3/2} \rho^3 (\# L_1)^{3/2}. 
    \end{align}
    Combining \eqref{eq: step 8 1} with \eqref{eq: pf broad 1}, we find 
    \begin{align} 
        |E_L| &\gtrapprox_{\epsilon_1} \rho^{\epsilon_1} \mathfrak c \lambda^3 r^3 (\# L_5/\#L)^3 \\
        &\gtrapprox_{\epsilon_1} \rho^{\epsilon_1 + 3c_2} \mathfrak{c} \lambda^3 r^3,
        \label{eq: step 8 2}
    \end{align}
    since $\# L_5 \gtrapprox \# L_1$. Since $\# L_5$ satisfies the 2-dimensional ball condition and is contained in a ball of radius $r$, and \eqref{eq: size of L}, we get 
    \begin{align}
        r^2 \geq \rho^{c_2} \rho^2 \# L. 
    \end{align}
    Combining this with \eqref{eq: step 8 2},
    \begin{align}
        |E_L| &\gtrapprox_\epsilon \rho^{\epsilon_1} r \mathfrak c \lambda^3 \rho^{c_2} \rho^2 \# L \\
        &\gtrapprox_\epsilon \rho^{\epsilon_1 + 3 c_2} \mathfrak c \lambda^3 \rho^a (\rho^2 \# L)^{1-c_2}.
    \end{align}
    If we choose $\epsilon_1 = \epsilon / 10$, then $\rho^{\epsilon_1 + 3 c_2} \geq \rho^{\epsilon / 2}$. If we choose $c_1(\epsilon)$ sufficiently small, which forces $\rho$ to be sufficiently small, the induction closes.
\end{proof}

\section{Proof of Theorem \ref{thm: main phicurved kakeya}}\label{sec: 7}

\subsection{Rescaled $\mathcal C(\phi)$ curves}

Given a translation invariant phase function 
\begin{align}
    \phi(x,t,\xi) = \langle x, \xi \rangle + \psi(t,\xi),
\end{align}
recall that the family $\mathcal C(\phi)$ consists of curves of the form 
\begin{align}
    \ell_{\xi, v}(t) &= (v - \nabla_\xi \psi(t,\xi), t) \\
    &= (\Phi(\xi, v, t), t),
\end{align}
where $\Phi(\xi, v, t) = v - \nabla_\xi \psi(t,\xi)$. 
In the next subsection, we will formulate a two-ends version of Theorem \ref{thm: main phicurved kakeya}. We will therefore need to work with a rescaled version of $\mathcal C(\phi)$ inside a radius $\tau$ ball $B_{\tau}$.  
\begin{definition}[Rescaled family $\mathcal C^{\tau}(t_0, \phi)$]
    Fix $t_0 \in B_1^1$ and $\tau \in (0,1]$. The family $\mathcal C^{\tau}(t_0, \phi)$ consists of curves 
    \begin{align}
        \ell_{\xi, v}(t) = (\tau^{-1} \Phi(\xi, \tau v, \tau t), t), 
    \end{align}
    where $t \in B(t_0,1)$. Also define 
    \begin{align}
        \Phi^\tau(\xi, v, t) = \tau^{-1} \Phi(\xi, \tau v, \tau t + t_0).
    \end{align}
\end{definition}

\begin{example}
    Consider the model family 
    \begin{align}
        \phi_{A,B}(x,t,\xi) = \langle x, \xi \rangle + \frac{1}{2} \langle (tA + t^2 B) \xi, \xi \rangle. 
    \end{align}
    Then $\ell_{\xi, v} \in \mathcal C^\tau (\phi_{A,B}, 0)$ is 
    \begin{align}
        \ell_{\xi, v}(t) = (v - [tA + \tau t^2 B]\xi, t). 
    \end{align}
    When $B = 0$, this is the standard family of lines which rescales well. In general, the family approaches the standard lines as $\tau \to 0$. 
\end{example} 

The next Lemma says that two intersecting curves in $\mathcal C^{\tau}(t_0, \phi)$ escape each other at a rate similar to straight lines. 
\begin{lemma}\label{lem: curve escape rate}
    Let $\phi$ be a translation invariant phase function satisfying H\"ormander's conditions. Take $\ell_1 = \ell_{\xi_1,v_1}, \ell_2 = \ell_{\xi_2, v_2} \in \mathcal C^\tau(t_0, \phi)$ which intersect at height $t$. Then 
    \begin{align}
        |\ell_1(t') - \ell_2(t')| \sim |t' - t| |\xi_1 - \xi_2|,
    \end{align}
    with implicit constant depending only on $\phi$. 
\end{lemma}

\begin{proof}
    We have 
    \begin{align}
        \Phi^\tau(\xi, v, t') = v - \tau^{-1} \nabla_\xi \phi(\tau t', \xi). 
    \end{align}
    Using $\Phi^\tau(\xi_1, v_1, t) = \Phi^\tau(\xi_2, v_2, t)$ to solve for $v_1 - v_2$ and two applications of the fundamental theorem of calculus, we find
    \begin{align}
        \Phi^\tau(\xi_1, v_1, t') - \Phi^\tau(\xi_2, v_2, t') &= \tau^{-1} [\nabla_\xi \phi(\tau t', \xi_2) - \nabla_\xi \phi(\tau t', \xi_1)] \\ &- \tau^{-1} [\nabla_\xi \phi(\tau t, \xi_2)- \nabla_\xi \phi(\tau t, \xi_1)] \\
        &= \int_0^1 \int_0^1 \partial_t \nabla_\xi^2 \psi(\tau (t + s(t'-t)), \xi_1 + (\xi_2 - \xi_1)r)dsdr \label{eq: lem: curve escape rate 1} \\&\cdot (\xi_2 - \xi_1) (t' - t).  \nonumber
    \end{align}
    The matrix $\partial_t \nabla_\xi^2 \psi(s, \xi)$ is nondegenerate by the H\"ormander conditions. The claim follows by taking absolute values of \eqref{eq: lem: curve escape rate 1}. 
\end{proof}

\subsection{Two-ends reduction}

We now discuss the standard two-ends reduction, which lets us assume that not too much mass is concentrated on one end of the tube. 
\begin{definition}[$\epsilon_0$-two-ends]\label{def: two ends}
    We call a shading $Y(\ell) \subset N_\delta(\ell)$ $\epsilon_0$-two-ends if for all $\x \in \R^3$ and $\delta \leq r \leq 1$, 
    \begin{align}
        |Y(\ell) \cap B_r(\x)| \leq r^{\epsilon_0} |Y(\ell)|. 
    \end{align}
\end{definition}
The next Lemma says that we can find a relatively large two-ends piece inside a $\lambda$-dense shading $Y(\ell)$. 

\begin{lemma}[Large two-ends piece \protect{\cite[Lemma 6]{taoTwoEnds}}]\label{lem: two ends piece}
    Let $Y(\ell)$ be a $\lambda$-dense shading and fix $0 < \epsilon_0 < 1$. There is a ball $B(\x, r)$ with $r \in [\lambda, 1]$ and 
    \begin{align}
        |Y(\ell) \cap B(\p, r)| \gtrsim r^{\epsilon_0} |Y(\ell)|. 
    \end{align}
    And for all $\x' \in \R^3$ and $\delta \leq r' \leq 1$ we have 
    \begin{align}
        |Y(\ell) \cap B(\x', r') \cap B(\x, r)| \lesssim (r'/ r)^{\epsilon_0} |Y(\ell) \cap B(\x, r)|. 
    \end{align}
\end{lemma}

We will reduce Theorem \ref{thm: main phicurved kakeya} to a $\epsilon_0$-two-ends version. Suppose we have passed to a $\tau$-ball in the variable $\x \in B^3$. When verifying that hairbrushes are coney later on, we need to assume $L$ is contained in a $c_1 \tau$-ball, where $c_1$ is a small constant depending only on $\epsilon_\phi$, $\epsilon_0$, and other constants depending only on $\phi$.
This will essentially allow us to ignore the terms in $\psi$ of order $O(|\xi|^3)$ when verifying our hairbrushes are twisty and coney in Section \ref{sec: 8}.
\begin{proposition}[Two-ends Version of Theorem \ref{thm: main phicurved kakeya}]\label{prop: two ends version of thm}
    Let $\phi$ be a phase function satisfying the hypotheses in Theorem \ref{thm: main phicurved kakeya}. Fix $\epsilon_0, \epsilon > 0$. There exists a constant $M = M(\epsilon)$ so that the following holds. Fix $\tau \in (0,1]$ and $t_0 \in B_1^1$. Let $(L,Y)$ be a collection of $\delta$-direction separated curves in $\mathcal C^\tau(t_0, \phi)$ with their associated $\lambda$-dense $\epsilon_0$-two-ends shadings. Assume that $L$ is contained in some $c_1 \tau$-ball. Then 
    \begin{align}
        |E_L| \gtrsim_{\epsilon_0, \epsilon} \delta^\epsilon (\tau \lambda)^M \delta^{\frac{a + 1}{2}} (\delta^2 \# L)^{3/4}.
    \end{align}
\end{proposition}

\begin{proof}[Proposition \ref{prop: two ends version of thm} implies Theorem \ref{thm: main phicurved kakeya}]
    
    This reduction is mostly standard, though we won't need to run it as efficiently since large losses in $\lambda$ are acceptable. Apply Lemma \ref{lem: two ends piece} to each $Y(\ell)$ and let $B(\x_\ell, r_\ell)$ be the resulting ball. Define the shading 
    \begin{align}
        Y'(\ell) = Y(\ell) \cap B(\x_\ell, r_\ell).
    \end{align}
    Thus $|Y'(\ell)| \gtrsim r_{\ell}^{\epsilon_0} \lambda \delta^2$. By pigeonholing, we may select $\tau \in [\delta, 1]$ and $L' \subset L$ satisfying  $\# L' \gtrapprox \# L$ and $r_\ell \sim \tau$ for each $\ell \in L'$. Since 
    $|Y'(\ell)| \lesssim \tau \delta^2$,
    we also have 
    \begin{align}
        \tau \gtrsim \lambda^{1/(1-\epsilon_0)}.
    \end{align}
    Let $\mathcal B$ be a cover of $\Bphi^3$ by boundedly overlapping balls of radius $10 \tau$. Each $\ell \in L'$ lies completely inside at least one of the balls, so we have a partition $L' = \bigsqcup_{B \in \mathcal B} L_B'$. Isolate the most popular ball $B \in \mathcal B$ with $\# L_B \gtrapprox \tau^{3} \#L$.

    Let $(x_B,t_B)$ be the center of $B$. After applying the isotropic rescaling 
    \begin{align}
        (x,t) \mapsto (\tau^{-1}(x-x_B), \tau^{-1}t),
    \end{align}
    the curves in $L_B'$ become curves in $\mathcal C^\tau(\tau^{-1}t_B, \phi)$ and the shadings become $\epsilon_0$-two-ends shadings by $\delta/\tau$-cubes. Moreover, 
    \begin{align}
        \tau^{-3} |Y'(\ell)| \gtrsim \lambda \tau^{\epsilon_0-1} (\delta/\tau)^2 \geq \lambda^{1+\epsilon_0} (\delta/\tau)^2.
    \end{align}
    Let $\tilde L_B'$ be the resulting collection of rescaled curves. We may extract a $(\delta/\tau)$-direction separated subset $\tilde L_B'' \subset \tilde L_B'$ with $\# \tilde L_B'' \gtrsim \tau^5 \#L$. 

    Finally, cover the parameter space $\mathcal C^\tau(\tau^{-1}t_B, \phi) \approx \R^4$ by $O((c_1 \tau)^{-4})$ balls of radius $c_1 \tau$, and choose the most popular one. This gives a collection
    $\tilde L_{B,B'}'' \subset \tilde  L_B''$ satisfying the hypotheses of Proposition \ref{prop: two ends version of thm} with $\# \tilde L_{B,B'}'' \gtrapprox\tau^{9} \# L$. 
    Applying Proposition \ref{prop: two ends version of thm} with $\delta$ replaced by $\delta / \tau$, $\lambda$ replaced by $\lambda^{1 + \epsilon_0}$, and $\epsilon$ replaced by $\epsilon/2$, we find 
    \begin{align}
        |E_L| &\geq 
        \tau^3 \Big |\bigcup_{\tilde \ell \in L_{B,B'}''} \tilde Y(\tilde \ell)\Big | \\
        &\gtrsim \tau^3 (\tau \lambda^{1 + \epsilon_0})^M (\delta / \tau)^{\frac{a + 1}{2}} ((\delta / \tau)^2 \# \tilde L_{B,B'}'')^{3/4} \\
        &\gtrsim \delta^{\epsilon/2} \lambda^{3M+10} \delta^{\frac{a + 1}{2}} (\delta^2 \# L)^{3/4}.
    \end{align}
    Absorbing the $\lessapprox$-loss into $\delta^{\epsilon/2}$ gives the result with $M'=3M+10$. 
\end{proof}

\subsection{The bush bound}

In Wolff's hairbrush reduction, we need to estimate a union of hairbrushes. This uses the following bush argument originally due to Bourgain \cite{BourgainBesicovitch}.

\begin{lemma}[Bush Bound for $\mathcal C(\phi)$]\label{lem: bush bound}
    Let $\phi: \Bphi^3 \times \Bphi^2 \to \R$ be a phase function satisfying H\"ormander's conditions. Let $L \subset \mathcal C^\tau(t_0, \phi)$, $E \subset \R^3$, and $0 < \sigma < \epsilon \leq 1$. Assume that 
    \begin{align}
        N_\sigma(\ell) \cap N_\sigma(\ell') \text{ implies } |\xi - \xi'| \geq C \epsilon,
    \end{align}
    and for all $\x \in \R^3$ and $\ell \in L$, 
    \begin{align}
        |N_\sigma(\ell) \cap E \cap (\R^3 \setminus B(\x, \sigma / \epsilon))| \gtrsim \lambda \sigma^2. 
    \end{align}
    Then 
    \begin{align}
        |E| \gtrsim \lambda \sigma^2 (\# L)^{1/2}. 
    \end{align}
\end{lemma}

\begin{proof}
    We have $\sum_{\ell \in L} |N_\sigma(\ell) \cap E| \gtrsim \lambda \sigma^2 \# L$, so there is some point $\x \in E$ belonging to $\gtrsim \lambda \sigma^2 \# L |E|^{-1}$ many $N_\sigma(\ell)$. Let $L_1$ be this collection of curves in the bush through $\x$. Define 
    \begin{align}
        T(\ell) = N_\sigma(\ell) \cap (\R^3 \setminus B(\x, \sigma / \epsilon))
    \end{align}
    for each $\ell \in L_1$. Next we will show that the sets $T(\ell)$ are pairwise disjoint. Choose $\ell_1, \ell_2 \in L_1$. In light of Lemma \ref{lem: metric on C(phi)} (or rather a small adjustment with the parameter $\tau$), we may assume that $\x = (x, t) \in \ell_1 \cap \ell_2$. By Lemma \ref{lem: curve escape rate}, 
    \begin{align}
        |\ell_1(t') - \ell_2(t')| \gtrsim (\sigma / \epsilon)(C \epsilon) = C \sigma
    \end{align}
    for $|t' - t| \geq \sigma / \epsilon$. As long as $C = O(1)$ is large enough, this shows $T(\ell_1) \cap T(\ell_2) = \emptyset$.

    By assumption, 
    \begin{align}
        |T(\ell) \cap E| \gtrsim \lambda \sigma^2. 
    \end{align}
    Therefore $|E| \gtrsim (\lambda \sigma^2) (\lambda \sigma^2 \# L|E|^{-1})$, which upon rearranging gives the lemma.
\end{proof}

\begin{remark}
    The phase function $\phi$ does not need to be translation invariant. Lemma \ref{lem: bush bound} is true under the assumption of H\"ormander's conditions only. 
\end{remark}

\subsection{Wolff's multiplicity dichotomy}

We will follow Wolff's lead to reduce Proposition \ref{prop: two ends version of thm} to an estimate on the volume of a hairbrush of curves, although we of course cannot use his method to estimate the hairbrush. We locate a low and high multiplicity dichotomy just as in \cite{wolffHairbrush}. 

\begin{definition}[Low and high multiplicities]
    Fix a number $\mu$ and consider the following possibilities: 
    \begin{itemize}
        \item \textbf{Case $I$ (low multiplicity):} There are at least $\# L /2$ curves $\ell \in L$ such that 
        \begin{align}
            |\{\x \in Y(\ell) : 
            \# \{\ell' \in L : \x \in N_\delta(\ell')\} \leq \mu\}| \geq \frac{\lambda}{2} \delta^2. 
        \end{align}
        \item \textbf{Case $II_\sigma$ (high multiplicity at angle $\sigma$): } There are at least $C^{-1}|\log \delta|^{-1} \# L$ curves $\ell \in L$ such that 
        \begin{align}
            |\{\x \in Y(\ell) : \#(\ell' \in L : \x\in N_\delta(\ell') \text{ and } \sigma \leq |\xi(\ell)-\xi(\ell')| \leq 2\sigma) 
			\geq (C |\log &\delta|)^{-1} \mu\}|&\label{eq: high mult tube property} \\
			\geq (C |\log \delta|)^{-1} \lambda \delta^2.& \nonumber
        \end{align}
        Write $L'$ for the set of these high-multiplicity curves.
    \end{itemize}
\end{definition}

\begin{lemma}[\protect{\cite[Lemma 3.2]{wolffHairbrush}}]\label{lem: mult dichotomy}
    There is a multiplicity $\mu$ for which both $I$ and $II_\sigma$ hold, for some $\sigma \in [\delta, (c_1 \tau)]$. 
\end{lemma}

The proof is exactly the same as in \cite{wolffHairbrush}, though we reproduce it here.
\begin{proof}
    Take the smallest $\mu$ for which $I$ holds. Then there are $\# L / 2$ many $\ell \in L$ for which 
    \begin{align}\label{eq: lem: mult dichotomy: 1}
        |\{\x \in Y(\ell) : 
            \# \{\ell' \in L : \x \in N_\delta(\ell')\} \geq \mu\}| \geq \frac{\lambda}{2} \delta^2.
    \end{align}
    For any $\ell$ as in \eqref{eq: lem: mult dichotomy: 1} and any $\x \in N_\delta(\ell)$ with $\# \{\ell' \in L : \x \in N_\delta(\ell')\} \geq \mu$, we may pigeonhole a common angle $\sigma$: 
    \begin{align}
        \#\{\ell' \in L : \x \in N_\delta(\ell'), \ \sigma \leq |\xi(\ell') - \xi(\ell)| \leq 2\sigma \} \gtrsim |\log \delta|^{-1} \mu. 
    \end{align}
    Then for any $\ell$ as in \eqref{eq: lem: mult dichotomy: 1} there is some $\sigma$ such that 
    \begin{align}
            |\{\x \in Y(\ell) : \#(\ell' \in L : \x\in N_\delta(\ell) \text{ and } \sigma \leq |\xi(\ell) - \xi(\ell')| \leq 2\sigma) 
			\geq (C |\log \delta|)^{-1} \mu\}|&\\
			\geq (C |\log \delta|)^{-1} \lambda \delta^2&,
        \end{align}
        for a constant $C = O(1)$. Therefore $II_\sigma$ holds with this choice of $\sigma$. 
\end{proof}

We make separate arguments in Case $I$ and $II_\sigma$. Case $I$ does not rely on any geometry. 

\begin{lemma}\label{lem: case I}
    If Case I holds then $|E_L| \gtrsim \lambda \mu^{-1} \delta^2 \#L$. 
\end{lemma}
The proof is just double counting. Case $II_\sigma$ will follow from Theorem \ref{thm: CT 3 param kak} after some work. 

\begin{lemma}\label{lem: case II}
    Assume the hypotheses on $\phi$ and $(L,Y)$ from Proposition \ref{prop: two ends version of thm} and that $II_\sigma$ holds. Fix a high multiplicity curve $\ell$ from Case $II_\sigma$. 
    Fix $\epsilon > 0$. There exists $M = M(\epsilon)$ such that 
    \begin{align}\label{eq: bound lem: case II}
        |N_\sigma(\ell) \cap E_L| \gtrsim_{\epsilon_0, \epsilon} \delta^\epsilon \sigma (\tau \lambda)^M \delta^{1 + a} \mu. 
    \end{align}
\end{lemma}

The remainder of the present section and the next are devoted to the proof of Lemma \ref{lem: case II}. Let us first see how Proposition \ref{prop: two ends version of thm} follows from Lemmas \ref{lem: case I}, \ref{lem: case II}, and \ref{lem: bush bound}. 

\begin{proof}[Proof of Proposition \ref{prop: two ends version of thm} from Lemmas \ref{lem: case I}, \ref{lem: case II}, and \ref{lem: bush bound}]
    By Lemma \ref{lem: mult dichotomy}, $I$ and $II_\sigma$ hold at a multiplicity $\mu$. Lemma \ref{lem: case I} says that 
    \begin{align}\label{eq: case I}
        |E_L| \gtrsim \lambda \mu^{-1} \delta^2 \#L. 
    \end{align}
    For each $\x_0 \in \R^3$, consider the new shading $(L,Y')$, where $Y'(\ell) = Y(\ell) \setminus B(\x_0, c_0)$. The shading $Y(\ell)$ is $\lambda$-dense and $\epsilon_0$-two-ends, so we may take $c_0 = c_0({\epsilon_0})$ small enough that $Y'(\ell)$ is $\lambda/2$-dense and still $\epsilon_0$-two-ends. Thus Lemma \ref{lem: case II} applies to give
    \begin{align}\label{eq: two-ends case II}
        |N_\sigma(\ell) \cap E_L \cap (\R^3 \setminus B(\x_0,c_0)| \gtrsim_{\epsilon, \epsilon_0} \delta^{\epsilon} \sigma (\tau \lambda)^{M(\epsilon)} \delta^{1 + a} \mu,
    \end{align}
    for each $\ell$ in the set of high multiplicity curves $L'$. 
    Let $C$ be the implicit constant in Lemma \ref{lem: bush bound}
    and refine $L'$ to a set $L''$ of $(C/c_0)\sigma$-direction separated curves, where $\# L'' \gtrsim (\delta/\sigma)^2 \#L' \gtrapprox (\delta / \sigma)^2 \# L$. We now apply Lemma \ref{lem: bush bound} with $\epsilon$ replaced by $\sigma / c_0$, $\sigma$ the same, $E$ replaced by $E_L$, $L$ replaced by $L''$, and $\lambda$ replaced by $\tilde \lambda = \delta^\epsilon \sigma^{-1} (\tau \lambda)^M \delta^{1 + a} \mu$:
    \begin{align}\label{eq: better with mu}
        |E_L| &\gtrsim (\delta^{\epsilon/2} \sigma^{-1} (\tau \lambda)^M \delta^{1+a} \mu \sigma^2 (\# L'')^{1/2} \\
        &\gtrsim \delta^{2 \epsilon} (\tau \lambda)^M \delta^{2+a} \mu (\# L)^{1/2}.
    \end{align}
    Taking the geometric mean of \eqref{eq: case I} and \eqref{eq: better with mu}, 
    \begin{align}
        |E_L| \gtrsim \delta^{\epsilon} (\tau \lambda)^{M} \delta^{\frac{a + 1}{2}} (\delta^2 \# L)^{3/4}.
    \end{align}
\end{proof}

\subsection{Reducing to an ideal hairbrush}

Let $\ell_0 = \ell_{\xi_0, v_0} \in L$ be a curve satisfying \eqref{eq: high mult tube property}. There are $\gtrapprox \lambda / \delta$ many $\delta$-cubes in $Y(\ell_0)$, each of which is incident to $\gtrapprox \mu$ many curves in $L$. Since each such curve makes an angle $\sim \sigma$ with the stem, it can be counted in at most $O(\sigma^{-1})$ of the $\delta$-cubes in $Y(\ell_0)$. Thus there is a set $H \subset L$ of $\gtrapprox \lambda \mu \delta^{-1} \sigma$ many distinct curves passing through $N_\delta(\ell_0)$. In light of Lemma \ref{lem: metric on C(phi)}, we may assume that the curves in $H$ pass directly through $\ell_0$. Write $(H,Y)$ for the collection of hairbrushes with their associated shadings. Since $|\xi(\ell) - \xi(\ell_0)| \sim \sigma$ for each $\ell \in H$, Lemma \ref{lem: curve escape rate} shows that 
\begin{align}
    E_H \subset N_{O(\sigma)}(\ell_0) \cap E_L,
\end{align}
so it is enough to prove \eqref{eq: bound lem: case II} with the left-hand side replaced by $|E_H|$.
We will now shear, cut, and stretch $E_H$ to obtain an ideal hairbrush.

\medskip

\noindent \textbf{Shearing:}

Apply the shear transformation 
\begin{align}
    k(x,t) = (x - \Phi(\xi_0, v_0, t), t)
\end{align}
to $E_H$. We get a new collection $H_1$ of curves $\ell_1 = k(\ell)$ with shadings $Y_1(\ell_1) = k(Y(\ell))$ by sheared $\delta$-cubes. Notice that $Y_1$ is still $\epsilon_0$-two-ends. 

\medskip

\noindent \textbf{Cutting out the center of $E_{H_1}$:}

For $\ell_1 \in H_1$, let $\theta(\ell_1)$ be the $t$-coordinate where $\ell_1$ intersects the stem (the $t$-axis). Since $Y_1(\ell_1)$ is $\epsilon_0$-two-ends, the set 
\begin{align}
    Y_2(\ell_1) = Y_1(\ell_1) \setminus B((0,0,\theta(\ell_1)), c_0)
\end{align}
satisfies $|Y_2(\ell_1)| \geq \frac{\lambda}{2} \delta^2$, as long as $c_0 = c_{\epsilon_0}$ is small enough. Define $H_2 = H_1$. We may now assume that for each $\ell_2 \in H_2$, $|t - \theta(\ell)| \geq c_0$. 

\medskip

\noindent \textbf{Radial rescaling of $E_{H_2}$:}

Define the map $j(x,t) =(\sigma^{-1}x,t)$. Consider the family of rescaled curves $H_3 = \{j(\ell_2) : \ell_2 \in H_2\}$, and define the \emph{measurable} shadings $(H_3, Y_3')$, where $Y_3'(\ell_3):=j(Y_2(\ell_2))$ and $\ell_3 = j(\ell_2)$. We have $|Y_3'(\ell_3)| =\sigma^{-2} |Y_2(\ell_2)| \geq \lambda (\sigma^{-1} \delta)^2$. Apply Lemma \ref{lem: meas to cubes} to $(H_3, Y_3')$ to obtain a $\Omega(\lambda)$-dense shading $(H_3, Y_3)$ by $\sigma^{-1} \delta$-cubes, where 
\begin{align}\label{eq: H3 vol comparison}
    |E_{(H_3, Y_3)}| \sim |E_{(H_3, Y_3')}|
\end{align}
and 
\begin{align}\label{eq: H3 neighborhood comparison}
    E_{(H_3, Y_3)} \subset N_{O(\sigma^{-1} \delta)}(E_{(H_3,Y_3')}).
\end{align}
Lemma \ref{lem: case II} is clearly true if $\lambda \lessapprox \sigma^{-1} \delta$, so we could assume $\lambda \gg \sigma^{-1} \delta$ in order to apply Lemma \ref{lem: meas to cubes}. We may continue to assume $|t - \theta(\ell)| \geq c_0$ as in the previous step because of \eqref{eq: H3 neighborhood comparison}. By \eqref{eq: H3 vol comparison} and the definition of the shadings $Y_3'$, 
\begin{align}
    |E_{(H_3, Y_3)}| \sim \sigma^{-2} |E_{(H_2,Y_2)}|.
\end{align}
The Kakeya estimate required to give Lemma \ref{lem: case II} is thus
\begin{align}
    |E_{H_3}| & 
    \gtrsim \delta^{\epsilon} \sigma^{-1} (\lambda \tau)^{M} \delta^{1 + a} \mu,
\end{align}
for some $M = M(\epsilon)$. 

\medskip

\noindent \textbf{The family containing $H_3$:}

We now write explicitly the family of curves containing $H_3$. The curves in the family are parameterized by where they intersect the stem, $\theta$, and their ``direction'' $\xi$. After applying the shearing step, these are parameterized by
\begin{align}\label{eq: param after sheet}
    \ell_{1, \xi, \theta}(t) = (\Phi_{\xi_0}^\tau(\xi, \theta, t), t), 
\end{align}
where
\begin{align}
    \Phi_{\xi_0}(\xi, \theta, t) &= \nabla_\xi \psi(t, \xi_0 + \xi) - \nabla_\xi\psi(\theta, \xi_0 + \xi) + \nabla_\xi\psi(\theta, \xi_0) - \nabla_\xi\psi(t, \xi_0) \\
    \Phi^\tau_{\xi_0}(\xi, \theta, t) &= \tau^{-1} \Phi_{\xi_0}(\xi, \tau \theta, \tau t).
\end{align}
Here $|\xi| \sim \sigma$, $\sigma \in [\delta, c_1 \tau]$, and $\xi_0 \in \Bphi^2$.
After cutting out the center, the range of $t$ is restricted to $|t - \theta| \geq c_0$, $t \in B(t_0, 1)$. After the radial rescaling steps, the curves are parameterized by 
\begin{align}\label{eq: ugly param curves}
    \ell_{3, \xi, \theta}(t) = (\Phi_{\xi_0}^{\tau,\sigma}(\xi, \theta, t), t),
\end{align}
where 
\begin{align}
    \Phi_{\xi_0}^{\tau, \sigma}(\xi, \theta, t) = \sigma^{-1} \tau^{-1} \Phi_{\xi_0}(\sigma \xi, \tau \theta, \tau t).
\end{align}
The points $\{\xi(\ell_3) : \ell_3 \in H_3\}$ are $\delta \sigma^{-1}$-separated and contained in an annulus of inner and outer radii $\sim 1$. 
Below we summarize the idealized hairbrush families of curves. 
\begin{definition}[Hairbrush family of curves $\mathcal H_{\xi_0}^{\sigma, \tau}(\phi)$]\label{def: hairbrush family}
    Fix $\epsilon_0 > 0$, $\tau \in (0, 1]$, $\sigma \in (0, c_1 \tau]$.
    Let $\xi_0 \in \Bphi^2$. The family $\mathcal H_{\xi_0}^{\sigma, \tau}(\phi)$ consists of curves $\ell_{\xi, \theta}$ parameterized by 
    \begin{align}
        \ell_{\xi, \theta}(t) = (\Phi_{\xi_0}^{\tau, \sigma}(\xi, \theta,t), t), 
    \end{align}
    for $t$ satisfying $t \in B(t_0, 1)$, $|t - \theta| \geq c_0$, and $\theta \in B(t_0, 1)$, $|\xi| \sim 1$. Write the domain 
    \begin{align}
        \Omega = \{(\xi, \theta, t) : |t - \theta| \geq c_0, \ t,\theta \in B(t_0,1), |\xi| \sim 1\}.
    \end{align}
\end{definition}

\noindent \textbf{Finishing the argument assuming conditions of Theorem \ref{thm: CT 3 param kak} are satisfied: }

The maps $\Phi_{\xi_0}^{\tau, \sigma}$ are the same as the defining functions in Definition \ref{def: 3 param family}, except that their domains are slightly more complicated. 
For the purposes of applying Theorem \ref{thm: CT 3 param kak} to estimate $|E_{H_3}|$, we can certainly lose $O_{\epsilon_0, \phi}(1)$ factors. So we may cover the domain of $\Phi_{\xi_0}^{\tau, \sigma}$ by $O(1)$ many balls of radius $\sim 1$, and apply Theorem \ref{thm: CT 3 param kak} in the ball with the most curves. There is therefore no issue in checking that $\Phi_{\xi_0}^{\tau, \sigma}$ satisfies the basic, coniness, and twistiness conditions at each point in $\Omega^{\tau, \sigma}$ directly. We will undertake these calculations in Section \ref{sec: 8} and show that:
\begin{lemma}\label{lem: hairbrushes satisfy conditions}
    Suppose that $\phi$ is translation invariant, non-positive-definite, and satisfies the open condition \eqref{eq: main kakeya open condition} in Theorem \ref{thm: main phicurved kakeya}. Then $\mathcal H_{\xi_0}^{\sigma, \tau}(\phi)$ satisfies (with implicit constants possibly depending on $\epsilon_0$ and $\phi$):
    \begin{itemize}
        \item Basic condition with $C_i \lesssim_i 1$,
        \item $\Omega(\tau^2)$-coniness, and
        \item $\Omega(\tau^4)$-twistiness.
    \end{itemize}
\end{lemma}

\begin{remark}
    Lemma \ref{lem: hairbrushes satisfy conditions} is the only point where \eqref{eq: main kakeya open condition} is used. 
\end{remark}

Once Lemma \ref{lem: hairbrushes satisfy conditions} has been verified, we apply Theorem \ref{thm: CT 3 param kak} with:
\begin{itemize}
    \item $\sigma^{-1} \delta$ in place of $\delta$,
    \item $E_{H_3}$ in place of $E$,
    \item $H_3$ in place of $L$,
    \item $\Omega(\tau^2)$ in place of $\mathfrak c$,
    \item $\Omega(\tau^4)$ in place of $\tau$, and
    \item $\Omega(\lambda)$ in place of $\lambda$
\end{itemize}
to give 
\begin{align}
    |E_{H_3}| &\gtrsim (\delta / \sigma)^\epsilon \tau^{4M + 4} \lambda^M (\delta / \sigma)^a ((\delta /\sigma)^2 \# H_3) \\
    &\gtrsim \delta^{\epsilon} \sigma^{-a} \sigma^{-1} (\lambda \tau)^{M'} \delta^{1 + a} \mu, \\
    &\gtrsim \delta^\epsilon \sigma^{-1}(\lambda \tau)^{M'} \delta^{1+a} \mu
\end{align}
where $M' = 4M + 4$, and we used $\sigma \lesssim 1$. This completes the proof of Lemma \ref{lem: case II} conditional on Lemma \ref{lem: hairbrushes satisfy conditions}.

\section{Proof of Theorem \ref{thm: main phicurved kakeya} Part 2: Verifying the 3-parameter family is twisty and coney (Lemma \ref{lem: hairbrushes satisfy conditions})}\label{sec: 8}

\subsection{An instructive special case}

The proof of Lemma \ref{lem: hairbrushes satisfy conditions} is motivated by a special case calculation for quadratic phase functions $\phi_{A,B}$, with $A \in \{I, I_-\}$ and $B$ symmetric: 
\begin{example}\label{ex: special case for pf of conetwist hairbrush}
    Consider the family $\mathcal H_{0}^{1,1}(\phi_{A,B})$ with defining function 
    \begin{align}\label{eq: special case cone det}
        \Phi^{1,1}_0(\xi, \theta, t) &= ((t-\theta)A + (t^2 - \theta^2) B)\xi. 
    \end{align}
    There is a quadratic form $\xi \mapsto Q(\xi)$ such that 
    \begin{align}\label{eq: special case M}
    (t - \theta) \det(\dot \omega_{\xi,\theta,t}(\theta), \ddot  \omega_{\xi,\theta,t}(\theta)) &= 2Q(\xi) + O(|(t,\theta)|), \\
        \det M_{\xi, \theta}(t) &= 4Q(\xi)^2 + O(|(t,\theta)|).
    \end{align}
    The form $Q$ is definite precisely when $A = I_-$ and 
    \begin{align}
        |B_{12}(B_{22} - B_{11})| > \frac{1}{2} |(B_{22} + B_{11})(B_{22} - B_{11})|, 
    \end{align}
    which is exactly \eqref{eq: main kakeya open condition} when $C = 0$. 
\end{example}
This follows from the Mathematica calculation in Listing \ref{listing: quadratic hairbrush} in Appendix \ref{appendix: B}. Consequently, $\mathcal H_0^{1,1}(\phi_{A,B})$ satisfies both twistiness and coniness everywhere precisely when the hypotheses of Theorem \ref{thm: main phicurved kakeya} are met. 
For general translation-invariant $\phi$, our goal is to show that \eqref{eq: special case cone det} and \eqref{eq: special case M} are true to leading order in $\xi$, though with a quadratic form $Q$ that depends on $A$, $B$, \emph{and} $C$. It is likely not a coincidence that the quadratic form and it's square describe the leading order of the twistiness determinant and the coniness determinant respectively. Though we suspect that this property won't persist if one tries to generalize to the non-translation invariant case. Below we run the calculations for translation invariant $\phi$. We cannot assume $\sigma = \tau = 1$, though we can reduce to this case by various applications of the chain rule. 

\subsection{Useful calculus tools}

It will be much more convenient to work with $\Phi_{\xi_0}$ than $\Phi_{\xi_0}^{\tau, \sigma}$. In what follows, we will do the main calculations with $\Phi_{\xi_0}$, and translate them to $\Phi^{\tau, \sigma}_{\xi_0}$ by the chain rule. Going forward we suppress the subscript $\xi_0$ for ease of calculation. Also let 
\begin{align}
    A &= \partial_t \nabla_\xi^2 \psi(0,0), \\
    B &= \tfrac{1}{2} \partial_t^2 \nabla_\xi^2 \psi(0,0), \\
    C &= \tfrac{1}{6} \partial_t^3 \nabla_\xi^2 \psi(0,0),
\end{align}
and define the skew-symmetric matrix 
\begin{align}
    a = 
    \begin{pmatrix}
        0 & 1 \\
        -1 & 0
    \end{pmatrix}.
\end{align}

The spread curves in Definition \ref{def: spread curves} will appear frequently, and we give the dictionary between spread curves for $\Phi^{\tau, \sigma}$ and $\Phi$ below. 
\begin{lemma}\label{lem: spread curve dictionary}
    Let $\omega_{\xi, \theta, t}(\theta')$ be the spread-curves associated to $\Phi$, and $\omega_{\xi, \theta, t}^{\tau, \sigma}(\theta')$ the spread-curves associated to $\Phi^{\tau, \sigma}$. Then 
    \begin{align}\label{eq: spread curve dictionary}
        \omega_{\xi, \theta, t}^{\tau, \sigma}(\theta') = \sigma^{-1} \omega_{\sigma \xi, \tau \theta, \tau t}(\tau \theta').
    \end{align}
\end{lemma}

\begin{proof}
    Define $\hat \xi_{\xi, \theta, t}(\theta')$ implicitely by 
    \begin{align}\label{eq: hat xi implicit}
        \Phi(\hat \xi_{\xi, \theta, t}(\theta'), \theta', t) = \Phi(\xi, \theta, t)
    \end{align}
    and $\hat \xi_{\xi, \theta, t}^{\tau, \sigma}(\theta')$ implicitely by 
    \begin{align}
        \Phi^{\tau, \sigma}(\hat \xi_{\xi, \theta, t}^{\tau, \sigma}(\theta'), \theta', t) = \Phi^{\tau, \sigma}(\xi, \theta, t). 
    \end{align}
    Thus $\hat \xi_{\xi,\theta, t}^{\tau, \sigma}(\theta') = \sigma^{-1} \hat \xi_{\sigma \xi, \tau \theta, \tau t}(\tau \theta')$, so by the chain rule
    \begin{align}
        \omega_{\xi, \theta, t}^{\tau, \sigma}(\theta') &= \sigma^{-1} \partial_t \Phi(\hat \xi_{\sigma \xi, \tau \theta, \tau t}(\tau \theta'), \tau \theta', \tau t) \\
        &= \sigma^{-1} \omega_{\sigma \xi, \tau \theta, \tau t}(\tau \theta'). 
    \end{align}
\end{proof}

When working with $\Phi$, the following calculus formulas will prove handy. 
\begin{lemma}\label{lem: calc formulae}
    Define
        \begin{align}
            K(\xi,\theta, t) &= \int_0^1 \partial_t \nabla_\xi^2 \psi(\theta, \xi_0 + s \xi) ds \\
            J(\xi,\theta,t) &= \int_0^1 \partial_t \nabla_\xi^2 \psi(\theta + s(t - \theta), \xixi) ds \\
            J'(\xi,\theta,t) &= \int_0^1 \int_0^1 s \partial_t^2 \nabla_\xi^2 \psi(\theta + rs(t-\theta), \xixi) dr ds \\
            J''(\xi,\theta,t) &= \int_0^1 \int_0^1 (1-s) \partial_t^2 \nabla_\xi^2 \psi(\theta + (s + r(1-s))(t-\theta), \xixi) dr ds \\
            J'''(\xi,\theta,t) &= \int_0^1 \int_0^1 \int_0^1 s^2 r \partial_t^3 \nabla_\xi^2 \psi(\theta + qrs(t-\theta),\xixi) dq dr ds. 
        \end{align}
        Then 
        \begin{align}
            \nabla_\xi \Phi &= (t - \theta) J \\
            \partial_\theta \Phi &= -K \xi \\
            J - \partial_t \nabla_\xi^2 \psi(\theta, \xixi) &= (t - \theta)J' \label{eq: J to J' identity}\\
            J - \partial_t \nabla_\xi^2 \psi(t, \xixi) &= -(t-\theta) J'' \label{eq: J to J'' identity}\\
            J' - \frac{1}{2}\partial_t^2 \nabla_\xi^2 \psi(\theta, \xixi) &= (t-\theta) J''' \label{eq: J' to J''' identity}. 
        \end{align}
\end{lemma}

The proof of Lemma \ref{lem: calc formulae} is many easy applications of the fundamental theorem of calculus. One upshot is that $K(\xi, \theta, t)$ and $J(\xi, \theta, t)$ have determinant $\sim 1$ for $(\xi,\theta, t) \in \Bphi^2 \times \Bphi \times \Bphi$, and all of the other terms are $O(1)$. 
Using Lemma \ref{lem: calc formulae}, we can give explicit formulas for $\dot \omega_{\xi, \theta, t}(\theta)$, $\ddot \omega_{\xi, \theta, t}(\theta)$, and $\dddot \omega_{\xi, \theta, t}(\theta)$. 

\begin{lemma}\label{lem: formulas for omega}
    \begin{align}
        \dot \omega_{\xi, \theta, t}(\theta) &= (t - \theta)^{-1} \partial_t \nabla_\xi^2 \psi(t, \xixi) J^{-1} K \xi \\
        &= O(|t - \theta|^{-1} |\xi|) \\
        \ddot \omega_{\xi,\theta, t}(\theta) &= (t - \theta)^{-1} \partial_t \nabla_\xi^2 \psi(t, \xixi) J^{-1}(\partial_\theta K + 2(t - \theta)^{-1} \partial_t \nabla_\xi^2 \psi(\theta, \xixi) J^{-1} K)\xi \\& \qquad+ O(|t - \theta|^{-1} |\xi|^2) \nonumber \\
        &= O(|t - \theta|^{-2} |\xi|)\\
        \dddot \omega_{\xi, \theta,t}(\theta) &= O(|t - \theta|^{-3} |\xi|). 
    \end{align}
\end{lemma}

\begin{proof}
    Start with $\hat \xi_{\xi,\theta, t}(\theta')$ given implicitely by 
    \begin{align}
        \Phi(\hat \xi_{\xi, \theta, t}(\theta'), \theta', t) = \Phi(\xi, \theta, t). 
    \end{align}
    As shorthand, we suppress the subscripts. 

    \medskip
    
    \noindent \textbf{Computing $\dot \omega_{\xi,\theta, t}(\theta)$: }
    
    After taking one derivative in $\theta'$, 
    \begin{align}
        \nabla_\xi \Phi(\hat \xi(\theta'), \theta', t) \partial_{\theta'}\hat \xi(\theta') + \partial_\theta \Phi(\hat \xi(\theta'), \theta', t)= 0. 
    \end{align}
    As further shorthand, we will suppress the arguments. So the above equation becomes 
    \begin{align}\label{eq: implicit one deriv}
        \nabla_\xi \Phi \cdot \partial_{\theta'} \hat \xi + \partial_\theta \Phi = 0. 
    \end{align}
    After evaluating the above equation at $\theta' = \theta$, using Lemma \ref{lem: calc formulae}, and solving for $\partial_{\theta'} \hat \xi$, we find 
    \begin{align}
        \partial_{\theta'} \hat \xi(\theta) = (t - \theta)^{-1} J^{-1} K \xi,
    \end{align}
    where the implicit arguments are $(\xi, \theta, t)$. Also 
    \begin{align}\label{eq: dot xi order}
        \partial_{\theta'} \hat \xi(\theta) = O(|t - \theta|^{-1} |\xi|),
    \end{align}
    which we will make use of in later calculations. 
    We now compute $\dot \omega(\theta)$. By the chain rule and the calculus formulas, 
    \begin{align}
        \dot \omega(\theta) &= \nabla_\xi \partial_t \Phi \cdot \partial_{\theta'} \hat \xi + \partial_\theta \partial_t \Phi \\
        &= \nabla_\xi \partial_t \Phi \cdot \partial_{\theta'} \hat \xi \\
        &= (t - \theta)^{-1} \nabla_\xi^2 \partial_t \psi(t, \xixi) J^{-1} K \xi.  
    \end{align}
    Then $\dot \omega(\theta)$ is of the order 
    \begin{align}\label{eq: ddot xi order}
        \dot \omega(\theta) = O(|t - \theta|^{-1} |\xi|). 
    \end{align}
    \noindent \textbf{Computing $\ddot \omega_{\xi, \theta, t}(\theta)$:}
    
    As one more notation simplification, define $\dot \xi = \partial_{\theta'} \hat\xi$,  $\ddot \xi = \partial_{\theta'}^2 \hat\xi$, and $\dddot \xi = \partial_{\theta'}^3 \hat\xi$.
    Now differentiate \eqref{eq: implicit one deriv} once more to find 
    \begin{align}\label{eq: implicit two deriv}
        \nabla_\xi \Phi \cdot \ddot \xi + \nabla_\xi^2 \Phi [\dot \xi, \dot \xi] + 2 \partial_\theta \nabla_\xi  \Phi \cdot \dot \xi + \partial_\theta^2 \Phi = 0. 
    \end{align}
    Rearranging \eqref{eq: implicit two deriv} and applying Lemma \ref{lem: calc formulae}, we get
    \begin{align}
        \ddot \xi &= -(t - \theta)^{-1} J^{-1} ((t - \theta) \nabla_\xi J[\dot \xi, \dot \xi] - 2 \partial_t \nabla_\xi^2 \psi(\theta, \xixi) \cdot \dot \xi - \partial_\theta K \cdot \xi) \\
        &= (t - \theta)^{-1} J^{-1} (\partial_\theta K + 2(t - \theta)^{-1} \partial_t \nabla_\xi^2 \psi(\theta, \xixi) J^{-1} K)\xi - J^{-1} \nabla_\xi J[\dot \xi, \dot \xi].
    \end{align}
    For future reference, $\ddot \xi$ has order 
    \begin{align}
        \ddot \xi = O(|t - \theta|^{-2} |\xi|).
    \end{align}
    We take two derivatives of $\omega(\theta')$ at $\theta' = \theta$ to find 
    \begin{align}
        \ddot \omega(\theta) &= \partial_t \nabla_\xi \Phi \cdot \ddot \xi + \partial_t \nabla_\xi^2 \Phi[\dot \xi, \dot \xi] \\
        &= \partial_t \nabla_\xi^2 \psi(t, \xixi) \cdot \ddot \xi + O(|t - \theta|^{-2} |\xi|^2) \\
        &= (t - \theta)^{-1} \partial_t \nabla_\xi^2 \psi(t, \xixi) J^{-1}(\partial_\theta K + 2(t - \theta)^{-1} \partial_t \nabla_\xi^2 \psi(\theta, \xixi) J^{-1} K)\xi \\ & \qquad+ \nabla_\xi \partial_t \nabla_\xi^2 \psi(t,\xixi)[\dot \xi, \dot \xi] - \partial_t \nabla_\xi^2 \psi(t,\xixi) J^{-1} \nabla_\xi J[\dot \xi, \dot \xi].\label{eq: err term omega ddot}
    \end{align}
    By \eqref{eq: J to J'' identity},
    \begin{align}
        \partial_t \nabla_\xi^2 \psi(t,\xixi) = J + (t-\theta)J''.
    \end{align}
    Therefore the last two terms in \eqref{eq: err term omega ddot} are equal to 
    \begin{align}
        (t-\theta) (\nabla_\xi J''[\dot \xi, \dot \xi]- J'' J^{-1} \nabla_\xi J[\dot \xi, \dot \xi]) = O(|t-\theta|^{-1} |\xi|^2),
    \end{align}
    where we used \eqref{eq: dot xi order}. Hence 
    \begin{align}
        \ddot \omega(\theta) &= (t - \theta)^{-1} \partial_t \nabla_\xi^2 \psi(t, \xixi) J^{-1}(\partial_\theta K + 2(t - \theta)^{-1} \partial_t \nabla_\xi^2 \psi(\theta, \xixi) J^{-1} K)\xi \\& \qquad+ O(|t - \theta|^{-1} |\xi|^2) 
    \end{align}
    Computing the order of $\ddot \omega$, we get 
    \begin{align}
        \ddot \omega(\theta) = O(|t-\theta|^{-2} |\xi|). 
    \end{align}

    \noindent \textbf{Estimating $\dddot \omega_{\xi, \theta, t}$:}
    
    By differentiating \eqref{eq: implicit two deriv}, we get 
    \begin{align}
        \nabla_\xi \Phi \cdot \dddot \xi + \partial_\theta^3 \Phi + \nabla_\xi^3 \Phi[\dot \xi, \dot \xi, \dot \xi] + 3(\nabla_\xi^2 \Phi[\dot \xi,\ddot \xi] + \partial_\theta \nabla_\xi \Phi \cdot \ddot \xi + \partial_\theta \nabla_\xi^2 \Phi[\dot \xi, \dot \xi] + \partial_\theta^2 \nabla_\xi \Phi \cdot \dot \xi) = 0.
    \end{align}
    By solving for $\dddot \xi$ and estimating term-wise using Lemma \ref{lem: calc formulae} and \eqref{eq: dot xi order}, \eqref{eq: ddot xi order}, we find that 
    \begin{align}
        \dddot \xi = O(|t - \theta|^{-3} |\xi|). 
    \end{align}
    Finally differentiate $\omega$ three times and evaluate at $\theta' = \theta$ to find 
    \begin{align}
        \dddot \omega &= \partial_t \nabla_\xi \Phi \cdot \dddot\xi + \partial_t \nabla_\xi^3 \Phi[\dot \xi, \dot \xi, \dot \xi] + 3 \partial_t \nabla_\xi^2 \Phi[\dot \xi, \ddot \xi] \\
        &= O(|t -\theta|^{-3} |\xi|). 
    \end{align}
\end{proof}

It will also be important to compute the projections under $\pi_\p$ of curves, which are described by the vector $\nabla_{\xi,\theta} \Phi^{\tau,\sigma}(\xi,\theta,t)^T a \gamma_{\xi,\theta}(t)$. Define 
\begin{align}
    \Gamma^{\tau,\sigma}_{\xi,\theta}(t) &= \nabla_{\xi} \Phi^{\tau, \sigma}(\xi,\theta,t)^T a \gamma_{\xi,\theta}^{\tau,\sigma}(t), \\
    \Lambda_{\xi, \theta}^{\tau, \sigma}(t) &= \partial_\theta \Phi^{\tau,\sigma}(\xi,\theta,t)^T a \gamma_{\xi,\theta}^{\tau,\sigma}(t),
\end{align}
and similarly define $\Gamma_{\xi,\theta}(t)$ and $\Lambda_{\xi,\theta}(t)$, but with the superscripts $\tau$ and $\sigma$ removed. 
Again, we have suppressed the subscript $\xi_0$ for easier notation. The following formulas are clear from the chain rule: 
\begin{lemma}\label{lem: Gamma Lambda change of vars formulas}
    We have the change of variables formulas 
    \begin{align}
        \Gamma_{\xi,\theta}^{\tau,\sigma}(t) &= \sigma^{-1} \Gamma_{\sigma\xi,\tau\theta}(\tau t) \\
        \Lambda_{\xi,\theta}^{\tau,\sigma}(t) &= \tau \sigma^{-2} \Lambda_{\sigma \xi, \tau \theta}(\tau t).
    \end{align} 
\end{lemma}
We can use Lemma \ref{lem: calc formulae} to understand $\Gamma_{\xi,\theta}$ and $\Lambda_{\xi,\theta}$:
\begin{lemma}\label{lem: Gamma Lambda formulas}
    \begin{align}
        \Gamma_{\xi,\theta}(t) &= Ja\partial_t\nabla_\xi \phi(t,\xixi) J^{-1}K \xi \\
        &=: \Gamma'_{\xi,\theta}(t) \xi \\
        \Lambda_{\xi,\theta}(t) &= -\xi^T K a J'' J^{-1} K \xi \\
        &=: -\xi^T \Lambda'_{\xi,\theta}(t) \xi
    \end{align}
\end{lemma}
\begin{proof}
    For the first equation, compute
    \begin{align}
        \Gamma_{\xi,\theta}(t) &= \nabla_\xi \Phi(\xi,\theta,t)^T a \gamma_{\xi,\theta}(t) \\
        &= ((t-\theta) J)a (t-\theta)^{-1} \partial_t \nabla_\xi^2 \psi(t,\xixi) J^{-1} K \xi) \\
        &= Ja \partial_t \nabla_\xi^2 \psi(t,\xixi) J^{-1}K\xi.
    \end{align}
    For the second equation, use the identity \eqref{eq: J to J'' identity}: 
    \begin{align}
        \Lambda_{\xi,\theta}(t) &= \partial_\theta \Phi(\xi,\theta,t)^T a \gamma_{\xi,\theta}(t) \\
        &= (-K \xi)^T a (t-\theta)^{-1} \partial_t \nabla_\xi^2 \psi(t,\xixi) J^{-1} K \xi) \\
        &= -(t -\theta)^{-1} \xi^T K a(I + (t-\theta) J''J^{-1}) K\xi\\
        &= -\xi^T Ka J'' J^{-1}K \xi. 
    \end{align}
    Here we used that $\xi^T KaK\xi=0$.
\end{proof}

\subsection{Verifying $\mathcal H_{\xi_0}^{\sigma, \tau}(\phi)$ satisfies the basic conditions}

\ \\
\noindent \textbf{Item (1)}

We may compute using the chain rule
\begin{align}
    \nabla_\xi \Phi^{\tau,\sigma}(\xi, \theta,t) &= \tau^{-1} \nabla_\xi \Phi(\sigma \xi, \tau \theta, \tau t) \\
    &= (t - \theta) J(\sigma \xi, \tau \theta, \tau t). 
\end{align}
We have $|t - \theta| \geq c_0$ and $J$ is non-degenerate, so 
$|\det \nabla_\xi \Phi^{\tau, \sigma}| \gtrsim 1$. Also compute 
\begin{align}
    \partial_\theta \Phi^{\tau, \sigma}(\xi,\theta,t) &= -K\xi.
\end{align}
Since $K=O(1)$ and $|\xi| \sim 1$, $|\partial_\theta \Phi^{\tau, \sigma}(\xi, \theta, t)| \lesssim 1$. For the last item of basic condition Part 1, compute 
\begin{align}
    \partial_{\xi_i, \xi_j} \Phi_k^{\tau, \sigma}(\xi,\theta,t) &=
    \sigma (t - \theta) \partial_{\xi_i} J_{jk}(\sigma \xi, \tau \theta, \tau t), 
\end{align}
which has size $\lesssim \sigma \leq 1$.  

\noindent \textbf{Item (2)}

Use Lemmas \ref{lem: spread curve dictionary} and \ref{lem: formulas for omega} to find 
\begin{align}
    \gamma_{\xi,\theta}^{\tau, \sigma}(t) = \dot \omega^{\tau, \sigma}_{\xi,\theta,t}(\theta) &= \sigma^{-1} \tau \omega_{\sigma \xi, \tau \theta, \tau t}(\tau \theta) \\
    &= \sigma^{-1} \tau (\tau t - \tau \theta)^{-1} \partial_t \nabla_\xi^2 \psi(\tau t, \xi_0 + \sigma \xi) J^{-1} K \cdot (\sigma \xi) \\
    &= (t - \theta)^{-1} \partial_t \nabla_\xi^2 \psi(\tau t, \xi_0 + \sigma \xi) J^{-1} K \xi,
\end{align}
where $K$ and $J$ are evaluated at $(\sigma \xi, \tau \theta, \tau t)$ and are nondegenerate on $\Omega$. Since $|t - \theta| \geq c_0$, we find
\begin{align}
    |\gamma_{\xi,\theta}^{\tau, \sigma}(t)| \sim 1.
\end{align}
This gives the second item. 

\noindent \textbf{Item (3)}

For the third bullet, a similar application of Lemma \ref{lem: spread curve dictionary} and \ref{lem: formulas for omega} shows 
\begin{align}
    |\ddot \omega^{\tau,\sigma}_{\xi,\theta,t}(\theta)|, |\dddot \omega^{\tau,\sigma}_{\xi,\theta,t}(\theta)| \lesssim 1. 
\end{align}
We also need to check that the $\partial_t^i \Phi^{\tau, \sigma}$ are bounded: 
\begin{align}
    \partial_t^i \Phi^{\tau,\sigma} &= \sigma^{-1} \tau^{i - 1} \partial_t^i \Phi(\sigma \xi, \tau \theta, \tau t) \\
    &= \tau^{i - 1} O_i(|\xi|) = O(1)
\end{align}
on $\Omega$, where the last line is clear from the definition of $\Phi$. It is also easy to see that $\Phi^{\tau, \sigma}$ is bounded. This finishes Item (3). 

\noindent \textbf{Item (4)}

By Lemma \ref{lem: Gamma Lambda change of vars formulas}, 
\begin{align}
    \partial_t^i \Gamma_{\xi,\theta}^{\tau, \sigma}(t) &= \tau^i \partial_t^i \Gamma'_{\sigma \xi,\tau \theta}(\tau t) \xi \\
    \partial_t^i \Lambda_{\xi,\theta}^{\tau,\sigma} &= -\tau^{i + 1} \xi^T \partial_t^i\Lambda_{\xi,\theta}'(\tau t) \xi. 
\end{align}
The functions $(\xi_0, \xi,\theta,t) \mapsto \partial_t^i \Gamma'_{\xi,\theta}(t)$, $(\xi_0, \xi, \theta, t) \mapsto \partial_t^i \Lambda'_{\xi,\theta}(t)$ are continuous on $\Bphi^6$, so by compactness they are $O_i(1)$ here. Also $|\xi| \sim 1$, so $\partial_t^i \Gamma_{\xi,\theta}^{\tau,\sigma}(t), \partial_t^i \Lambda_{\xi,\theta}^{\tau,\sigma}(t) = O_i(1)$ on $\Omega$. This completes the verification of the basic conditions. 

\subsection{Verifying $\mathcal H_{\xi_0}^{\sigma, \tau}(\phi)$ is $\Omega(\tau^2)$-coney}

We need to show the following: 
\begin{claim}\label{claim: coniness check}
    For each $(\xi, \theta, t) \in \Omega$, 
\begin{align}\label{eq: coniness check: det lower bd}
    |\det(\dot \omega^{\tau,\sigma}_{\xi,\theta,t}(\theta), \dot \omega^{\tau,\sigma}_{\xi,\theta,t}(\theta))| \gtrsim \tau^2
\end{align}
and 
\begin{align}\label{eq: coniness check: bd dddot omega}
    |\dddot \omega^{\tau, \sigma}_{\xi,\theta,t}(\theta)| \lesssim 1. 
\end{align}
\end{claim}

The bound \eqref{eq: coniness check: bd dddot omega} follows immediately from Lemmas \ref{lem: spread curve dictionary} and \ref{lem: formulas for omega}. The bound \eqref{eq: coniness check: det lower bd} would follow from the claim below.
\begin{claim}\label{claim: coniness check: det lower bd reduced}
    Suppose that $(\xi_0,\xi,\theta,t) \in \Bphi^6$ and $|\xi| \leq c |t - \theta|$ with $\epsilon_\phi$ and $c = c_\phi$ chosen small enough. Then
    \begin{align}
    |\det(\dot \omega_{\xi,\theta,t}(\theta), \ddot \omega_{\xi,\theta,t}(\theta))| \gtrsim |t - \theta|^{-1} |\xi|^2.
\end{align}
\end{claim}

Let us first see why \eqref{eq: coniness check: det lower bd} follows from Claim \ref{claim: coniness check: det lower bd reduced}. 

\begin{proof}[Claim \ref{claim: coniness check: det lower bd reduced} implies \eqref{eq: coniness check: det lower bd}]
    By Lemma \ref{lem: spread curve dictionary}, 
    \begin{align}
        \det(\dot \omega^{\tau,\sigma}_{\xi, \theta, t}(\theta), \ddot \omega^{\tau,\sigma}_{\xi, \theta, t}(\theta)) &= 
        \sigma^{-2} \tau^3 \det(\dot \omega_{\sigma \xi, \tau \theta, \tau t}(\tau \theta), \ddot \omega_{\sigma \xi, \tau \theta, \tau t}(\tau \theta)).
    \end{align}
    Applying Claim \ref{claim: coniness check: det lower bd reduced}, we find 
    \begin{align}
        |\det(\dot \omega^{\tau,\sigma}_{\xi, \theta, t}(\theta), \ddot \omega^{\tau,\sigma}_{\xi, \theta, t}(\theta))| &\gtrsim \tau^2 |t - \theta|^{-1} |\xi|^2 \\
        &\geq c_0 \tau^2 / 4,
    \end{align}
    as long as $(\xi_0, \sigma \xi, \tau \theta, \tau t) \in \Bphi^6$ (which is true since $\sigma, \tau \leq 1$), and $|\sigma \xi| \leq c |\tau (t - \theta)|$. The second inequality is true as long as $\sigma \leq c_0 c/2$, which is guaranteed by setting $c_1 = c_0 c/ 2$ (after setting $c_0$ at the beginning).
\end{proof}

\begin{proof}[Proof of Claim \ref{claim: coniness check: det lower bd reduced}]
Applying the formulas in Lemma \ref{lem: formulas for omega}, 
\begin{align}
    \det(\dot \omega, \ddot \omega) &= 
    (t - \theta)^{-2} \underbrace{ \det(\partial_t \nabla_\xi^2 \psi(t, \xixi) J^{-1})}_{(I)} \label{eq: full det I II expr} \\& \qquad\cdot \underbrace{\det(K\xi, (\partial_\theta K+2(t-\theta)^{-1}\partial_t \nabla_\xi^2 \psi(\theta, \xixi)J^{-1}K)\xi)}_{(II)}. \nonumber \\
    &+ O(|t-\theta|^{-2} |\xi|^3)
\end{align}

Since $\partial_t \nabla_\xi^2 \psi(t, \xixi)$ and $J^{-1}$ are both nondegenerate, 
\begin{align}
    (I) \sim 1. 
\end{align}
We now work on estimating $(II)$. The identity \eqref{eq: J to J' identity} gives 
\begin{align}
    (t - \theta)^{-1}\partial_t \nabla_\xi^2 \psi(\theta, \xixi) J^{-1} K = (t - \theta)^{-1} K - J' J^{-1} K.
\end{align}
Plugging this into $(II)$ and regrouping, we can remove the singularity at $t = \theta$:
\begin{align}
    (II) &= \det(K\xi, (\partial_\theta K - 2J' J^{-1} K)\xi + 2(t - \theta)^{-1} K\xi) \\
    &= \det(K\xi, (\partial_\theta K - 2 J'J^{-1}K)\xi).
\end{align}
We now expand $K$ and $\partial_\theta K - 2 J'J^{-1}K$ in $\xi$ near 0: 
\begin{align}
    K(\xi,\theta,t) &= \partial_t \nabla_\xi^2 \psi(\theta, \xi_0) + O(|\xi|) \label{eq: K expand in xi}\\
    \partial_\theta K - 2 J'J^{-1}K &= \partial_t^2 \nabla_\xi^2 \psi(\theta, \xi_0) - 2J'(0,\theta,t)J^{-1}(0,\theta, t) \partial_t \nabla_\xi^2 \psi(\theta, \xi_0) + O(|\xi|) \label{eq: resolve eqn again},
\end{align}
where the error terms are uniform in $(\xi_0, \theta, t) \in \Bphi^4$. By applying \eqref{eq: J to J' identity} and $\eqref{eq: J' to J''' identity}$, the main term in \eqref{eq: resolve eqn again} can be turned into 
\begin{align}\label{eq: fully resolved}
    2(t - \theta)(J' J^{-1} J' - J''')\mid_{\xi = 0}. 
\end{align}
Plugging \eqref{eq: K expand in xi} and \eqref{eq: fully resolved} back into $(II)$,
\begin{align}
    (II) &= 2(t - \theta) \underbrace{\det(\partial_t \nabla_\xi^2 \psi(\theta, \xi_0) \xi, (J' J^{-1} J' - J''')\mid_{\xi = 0} \cdot\xi)}_{(II)_{\rm{main}}} + O(|\xi|^3).
\end{align}
We have finally resolved $(II)$ enough to expand each of $\partial_t \nabla_\xi^2 \psi(\theta, \xi_0)$ and $(J'J^{-1} J'-J''')\mid_{\xi = 0}$ around $(\xi_0, \theta,t) = 0$. We obtain 
\begin{align}
    (II) = 2(t - \theta) \det(A \xi, (BA^{-1}B - C) \xi) + O(|(\xi_0, \theta, t)| |t - \theta| |\xi|^2 + |\xi|^3).
\end{align}
Finally by combining with \eqref{eq: full det I II expr}, 
\begin{align}
    \det(\dot \omega, \ddot \omega) &= 2(t - \theta)^{-1} \cdot (I) \cdot \underbrace{\det(A \xi, (B A^{-1} B - C)\xi)}_{Q(\xi)} \\ & \qquad+ \underbrace{O(|t-\theta|^{-1} |(\xi_0, \theta,t)| |\xi|^2 + |t-\theta|^{-2} |\xi|^3)}_{\rm{Error}}.
\end{align}

\noindent \textbf{Quadratic form $Q(\xi)$ is definite:}

The leading order of $\det(\dot \omega, \ddot \omega)$ is a quadratic form. 
We have the identity $\det(v,w) = v^T a w$, so 
\begin{align}
    Q(\xi) &= \xi^T \underbrace{A a(B A^{-1} B -C)}_M \xi \\
    &= \xi^T \sym(Aa(B A^{-1} B - C))\xi,
\end{align}
where $\sym(M)$ is the symmetrization of the matrix $A$. The quadratic form $Q$ is definite exactly when $\det \sym(M) > 0$. We calculate $\det \sym(M)$ when $A = I$ and $A = I_-$ separately: 

\noindent \textbf{$A = I$ Case:}

We have $M = a(B^2 - C)$. The matrix $B^2 - C$ is symmetric, say 
\begin{align}
    B^2 - C = 
    \begin{pmatrix}
        x & y \\
        y & z 
    \end{pmatrix}.
\end{align}
Then
\begin{align}
    \det \sym(M) &= - (y^2 + (\frac{z-x}{2})^2) \leq 0,
\end{align}
so $Q$ is never definite when $A = I$, that is when $\phi$ is a positive-definite phase function. 

\noindent \textbf{$A = I_-$ Case: }

In this case, we compute 
\begin{align}
    M = (I_-)a(B I_-^{-1} B - C) &= 
    \begin{pmatrix}
        B_{12}(B_{11} - B_{22}) - C_{12} & B_{12}^2 - B_{22}^2 - C_{22} \\
        B_{11}^2 - B_{12}^2 - C_{11} & B_{12}(B_{11} - B_{22}) - C_{12}
    \end{pmatrix}, \\
    \sym(M) &= 
    \begin{pmatrix}
        B_{12}(B_{11} - B_{22}) - C_{12} & (B_{11}^2 - B_{22}^2 - C_{11} - C_{22})/2 \\
        (B_{11}^2 - B_{22}^2 - C_{11} - C_{22})/2 & B_{12}(B_{11} - B_{22}) - C_{12}
    \end{pmatrix} \\
    \det \sym(M) &= (B_{12}(B_{11} - B_{22}) - C_{12})^2 - ((B_{11}^2 - B_{22}^2 - C_{11} - C_{22})/2)^2.
\end{align}
Thus $Q(\xi)$ is definite exactly when 
\begin{align}
    |B_{12}(B_{22} - B_{11}) + C_{12}| > \frac{1}{2} |B_{22}^2 - B_{11}^2 + C_{11} + C_{22}|,
\end{align}
which is the condition \eqref{eq: main kakeya open condition}. 

\noindent \textbf{Finishing the proof of Claim \ref{claim: coniness check: det lower bd reduced}:}

Since $Q$ is definite and $(I) \sim 1$, there is $c = c_\phi$ so that 
\begin{align}
    |2(t - \theta)^{-1}(I) Q(\xi)| \geq c |t - \theta|^{-1} |\xi|^2.
\end{align}
As long as $|(\xi_0, \theta, t)| \ll_\phi 1$ and $|\xi| \ll_\phi |t - \theta|$ we get $\rm{Error} \leq \frac{1}{2} c |t- \theta|^{-1} |\xi|^2$. Thus 
\begin{align}
    |\det(\dot \omega, \ddot \omega)| \gtrsim |t - \theta|^{-1} |\xi|^2. 
\end{align}
This proves the claim. 
\end{proof}

\subsection{Verifying $\mathcal H_{\xi_0}^{\tau, \sigma}(\phi)$ is $\Omega(\tau^4)$-twisty}

The tangency matrix $M^{\tau,\sigma}_{\xi, \theta}$ from \eqref{eq: tangency matrix} is 
\begin{align}
    M^{\tau, \sigma}_{\xi, \theta}(t) = 
    \begin{pmatrix}
        \Gamma_{\xi, \theta}^{\tau, \sigma} &\Lambda_{\xi,\theta}^{\tau, \sigma} \\
        \dot \Gamma_{\xi, \theta}^{\tau, \sigma} &\dot \Lambda_{\xi,\theta}^{\tau, \sigma} \\
        \ddot \Gamma_{\xi, \theta}^{\tau, \sigma} &\ddot \Lambda_{\xi,\theta}^{\tau, \sigma}.
    \end{pmatrix}
\end{align}
Define $M_{\xi,\theta}(t)$ analogously. By Lemma \ref{lem: Gamma Lambda change of vars formulas},
\begin{align}\label{eq: twist matrix after chain rule}
    M^{\tau, \sigma}_{\xi,\theta}(t) = 
    \begin{pmatrix}
        1 & 0 & 0 \\
        0 & \tau & 0 \\
        0 & 0& \tau^2
    \end{pmatrix}
    M_{\sigma \xi, \tau \theta}(\tau t) 
    \begin{pmatrix}
        \sigma^{-1} & 0 & 0\\
        0 & \sigma^{-1} & 0\\
        0 & 0 & \sigma^{-2} \tau
    \end{pmatrix}.
\end{align}

It is enough to prove that 
\begin{claim}\label{claim: twist calculation}
    $|\det M_{\xi,\theta}(t)| \gtrsim |\xi|^4$, as long as $(\xi_0, \xi, \theta, t) \in \Bphi^6$, for $\epsilon_\phi$ chosen small enough. 
\end{claim}

Indeed, upon taking the determinant of \eqref{eq: twist matrix after chain rule} we would get 
\begin{align}
    |\det M_{\xi,\theta}^{\tau, \sigma}(t)| = \tau^4 \sigma^{-4} |\det(M_{\sigma\xi,\tau\theta}(\tau t))| \gtrsim \tau^4 |\xi|^4 \gtrsim \tau^4
\end{align}
on $\Omega$. We now prove Claim \ref{claim: twist calculation}, and along the way will find the square of the quadratic form $Q(\xi)$ from the coniness condition calculation. 

\begin{proof}[Proof of Claim \ref{claim: twist calculation}]
    We will compute $\Gamma_{\xi,\theta}(t), \Lambda_{\xi,\theta}(t)$, and their first two derivatives in $t$ at the origin. After subtracting terms independent of $t$ or affine in $\xi$, which do not affect $\Phi$, we may use the expansion
    \begin{align}
        \psi(t,\xi) &= \frac{1}{2} \langle (tA + t^2 B + t^3 C + t^4 D) \xi, \xi \rangle + O(|t|^5 |\xi|^2 + |t| |\xi|^3),
    \end{align}
    where $A,B,C,D$ are symmetric $2 \times 2$ matrices, and $A \in \{I, I_-\}$. 

    \medskip
    
    \noindent \textbf{Calculation of $\Gamma_{\xi,\theta}(t)$ and $t$-derivatives:}
    
    Use the formula for $\Gamma'_{\xi,\theta}(t)$ in Lemma \ref{lem: Gamma Lambda formulas} and evaluate at $(\xi_0,\xi,\theta,t) = (0,0,0,t)$:
    \begin{align}
        \Gamma'_0(t) &= [A  + tB + t^2 C]a[A + 2tB + 3t^2 C][A + tB + t^2C]^{-1} A + O(|t|^3). 
    \end{align}
    To simplify this expression, we will use the following identities
    \begin{align}
        M a M &= \det(M)a \\
        AaB + BaA &= c(A,B)a, \text{ $c(A,B) := A_{11}B_{22} + A_{22} B_{11} - 2A_{12} B_{12}$},
    \end{align}
    which hold when $M,A,B$ are symmetric. After expanding the matrix $[A+tB+t^2C]^{-1}$ to order $2$ in $t$ and grouping,
    \begin{align}
        \Gamma'_0(t) &= \det(A)a + tc(A,B)a + t^2[CaA + 2AaC - AaBA^{-1} B + \det(B)a] + O(|t|^3). 
    \end{align}
    We may now read off the derivatives in $t$ at the origin and give expansions to constant order with $(\xi_0,\xi,\theta,t) \in \Bphi^6$: 
    \begin{align}
        \Gamma'_{\xi,\theta,t}(t) &= \det(A) a + O(\epsilon_\phi) \\
        &=: E_1 + O(\epsilon_\phi) \\
        \dot \Gamma'_{\xi,\theta,t}(t) &= c(A,B) a + O(\epsilon_\phi)\\
        &= E_2 + O(\epsilon_\phi)\\
        \ddot \Gamma'_{\xi,\theta,t}(t) &= 2(CaA + 2AaC - AaBA^{-1}B + \det(B)a) + O(\epsilon_\phi) \\
        &=: E_3 + O(\epsilon_\phi). 
    \end{align}

    \noindent \textbf{Calculation of $\Lambda_{\xi,\theta,t}(t)$ and $t$-derivatives:}

    Use the formula for $\Lambda_{\xi,\theta}'(t)$ in Lemma \ref{lem: Gamma Lambda formulas} and evaluate at $(\xi_0,\xi,\theta,t) = (0,0,0,t)$: 
    \begin{align}
        \Lambda_{0}'(t) &= Aa[B + 2tC + 3t^2 D][A + tB + t^2 C]^{-1} A + O(|t|^3). 
    \end{align}
    Once again expand $[A + tB  +t^2 C]^{-1}$ to order 2 in $t$ and group terms: 
    \begin{align}
        \Lambda_{0}'(t) &= AaB + tAa(2C - BA^{-1}B) \\&\qquad + t^2Aa(3D + B(A^{-1} B)^2 - BA^{-1} C - 2CA^{-1} B) + O(|t|^3).\nonumber
    \end{align}
    Again, we may read off the derivatives in $t$ to get 
    \begin{align}
        \Lambda_{\xi,\theta}'(t) &= AaB + O(\epsilon_\phi) \\
        &=: F_1 + O(\epsilon_\phi) \\
        \dot \Lambda_{\xi,\theta}'(t) &= Aa(2C - BA^{-1}B) + O(\epsilon_\phi) \\
        &=: F_2 + O(\epsilon_\phi) \\
        \ddot \Lambda'_{\xi,\theta}(t) &= 2Aa(3D + B(A^{-1} B)^2 - BA^{-1} C - 2CA^{-1} B) + O(\epsilon_\phi)\\
        &=: F_3 + O(\epsilon_\phi).
    \end{align}

    \noindent \textbf{Calculation of $\det M_{\xi,\theta}(t)$:}
        
    Applying the formulas in the above two steps, we have 
    \begin{align}\label{eq: twist phi final expr}
        M_{\xi,\theta}(t) &= 
        \underbrace{
        \begin{pmatrix}
            E_1 \xi & -\xi^T F_1 \xi \\
            E_2 \xi & -\xi^T F_2 \xi \\
            E_3 \xi & -\xi^T F_3 \xi \\
        \end{pmatrix}}_{M(\xi)}
        + O(\epsilon_\phi(|\xi|, |\xi|^2)).
    \end{align}

    \begin{claim}\label{claim: M quadratic form}
        We have $\det M(\xi) = -4 Q(\xi)^2$, where $Q$ is the quadratic form from the coniness step,
        \begin{align}
            Q(\xi) = \xi^T 
            \begin{pmatrix}
                B_{11} B_{12}-B_{12} B_{22}  + C_{12} & (B_{11}^2 - B_{22}^2 - C_{11} - C_{22})/2 \\
                (B_{11}^2 - B_{22}^2 - C_{11} - C_{22})/2 & 
                B_{11} B_{12} - B_{12} B_{22} - C_{12}
            \end{pmatrix}
            \xi
        \end{align}
    \end{claim}
    Once we establish Claim \ref{claim: M quadratic form}, the work in the coniness step shows that $\det M(\xi) \gtrsim |\xi|^4$ for all $\xi$ exactly when Condition \eqref{eq: main kakeya open condition} is met. By taking $\epsilon_\phi$ small enough, the error in \eqref{eq: twist phi final expr} becomes small enough and Claim \ref{claim: twist calculation} is established. 

    \medskip

    \noindent \textbf{Establishing Claim \ref{claim: M quadratic form}: }

    Start by expanding $M(\xi)$ down the third column: 
    \begin{align}
        \det M(\xi) &= -(\xi^T F_1 \xi) \det(E_2 \xi, E_3 \xi) + (\xi^T F_2 \xi) \det(E_1 \xi, E_3 \xi) - (\xi^T F_3 \xi) \det(E_1 \xi, E_2 \xi) \\
        &= -(\xi^T F_1 \xi) \det(E_2 \xi, E_3 \xi) + (\xi^T F_2 \xi) \det(E_1 \xi, E_3 \xi), \label{eq: M quadratic form 1}
    \end{align}
    since 
    \begin{align}
        \det(E_1 \xi, E_2 \xi) = \det(A) c(A,B) \det(a\xi,a\xi) = 0.
    \end{align}
    We can factor \eqref{eq: M quadratic form 1} as follows: 
    \begin{align}
        \det M(\xi) &= \xi^T [\det(A) F_2 - c(A,B) F_1] \xi \cdot \det(a\xi, E_3 \xi) \\
        &= (\xi^T \underbrace{[\det(A) F_2 - c(A,B) F_1]}_{Q_1}\xi) \cdot (\xi^T \underbrace{E_3}_{Q_2} \xi). 
    \end{align}
    Since $A=I_-$, a straightforward calculation shows that
    \begin{align}
        \sym(Q_1) &= -\sym(Q_2) \\ &= 2
        \begin{pmatrix}
                B_{11} B_{12}-B_{12} B_{22}  + C_{12} & (B_{11}^2 - B_{22}^2 - C_{11} - C_{22})/2 \\
                (B_{11}^2 - B_{22}^2 - C_{11} - C_{22})/2 & 
                B_{11} B_{12} - B_{12} B_{22} - C_{12}
            \end{pmatrix}.
    \end{align}
    Therefore $\det M(\xi) = - 4Q(\xi)^2$, proving Claim \ref{claim: M quadratic form} and completing the proof of Claim \ref{claim: twist calculation}.
\end{proof}

\appendix

\section{Proofs of Theorems \ref{thm: C 3 param kak} and \ref{thm: PT 3 param kak}}\label{appendix: A}

In this section we give the proof of Theorem \ref{thm: C 3 param kak} and sketch the proof of Theorem \ref{thm: PT 3 param kak}. We leave them here in the appendix since they are not used to give the main result, but they extend the possible applications of the coniness and twistiness conditions.

\subsection{Proof of Theorem \ref{thm: C 3 param kak}}

Define the shadings $Y(\ell) = N_\delta(\ell) \cap E$. 
By cutting the parameter space into $1 \times 1 \times \mathfrak c / |\log \delta|^3$-slabs, it is enough to show that 
\begin{align}\label{eq: C 3 param kak main}
    |E_L| \gtrsim_\epsilon  \delta^\epsilon \mathfrak c^{1/2} \lambda^2 \delta^{1/2} (\delta^2 \# L), 
\end{align}
where $\theta(L)$ has diameter at most $\mathfrak c / |\log \delta|^3$. We now prove \eqref{eq: C 3 param kak main}

\noindent \textbf{Finding a typical angle $r$ and easy bound}:

We may assume that $Y(\ell) \subset \{\x : \# L(\x) \geq \delta^{-\epsilon}\}$ by cutting $Y(\ell)$ into the low and high multiplicity parts as in Section \ref{sec: 6}. As in Step 2 of Section \ref{sec: 6}, we can pigeonhole a $r = |\log \delta|^{-j}$ and a set of curves with their shadings $(L_1, Y_1)$ with the following properties.
\begin{itemize}
    \item $L_1 \subset L$ and $\#L_1 \gtrapprox \# L$, 
    \item $|Y_1(\ell)| \gtrapprox \lambda \delta^2$, and 
    \item 
    \begin{align}
        Y_1(\ell) \subset \{\x \in Y(\ell) : r_{L}(\x) \sim r \text{ and } \theta(\ell) \in J_{L}(\x)\}. \label{eq: C kak shading cond}
    \end{align}
\end{itemize}
The set of curves through a point is one dimensional, so $\# L_1(\x) \lesssim r/\delta$. This gives the easy bound
\begin{align}\label{eq: C kak easy bd}
    |E_L| \geq
    |E_{L_1}| &\gtrsim (\delta / r) \sum_{\ell \in L_1} |Y_1(\ell)| \\
    &\gtrapprox \lambda (\delta / r) (\delta^2 \# L). 
\end{align}

\noindent \textbf{The curves are contained in disjoint $1 \times r$-tubes and applying multilinear Kakeya inside}:

Let $\mathcal B$ be a set of $O(r^{-3})$ boundedly overlapping balls whose union covers $L_1$. For each $B \in \mathcal B$ let $L_1^B = L_1 \cap B$. Because of \eqref{eq: C kak shading cond}, there are $O(1)$ balls $B$ with $\x \in E_{L_1^B}$. Thus 
\begin{align}
    |\bigcup_{\ell \in L_1} Y_1(\ell)| \gtrsim \sum_{B \in \mathcal B} |\bigcup_{\ell\in L_1^B} Y_1(\ell)|.
\end{align}
By the same application of the multilinear Kakeya theorem as in Step 8 of Section \ref{sec: 6}, we find 
\begin{align}\label{eq: C kak mlk in ball}
    |\bigcup_{\ell \in L_1^B} Y_1(\ell)| \gtrapprox_\epsilon \delta^\epsilon \mathfrak c \lambda^3 r (\delta^2 \# L_1^B). 
\end{align}
Summing over $B \in \mathcal B$, 
\begin{align}\label{eq: C kak mlk bd}
    |E_L| \gtrapprox_\epsilon \delta^\epsilon \mathfrak c \lambda^3 r (\delta^2 \# L). 
\end{align}
Upon taking the geometric mean of \eqref{eq: C kak easy bd} and \eqref{eq: C kak mlk bd}, we get the claim \eqref{eq: C 3 param kak main}. 

\subsection{Sketch of proof of Theorem \ref{thm: PT 3 param kak}}

Steps 1 through 7 in Section \ref{sec: 6} would yield Theorem \ref{thm: PT 3 param kak} if one could take $\tilde \rho = (\rho r)^{1/2}$, instead of $\min((\rho r)^{1/2}, \rho / r)$. That would precisely be the proof contained in \cite{KWZ}. The only obstruction is Lemma \ref{lem: taylor expand C(X,p)}, which forces the prisms to be too narrow. By assuming that $X$ is planey, we have a better expansion of $\ell_{s,\theta'}$ as expected.

\begin{lemma}[Planey version of Lemma \ref{lem: taylor expand C(X,p)}]\label{lem: planey taylor expand C(X,p)}
    Fix $s_0 \in B_1^1$.
    Suppose that $X$ satisfies the basic conditions and is planey. Then the curve $\ell_{s,\theta'} \in \mathcal C(X,\p)$ has the expansion 
    \begin{align}
        \ell_{s,\theta'}(t) - \ell_\p(t) &= ((t-s)(\theta' - \theta + \eta_{\p,s}(\theta')) \gamma_\p(s_0) \\ & \qquad+ O(|\theta' - \theta| |t - s|^2 + |\theta' - \theta| |t -s| |s -s_0|),0),\nonumber
    \end{align}
    where $\eta_{\p,s}(\theta') = O(|\theta' - \theta|^2)$ is a scalar function.
\end{lemma}

\begin{proof}
    With $Y$ as in the proof of Lemma \ref{lem: taylor expand C(X,p)}, Taylor expansion in $t$ gives
    \begin{align}
        Y(\theta',s,t) &= (t-s)(\omega_{\p,s}(\theta') - \omega_{\p,s}(\theta)) + O(|\theta'-\theta||t-s|^2).
    \end{align}
    Planiness gives $\det(\dot \omega_{\p,s}(\theta'), \ddot \omega_{\p,s}(\theta')) = 0$ for every $\theta'$. Consequently, 
    \begin{align}
        \omega_{\p,s}(\theta') - \omega_{\p,s}(\theta) = (\theta'-\theta + \eta_{\p,s}(\theta')) \gamma_{\p}(s),
    \end{align}
    where $\eta_{\p,s}(\theta') = O(|\theta'-\theta|^2)$. Expanding $\gamma_\p(s) = \gamma_\p(s_0) + O(|s-s_0|)$ proves the lemma. 
\end{proof}

We now describe the better grain structure as a consequence of Lemma \ref{lem: planey taylor expand C(X,p)}. Fix $r \geq \delta$, $\rho \in [\delta, (\delta r)^{1/2}]$, and assume $|\theta' - \theta| \sim r$. Then 
\begin{align}
    |\ell_{s,\theta'}(t) - \ell_\p(t)| \sim r|t - s|
\end{align}
as long as $|\theta' - \theta|, |s -s_0|,|t-s| \ll 1$. So $\ell_{s,\theta'}$ is inside $N_\rho(\ell_\p)$ for only $t$ in the range 
\begin{align}
    |t - s| \lesssim \rho / r.
\end{align}
Fix $s_0$ and consider $|s - s_0| \leq (\delta / r)^{1/2}$. For this range of $s$ and inside $N_\rho(\ell)$, we have 
\begin{align}\label{eq: planey 1}
    \ell_{s,\theta}(t) - \ell_p(t) = ((t - s) (\theta' - \theta + \eta_{\p,s}(\theta'))\gamma_\p(s_0) + O(\delta) ,0).
\end{align}
Thus the image of the right-hand side of \eqref{eq: planey 1} inside $N_\rho(\ell)$ and for $|t - s| \leq \rho / r$ is contained in a $\delta \times \rho \times (\delta / r)^{1/2}$-slab. The analogues of Lemmas \ref{lem: prism filled up} and \ref{lem: add prisms fat tube}, and Proposition \ref{prop: narrow prop} continue to hold with $\rho \in [\delta, (\delta r)^{1/2}]$, so we may run Steps 1 through 7 without any $\rho^a$ loss. 

\section{Mathematica calculations}\label{appendix: B}

The below Mathematica function in Listing \ref{listing:coneTwistExpressions} computes the determinant in the coniness/planiness condition, and the matrix $M_{\p,t}$ in the twistiness, 1-flat, and 2-flat conditions. First we take some implicit derivatives of $\hat p_{\p,t}(\theta)$ from Definition \ref{def: spread curves} by hand,
\begin{align}
    \dot p &:= \partial_\theta \hat p_{\p,t}(\theta) = -(\nabla_p X)^{-1} \partial_\theta X 
    \\
    \ddot p &:= \partial_\theta^2 \hat p_{\p,t}(\theta) = -(\nabla_p X)^{-1} \cdot (\nabla_p^2 X[\dot p, \dot p] + 2\partial_\theta \nabla_p X \cdot \dot p + \partial_\theta^2 X).
\end{align}
The relevant derivatives of $\omega_{\p,t}$ can then be computed: 
\begin{align}
    \gamma_{\p}(t) =\dot \omega &:= \dot \omega_{\p,t}(\theta) = \nabla_p \partial_t X \cdot \dot p + \partial_\theta \partial_t X,\\
    \ddot \omega &:= \ddot \omega_{\p,t}(\theta) = \nabla_p \partial_t X \cdot \ddot p + \nabla_p^2 \partial_t X[\dot p, \dot p] + 2 \nabla_p \partial_\theta \partial_t X \cdot \dot p + \partial_\theta^2 \partial_t X.
\end{align}
These formulas feature in Listing \ref{listing:coneTwistExpressions}. Listing \ref{listing: quadratic hairbrush} contains the calculation of the coniness determinant and tangency matrix for hairbrushes from the quadratic phase functions $\phi_{A,B}$. 

\begin{lstlisting}[language=Mathematica, caption={Coniness and Twistiness quantity calculations\label{listing:coneTwistExpressions}}, captionpos=t,basicstyle=\scriptsize\ttfamily]
coneTwistExpressions[XFnc_, pVars_List, thVar_, tVar_] :=
  Module[
   {
    (* Intermediate symbols *)
    thX, pX, ppX, pthX, ththX, ptX, pptX, thtX, pthtX, ththtX, dotP,
    ddotP, dotOmega, ddotOmega, coneDet, proj, tangMatrix
    },
   X = XFnc @@ Join[pVars, {thVar, tVar}];
   (* Required derivatives *)
   thX = D[X, thVar];
   pX = D[X, {pVars, 1}];
   ppX = D[X, {pVars, 2}];
   pthX = D[pX, thVar];
   ththX = D[thX, thVar];
   ptX = D[pX, tVar];
   pptX = D[ppX, tVar];
   thtX = D[thX, tVar];
   pthtX = D[pthX, tVar];
   ththtX = D[ththX, tVar];
   
   (* Compute dotP, ddotP*)
   dotP = -Inverse[pX] . thX;
   ddotP = -Inverse[
       pX] . (ppX . dotP . dotP + 2*pthX . dotP + ththX);
   (* Compute dotOmega, ddotOmega *)
   dotOmega = ptX . dotP + thtX;
   ddotOmega = 
    ptX . ddotP + pptX . dotP . dotP + 2*pthtX . dotP + ththtX;
   (*Coniness*)
   coneDet = Simplify[Det[{dotOmega, ddotOmega}]];
   (*Twistiness*)
   a = {{0, 1}, {-1, 0}};
   proj = 
    Join[Transpose[pX] . a . 
      dotOmega, {Transpose[thX] . a . dotOmega}];
   tangMatrix = 
    Simplify[{proj, D[proj, {tVar, 1}], D[proj, {tVar, 2}]}];
   (* Return *)
   {
    coneDet,
    tangMatrix,
    }
   ];
\end{lstlisting}




\begin{lstlisting}[language=Mathematica, caption={Calculation for Example \ref{ex: special case for pf of conetwist hairbrush} \label{listing: quadratic hairbrush}}, captionpos=t,basicstyle=\scriptsize\ttfamily]
(* Input *)
XFnc[xi1_,xi2_,th_,t_] := ((t-th)*{{1,0},{0,-1}} + (t^2-th^2)*{{B11,B12},{B12,B22}}).{xi1, xi2};

{coneDet, tangMatrix} = coneTwistExpressions[XFnc, {xi1,xi2}, th, t];
(Factor[Simplify[coneDet]]*(t-th)) /. {t -> 0, th -> 0}
Factor[Simplify[Det[tangMatrix]]] /. {t -> 0, th -> 0}

(* Output *)
-2 (-B11 + B22) (B12 xi1^2 + B11 xi1 xi2 + B22 xi1 xi2 + B12 xi2^2)
-4 (-B11 + B22)^2 (B12 xi1^2 + B11 xi1 xi2 + B22 xi1 xi2 + B12 xi2^2)^2
\end{lstlisting}

\bibliographystyle{alpha}
\bibliography{reference}

\end{document}